\documentclass[12pt]{article}
\usepackage[all]{xy}
\usepackage[dvips]{graphicx}
\usepackage{enumerate,hyperref}
\usepackage{ifthen,color}
\usepackage{amsmath,amsfonts,amssymb,amsthm,dsfont}
\usepackage{mathrsfs}

\newtheorem{thm}{Theorem}[section]

\newtheorem{coro}[thm]{Corollary}

\newtheorem{exa}[thm]{Example}

\newtheorem{de}[thm]{Definition}

\newcommand{\comp}{ \,{\scriptstyle \stackrel{\circ}{}}\, }

\newcommand{\Acal}{\mathcal{A}}
\newcommand{\calA}{\mathcal{A}}
\newcommand{\Bcal}{\mathcal{B}}

\newcommand{\Hcal}{\mathcal{H}}
\newcommand{\calH}{\mathcal{H}}
\newcommand{\Ncal}{\mathcal{N}}
\newcommand{\calN}{\mathcal{N}}

\newcommand{\Tcal}{\mathcal{T}}
 
\newcommand{\RN}{ \mathbb{R} } 
 
\newcommand{\CN}{ \mathbb{C} } 
 
\newcommand{\HN}{\mathbb{H}}
\newcommand{\NN}{\mathbb{N}}

\newcommand{\eps}{\epsilon}

\newcommand{\MatA}{\text{M}_{\Acal}}
\newcommand{\Amat}{\text{M}_{\Acal}}
\newcommand{\EndA}{\mathcal{R}_{\Acal}}
\newcommand{\EndB}{\mathcal{L}_{\Bcal}}

\newcommand{\Fun}{\text{C}_{\Acal}}
\newcommand{\MM}{\mathbf{M}}
\newcommand{\Acalx}{\mathcal{A}^{\times}}
\newcommand{\Acalzd}{\mathbf{zd}(\mathcal{A})}

\newcommand{\Bcalx}{\mathcal{B}^{\times}}
\newcommand{\Bcalzd}{\mathbf{zd}(\mathcal{B})}
\newcommand{\ds}{\displaystyle}

\newcommand{\Abound}{m_{\Acal}}

\begin{document}
\begin{titlepage}

\title{{\sc Introduction to $\Acal$-Calculus}}
\author{ James S. Cook \\ jcook4@liberty.edu \\ Liberty University \\ Department of Mathematics}
\date{\today}
\maketitle

\begin{abstract}
Let $\Acal$ denote an $n$-dimensional associative algebra over $\RN$. 
This paper gives an introductory exposition of calculus over $\Acal$. An $\Acal$-differentiable function $f: \Acal \rightarrow \Acal$ is one for which the differential is right-$\Acal$-linear. The basis-dependent correspondence between right-$\Acal$-linear maps and the regular representation of real matrices is discussed in detail. The requirement that the Jacobian matrix of a function fall in the regular representation of $\Acal$ gives $n^2-n$ generalized $\Acal$-CR equations. In contrast, some authors use a deleted-difference quotient to describe differentiability over an algebra. We compare these concepts of differentiability over an algebra and prove they are equivalent in the semisimple commutative case. We also show how difference quotients are ill-equipt to study calculus over a nilpotent algebra.  \\

Cauchy Riemann equations are elegantly captured by the Wirtinger calculus as $\frac{\partial f}{\partial \bar{z}}=0$. We generalize this to any commutative unital associative algebra. Instead of one conjugate, we need $n-1$ conjugate variables.
Our construction modifies that given by Alvarez-Parrilla, Fr{\' i}as-Armenta, L{\' o}pez-Gonz{\' a}lez and Yee-Romero in \cite{pagr2012}. \\

We derive Taylor's Theorem over an algebra. Following Wagner, we show how Generalized Laplace equations are naturally seen from the multiplication table of an algebra. We show how d'Alembert's solution to the wave-equation can be derived from the function theory of an appropriate algebra. The inverse problem of $\Acal$-calculus is introduced and we show how the Tableau for an $\Acal$-differentiable function has a rather special form. The integral over an algebra is also studied. We find the usual elementary topological results about closed and exact forms generalize nicely to our integral. Generalizations of the First and Second Fundamental Theorems of Calculus are given. \\

This paper should be accessible to undergraduates with a firm foundation in linear algebra and calculus. Some background in complex analysis would also be helpful. 
 \end{abstract}


\end{titlepage}
\section{Introduction and overview}
We use $\Acal$ to denote a real associative algebra of finite dimension. Elements of $\Acal$ are known as $\Acal$-numbers. Our program of study is to describe a calculus where real numbers have been replaced by $\Acal$-numbers. The resulting calculus we refer to as $\Acal$-calculus. The goal of this paper is to explain the basic differential and integral features of $\Acal$-calculus. We hope this paper serves as a primer for those who wish to investigate the many open questions in $\Acal$-calculus.  \\

\noindent
In Section \ref{sec:history} we give an overview of the literature which is connected to $\Acal$-calculus. The basic ideas of this research have been known since around 1890. However, the development is largely disconnected. As Waterhouse puts it on page 353-354 of \cite{waterhouseII}:
\begin{quote}
The theory has been developed only in fits and starts, and some of the simple results seem not to be on record at all, while others are only recorded in a language that is now hard to understand.
\end{quote}
Indeed, one major goal of this work is to provide some updates on proofs which were given without the benefit of modern linear algebra. \\

\noindent
Section \ref{sec:calculusonR} introduces the calculus of real normed linear spaces. You can consult \cite{edwards}, \cite{Dmaster} or \cite{zorich} for further details. \\

\noindent
An introduction to associative algebras is given in Section \ref{sec:algebraintro}. We explain how the algebra $\Acal$ is naturally interchanged with linear transformations on $\Acal$ which are right-$\Acal$-linear. In addition, given a choice of basis $\beta$ for $\Acal$ we explain how the regular representation $\MatA( \beta)$ gives another object which is isomorphic to $\Acal$. Isomorphism here requires we preserve both the structure of addition and multiplication. We use the isomorphism of $\Acal$ and $\MatA( \beta)$ to find a number of interesting results about zero-divisors and units in the algebra. In particular, we argue that an invertible basis exists in any unital associative algebra. We also note the group of units $\Acalx$ form a dense subset of $\Acal$. \\

\noindent
In Section \ref{sec:adifffunctions} we discuss how the norm on $\Acal$ is generally {\it submultiplicative}. So far as we know, the estimate derived in Equation \ref{eqn:submultiplicative} has not appeared elsewhere in the generality we offer here. The $\Acal$-differentiability of a function is defined via the representation theory developed in Section \ref{sec:algebraintro}. We prove the derivative over $\Acal$ enjoys many of the same elementary properties as the usual derivative in $\RN$. For example, $\frac{d}{d\zeta} \zeta^n = n \zeta^{n-1}$ for any $n \in \NN$. The sum, difference, scalar produce and suitably qualified composite of $\Acal$-differentiable functions are once more $\Acal$-differentiable. In the case $\Acal$ is commutative the $\Acal$-differentiability of $f$ and $g$ imply $f \star g$ is likewise $\Acal$-differentiable. However, in the noncommutative case it is possible the product of $\Acal$-differentiable functions is not once more $\Acal$-differentiable. Example \ref{exa:uppertriangular} provides functions $f$ and $g$ for which $f \star g$ is $\Acal$-differentiable and yet $g \star f$ is not $\Acal$-differentiable. Finally, we conclude the section by studying differentiable functions over isomorphic algebras are connected. We later use that result to derive d'Alembert's solution to the wave equation.

\noindent
The generalized Cauchy-Riemann Equations are studied in detail in Section \ref{sec:CReqns}. There are several ways to understand the $\Acal$-CR Equations: if $\Acal$ is unital and $\mathds{1}$ is the first basis element of $\beta = \{ v_1, \dots , v_n\}$ then suppose $U$ is open and $f: U \subseteq \Acal \rightarrow \Acal$ we have $f$ is $\Acal$-differentiable on $U$ if $f$ is $\RN$-differentiable on $U$ and for each point in $U$: 
\begin{quote}
\begin{enumerate}[{\bf (i.)}]
\item the Jacobian matrix fits in the regular representation; 
$
\ds J_f \in \MatA (\beta)
$
\item the generalized $\Acal$-equations hold; $\ds \frac{\partial f}{\partial x_j} = \frac{\partial f}{\partial x_1} \star v_j$ for each $j=1, \dots, n$.
\item $f = f( \zeta)$ or $f$ is $\Acal$-holomorphic; $ \ds \frac{\partial f}{\partial \overline{\zeta}_j}=0$ for $j=2, \dots, n$.
\end{enumerate}
\end{quote}
We spend some effort to developing (iii.). The catalyst for our construction was given by Alvarez-Parrilla, Fr{\' i}as-Armenta, L{\' o}pez-Gonz{\' a}lez and Yee-Romero in \cite{pagr2012}. However, we have to differ with \cite{pagr2012} in their construction of $\frac{\partial}{\partial \zeta}$.  Our Theorem \ref{thm:wirtingerworks} is clear evidence that our construction of the Wirtinger calculus over $\Acal$ should be prefered to that in \cite{pagr2012}. Example \ref{exa:trihyperbolicidentities} illustrates how $\Acal$-differentiable functions provide natural solutions to the analog of Laplace's equation for $\Acal$ all in the language of the Wirtinger calculus. Finally, we give a roadmap towards converting a problem of real calculus to one of $\Acal$-calculus. The problem of knowing which $\Acal$ will {\it work} is a difficult and open problem. \\

\noindent
Section \ref{sec:deleteddiffquotients} provides a bridge to the other main definition used to define differentiability over an algebra. In particular, many authors use a limiting process on difference quotients on the algebra. However, the limit is done modulo the zero divisors of the algebra so we might call it a {\bf deleted-difference-quotient} (D1). In contrast, we define differentiability on the algebra through the differential with the healthy backdrop of representation theory in Section \ref{sec:adifffunctions}(D2). Naturally, we wonder if there is a distinction between these approaches.  At a point they are inequivalent. On an open set, our D2 approach is more general than the deleted-difference quotient. Example \ref{exa:nowhereD1} gives an example which is everywhere D2 yet nowhere D1 differentiable. It is no accident that Example \ref{exa:nowhereD1} concerns a nilpotent algebra. We prove in Theorem \ref{thm:D2sameasD1forsemisimple} that D1 and D2 are equivalent characterizations of differentiability on an open set for a commutative semisimple algebra of finite dimension. As far as we know, the results we present here are original. That said, we must give credit to \cite{gadeaD1vsD2} and \cite{price} whose work provided much inspiration for this Section. \\

\noindent
Section \ref{sec:higherderivatives} contains a story which seems too good to be true. In particular, we find yet another isomorphism with $\Acal$ which allows us to circumvent all the usual symmetric-multilinear algebra which is required for higher derivatives of real maps on $\RN^n$. In particular, this leads to Theorem \ref{thm:partialsproducts} where we learn:
\begin{equation} \label{eqn:cool}
 \frac{\partial^k f}{\partial x_{i_1}\partial x_{i_2} \cdots \partial x_{i_k}} = \frac{\partial^k f}{\partial x_1^k} \star v_{i_{1}} \star v_{i_{2}} \star \cdots \star v_{i_{k}}.
 \end{equation}
 Again, in the context of $\Acal$-differentiable functions we find this curious identity which allows us to trade partial derivatives of various coordinates for the one coordinate which is paired with the unity in $\Acal$. It is then a simple matter to reproduce the wonderful result of Wagner from 1948 which says an equation in the algebra corresponds naturally to a generalized Laplace equation. We hope our proof is an improvement to that which exists in the current literature. To illustrate the beauty of Wagner's result we show how the wave equation appears as the generalized Laplace equation of an appropriate algebra. We then derive d'Alambert's solution through an explicit isomorphism to the direct product algebra. Finally, we find a relatively effortless derivation of Taylor's Theorem for $\Acal$-calculus which follows almost as a Corollary to Equation \ref{eqn:cool}. \\
 
 \noindent
The eventual goal of Section \ref{sec:inverseproblem} is recognizing when and how $\Acal$-calculus can be used to gain deeper insight into existing problems of real calculus. It is mostly an invitation to think on the problem for future researchers. Example \ref{exa:AODEbackwards} is given as an inverse to the inverse problem. We can easily create more such examples. It is not difficult to solve a problem in $\Acal$-calculus then convert it to a corresponding problem of real calculus. The problem we would like to gain further insight in future work is the inverse problem; given a problem of real calculus, can we convert it to an more lucid problem of $\Acal$-calculus? We hope Theorem \ref{thm:gaussAmap} takes us a step closer to solving the inverse problem. Perhaps the reader will be amused that once more Equation \ref{eqn:cool} is the core of the proof for Theorem \ref{thm:gaussAmap}. \\

\noindent
Integration over an commutative unital algebra of finite dimension over $\RN$ is studied in Section\ref{sec:integration}. Our integral is a natural generalization of that which is studied in the usual complex analysis. Essentially the same integral can be found in the 1928 Thesis of Ketchum. 
If $C$ is a curve of length $L$ and $f$ is an $\Acal$-differentiable function bounded by $M>0$ on $C$ then Theorem \ref{thm:subML} states:
\begin{equation}
 \bigg{|}\bigg{|} \int_C f( \zeta) \star d\zeta \bigg{|}\bigg{|} \leq \Abound ML. 
\end{equation}
Here $\Abound$ is the constant\footnote{derived in Theorem \ref{thm:submultiplicative}} for which $|| v \star w|| \leq \Abound ||v|| \, ||w||$ for all $v,w \in \Acal$. The Fundamental Theorem of Calculus part II is generalized in Theorem \ref{thm:FTCforalg}. Theorem \ref{thm:topologicalAintegral} provides the equivalence of path-independence, trivial closed integrals and the existence of an antiderivative in a connected subset of $\Acal$.  Theorem \ref{thm:exactAdiff} shows $\Acal$-differentiability of $f$ implies exactness of $f \star d\zeta$. We obtain Cauchy's Integral Theorem for $\Acal$ as Corollary \ref{coro:adiffloopszero}.  An analog to The Fundamental Theorem of Calculus part I is given in Theorem \ref{thm:FTCIforalg}. Finally, we indicate a logarithm in $\Acal$ can be define via the integral in Example \ref{exa:logdefined}. \\

\noindent
In our final Section we indicate some directions for future work as well as the unpublished work in preparation by D. Freese and N. BeDell.

\section{History} \label{sec:history}
\noindent
Any history of generalized calculus is necessarily incomplete. Here we mention primarily those authors whose work seems to precede our own, or, descends directly from our closest mathematical relatives. \\

\noindent
Scheffers, a student of Lie, is usually credited with initiating the program of hypercomplex analysis in \cite{scheffers}. Then, Segre \cite{Segre}, Hausdorff \cite{haus}, Spaminato \cite{Spam}, Ringleb \cite{ringleb} and Ketchum \cite{ketchum} joined in Scheffer's program of analysis. See \cite{ward1940} for additional references of early authors in hypercomplex analysis. In \cite{ketchum} we find many of the usual theorems of complex function theory set forth for a general commutative alegbra over $\CN$.  Ketchum's later work on polygenic functions in \cite{ketchum2}  and \cite{ketchum3} deal with an infinite dimensional algebra which provides a function theory for the three dimensional wave equation. Several decades pass until Synder's 1958 thesis \cite{Synder} continues Ketchum's investigations and then Kunz' \cite{kunz} gives further insight. \\

\noindent
It seems most of the work before Ward \cite{ward1940} focused on calculus over and algebra with a complex base field. Ward, a student of MacDuffee at the University of Wisconsin, used matrix arguments  and some of the insight brought from Jacobson's enveloping algebra to define {\it analytic} functions over a noncommutative algebra in his 1939 thesis \cite{wardthesis}. However, some of Ward's theorems are given for functions whose derivative falls inside the regular representation. Wagner, also a student of MacDuffee, was successful in refining Ward's work in the specific context of commutative algebras. Wagner found matrix arguments which showed how to produce {\it generalized Laplace equations} \cite{wagner1948}. After serving in World War II, Ward produced \cite{ward1952} in 1952 which echoed the improvements of Wagner. Ward showed in that given a particular system of PDEs, there exists an algebra for which the system of PDEs form the generalized Cauchy Riemann equations of the constructed algebra. Ward and Wagner's works are a joy. Unfortunately, they are given in a {\it matric} language which probably obscures their logic for most modern readers. Most of their work involves explicit manipulation of structure constants and component equations. To understand their work one has to understand that $C_{ijk}$ is either the components of the first regular representation $R_i$, or second regular representation $S_j$ or the paraisotropic matrix $Q_k$. In some sense, much of what we present in this paper is logically implicit in \cite{ward1940} and \cite{wagner1948}. The works of Ward and Wagner are impressive in that they were able to see as much as they did without the benefit of modern linear algebra. We hope that we can exposit Ward and Wagner's results in a language which exposes how simple and natural they truly are. \\

\noindent
Our work is probably more inline with that of Wagner to be honest. Ward's analytic functions for noncommutative algebras have been studied by Trampus \cite{Trampus} and Rineheart \cite{Rinehart}. \\

\noindent
However, as a point of history, the author only found Ward and Wagner's works after much of the theory of $\Acal$-calculus was independently uncovered through conversations with Nguyen, Leslie and Zhang in the Fall 2012 Semester at Liberty University. In fact,   Vladimirov and Volovich's work in \cite{VV1984} served as the catalyst for this whole project. \\

\noindent
It may be useful to demarcate two major directions in generalized calculus over an algebra:
\begin{quote}
\begin{enumerate}[{\bf (1.)}]
\item function theory of generalized complex numbers
\item function theory of generalized hyperbolic numbers
\end{enumerate}
\end{quote}
\noindent
For {\bf (1.)}, we have the continuation of Scheffer's work on algebras with a complex base field. These take inspiration from Segre's 1892 introduction of the bicomplex numbers. More recently, Price's text {\it Multicomplex Functions and Spaces} discusses the integral and differential calculus of bicomplex and multicomplex numbers \cite{price}.  M.E. Luna-Elizarrar{\'a}s, M. Shapiro, D.C. Struppa and A. Vajiac provide an update to Price's work on the bicomplex numbers in \cite{luna}.
Bicomplex numbers are also related to the {\it complicated numbers} studied by Good and clarified by Waterhouse \cite{waterhouseI}. Good noticed primes of the form $8n+1$ have a particular connection with complicated numbers. In \cite{waterhouseII}Waterhouse gives a reasonably lucid derivation of Wagner's generalized Laplace equations and he gives some evidence that $\Acal$ is well-suited to bring new insights into the analysis of a PDE. Pedersen \cite{pedersenI} \cite{pedersenII} and Jonasson \cite{jonasson} applied Waterhouse's calculus to study polynomial solutions to a large class of PDEs.   Particular algebras have been studied in great depth, for example, see \cite{Dicky} for the state of the art on calculus with Cayley numbers. Plaksa and Pukhtaievych show in \cite{plaksa} the Cauchy integral formula, Morera's theorem and Cauchy's integral formula generalize to certain algebras over $\CN$. I should once more emphasize the paper by Ketchum in 1928 is very related to all the works above as Ketchum explains the basic differential and integral calculus over a complex algebra in \cite{ketchum} \\.

\noindent
For {\bf (2.)}, it seems hyperbolic numbers were introduced as the {\it real tessarines} by James Cockle in 1848 \cite{cockle}. However, over the years these are discovered or invented by numerous disconnected authors. This reality is reflected in the abundance of names which have been given to hyperbolic numbers: 
$$
\begin{array}{|c|c|c|} 
\text{hyperbolic numbers} & 
\text{(real) tessarines,} &
\text{split-complex numbers}  \\ \hline
\text{algebraic motors} & 
\text{bireal numbers,} &
\text{approximate numbers}  \\ \hline
\text{countercomplex numbers} & 
\text{double numbers} &
\text{anormal-complex numbers}  \\ \hline
\text{perplex numbers} & 
\text{Lorentz numbers} &
\text{paracomplex numbers}  \\ \hline
\text{semi-complex numbers} & 
\text{split binarions} &
\text{spacetime numbers}  \\ \hline
\text{Study numbers} & 
\text{twocomplex numbers} &
  \\ \hline
\end{array}
$$
Fox's 1949 thesis \cite{fox} on the bireal numbers gives some interesting pictures to illustrate the analog of conformal mapping in the case of hyperbolic numbers. In 1998, Motter and Rosa study hyperbolic integral analysis and provide some discussion of the analog of a Riemann surface for hyperbolic calculus \cite{MotterRosa}. In 1999, Konderak used Lorentz numbers in \cite{konderak} to prove a result about immersed Lorentz surfaces in $\mathbb{R}^3_1$.  The 2003 careful analysis by Gadea, Grifone and Masqu{\'e} in \cite{gadeaD1vsD2} on double numbers is very helpful towards understanding the distinction between differentiability defined with or without the concept of a difference quotient, although, the paper is largely concerned with the construction of manifolds over the double numbers. In 2005, Khrennikov and Segre give a broad introduction to hyperbolic analysis which includes an analysis of how hyperbolic analysis fits into Clifford analysis \cite{khrennikov}. In 2008, Kravchenko et. al. show in \cite{kravchenkoI} and \cite{kravchenkoII} how hyperbolic numbers to analyze certain solutions to the Klein-Gordon equation. In 2014, Terlizzi,  Konderak and Lacirasella provide many explicit formulae for functions over Lorentz numbers which they use to study manifolds modelled over Lorentz numbers \cite{Terlizzi}. \\

\noindent
In contrast, this paper is written in the tradition of Ward and Wagner which is not specific to a particular algebra. Our general approach also allows some nilpotent examples not widely considered elsewhere. Gadea and Masqu{\'e}'s work in \cite{gadeaAanalytic} is also quite general, they give some picture of how to build manifolds which locally support an $\Acal$-calculus.  Rosenfeld studied differentiability for many noncommutative algebras in \cite{Rosenfeld}. Rosenfeld showed there were very few differentiable functions for a large class of simple algebras, in contrast, we show there are algebras with nilpotent elements which support many differentiable functions\footnote{the author is thankful to Robert Bryant for his comment at \href{http://mathoverflow.net/q/191088/24854}{mathoverflow question 191088} and private correspondence which helped clarify this point.}. Lorch's work \cite{lorch} is also interesting. Lorch's concept of differentiation over an algebra was foundational for the recent papers \cite{pagr2012} and \cite{pagr2015} where it is shown how to solve certain ordinary differential equations via an algebra substitution. In particular, the concept of multiple conjugates as shown in \cite{pagr2012} is very helpful and it was the inspiration for the generalized Wirtinger calculus we present in this paper. \\

\noindent
We apologize to the authors we have inadvertently slighted. Please see \cite{hypercomplexbook2014} to see some other directions which are currently under investigation. In particular, quaternionic analysis continues to be a source of interesting questions.

\section{Calculus on a normed linear space} \label{sec:calculusonR}
Let $V$ and $W$ be finite dimensional normed linear spaces over $\RN$. We denote the norm of $x \in V$ by $||x||$. If $F: U \subset V \rightarrow W$ is a function then we say $F$ is differentiable at $p$ if there exists a $\RN$-linear map $d_pF: V \rightarrow W$ for which
\begin{equation} \label{eqn:diffwrtR}
 \lim_{h \rightarrow 0}\frac{F(p+h)-F(p)-d_pF(h)}{||h||} = 0. 
\end{equation}
Likewise, if $F$ is differentiable for each $p \in U$ then we say $F$ is {\bf differentiable on $U$}. If $\beta  = \{ v_1, \dots , v_n \}$ is a basis for $V$ with coordinate functions $x_1, \dots , x_n$ then we define
\begin{equation}
\frac{\partial F}{\partial x_i}(p) = \lim_{ t \rightarrow 0} \frac{F(p+tv_i)-F(p)}{t}. 
\end{equation}
If the map $p \rightarrow \frac{\partial F}{\partial x_i}(p)$ is continuous at $p$ for all $i=1, \dots , n$ then we say $F$ is {\bf continuously differentiable at $p$} and write $F \in C^1(p)$. Let $U$ be an open set. If $F$ is continuously differentiable on $U$ then $F$ differentiable on $U$ and we write $F \in C^1(U)$. If $F \in C^1(U)$ then we are free to construct the differential of $F$ at $p$ from the partial derivatives of $F$ at $p$; for $h = h_1v_1 + \cdots +h_nv_n$,
\begin{equation} \label{eqn:partialderNLS}
d_pF(h) = \sum_{i=1}^n h_i \frac{\partial F}{\partial x_i}(p) \qquad \text{or, as is often useful,} \qquad \frac{\partial F}{\partial x_i}(p) = d_pF(v_i) . 
\end{equation}
The differential $d_pF: V \rightarrow W$ is a linear transformation of vector spaces hence we can associate, given a choice of bases $\beta$ for $V$ and $\gamma$ for $W$, a matrix $[d_pF]_{\beta, \gamma} \in \RN^{m \times n}$. In particular, the {\bf Jacobian matrix} for $F \in C^1(p)$ is given by:
\begin{equation}
[d_pF]_{\beta, \gamma} = [[d_pF(v_1)]_{\gamma}| \cdots | [d_pF(v_n)]_{\gamma}] = \left[\left[\frac{\partial F}{\partial x_1}\right]_{\gamma} \bigg{|} \cdots \bigg{|} \left[\frac{\partial F}{\partial x_n}\right]_{\gamma} \right] 
\end{equation}
where $[y_1w_1+ \cdots + y_mw_m]_{\gamma} = (y_1, \dots , y_m) \in \RN^m$ is the usual $\gamma$-coordinate map. In many of the applications we study the vector spaces $V$ and $W$ are $\RN^n$ and it is our custom to use $e_1, \dots, e_n$ to denote the {\bf standard basis} where $(e_i)_j = \delta_{ij}$. In this special case, we note the Jacobian matrix simply by the standard matrix of the differential;
\begin{equation}
 [d_pF] = [d_pF(e_1)| \cdots | d_pF(e_n)] = \left[\frac{\partial F}{\partial x^1}\bigg{|} \cdots \bigg{|} \frac{\partial F}{\partial x^n} \right]. 
 \end{equation}

\section{Real linear associative algebras} \label{sec:algebraintro}
A vector space paired with a multiplication forms an {\it algebra}.

\begin{de} \label{defn:algebra}
Let $\Acal$ be a finite-dimensional real vector space paired with a function $\star: \Acal \times \Acal \rightarrow \Acal$ which is called {\bf multiplication}. In particular, the multiplication map satisfies the properties below: 
\begin{enumerate}[{\bf (i.)}]
\item {\bf bilinear:} $ (cx + y ) \star z = c( x \star z)+ y \star z$ and $x \star (cy+ z) = c(x \star y)+ x \star z $ for all $x,y,z \in \Acal$ and $c \in \RN$,
\item {\bf associative:} for which $x \star (y  \star z) = (x \star y)  \star z$  for all $x,y,z \in \Acal$ and,
\item {\bf unital:} there exists $\mathds{1} \in \Acal$ for which $\mathds{1} \star x = x$ and $x \star \mathds{1} = x$.
\end{enumerate}
We say $x \in \Acal$ is an {\bf $\Acal$-number}. If $x \star y = y \star x$ for all $x,y \in \Acal$ then $\Acal$ is {\bf commutative}.
\end{de}

\noindent
When there is no ambiguity we use $1 = \mathds{1}$ and we replace $\star$ with juxtaposition; $xy = x \star y$. We assume $\Acal$ is an associative algebra of finite dimension over $\RN$ throughout the remainder of this paper. In the commutative case there is no need to distinguish between {\it left} and {\it right} properties. However, we allow the possibility that $\Acal$ be {\bf noncommutative} at this point in our development. \\

\noindent
If $\alpha \in \Acal$ then $\ell_{\alpha}(x) = \alpha \star x$ is a {\bf left-multiplication} map on $\Acal$. It is a {\bf right-$\Acal$-linear} as:
\begin{equation}  \label{eqn:rightlinearityofleftmult}
\ell_{\alpha}(x \star y) = \alpha \star (x \star y) = (\alpha \star x ) \star y = \ell_{\alpha}(x) \star y.
\end{equation}
Likewise, $r_{\alpha}(x) =  x \star \alpha$ is a {\bf right-multiplication} map on $\Acal$. It is a {\bf left-$\Acal$-linear} as:
\begin{equation} 
r_{\alpha}(x \star y) = (x \star y) \star \alpha = 
x \star (y \star \alpha ) = x \star r_{\alpha}(y). 
\end{equation}
Notice, associativity of $\star$ is given by $\ell_{\alpha} \comp r_{\beta}  = r_{\beta} \comp \ell_{\alpha}$ for all $\alpha, \beta \in \Acal$. From Equation \ref{eqn:rightlinearityofleftmult} we see that every left-multiplication map is right-$\Acal$-linear. In fact, if $\mathds{1} \in \Acal$ then every right-$\Acal$-linear map on $\Acal$ is a left multiplication by a particular element of $\Acal$.

\begin{thm} \label{thm:allleftlinearmapsaremultiplyonleft}
If $T: \Acal \rightarrow \Acal$ is a right-$\Acal$-linear map then there exists a unique $\alpha \in \Acal$ for which $T = \ell_{\alpha}$.
\end{thm}

\noindent
{\bf Proof:} let $T \in \EndA$ and consider $T(x) = T(\mathds{1} \star x) = T(\mathds{1}) \star x$ for each $x \in \Acal$. Therefore, $T = \ell_{T(\mathds{1})}$. $\Box$ \\

\noindent
The simple calculation above is key to understanding generalized Cauchy Riemann equations.

\begin{de} \label{defn:leftlineartrans}
Let $\EndA$ define the set of all right-$\Acal$-linear transformations on $\Acal$. If $T \in \EndA$ then $T: \Acal \rightarrow \Acal$ is a $\RN$-linear transformation for which
$T(x \star y) = T(x) \star y$ for all $x, y \in \Acal$
\end{de}

\noindent
Recall the sum, scalar multiple and composition of endomorphisms is once more an endomorphism. Moreover, addition, scalar multiplication and composition of transformations are known to be bilinear and associative. Furthermore, if $Id(x) = x$ for all $x \in V$ then observe $Id = \mathds{1}$ for the operation of composition. In summary,  $\text{gl}(V) = \{ T: V \rightarrow V \ | \ T \ \text{linear transformation} \}$ forms the {\bf general linear algebra on $V$}. In the context of $V = \Acal$, the subset $\EndA \subseteq \text{gl}(\Acal)$ is special: 

\begin{thm} \label{thm:subalgebraofleftlinearmaps}
The set $\EndA$ is a subalgebra of $\text{gl}(\Acal)$; that is, $\EndA \leq \text{gl}( \Acal)$\footnote{we intend the notation $\mathcal{A} \leq \mathcal{B}$ to indicate $\mathcal{A}$ is a {\bf subalgebra} of $\mathcal{B}$.}.
\end{thm}

\noindent
{\bf Proof:} notice that $Id(x \star y) = x\star y = Id(x) \star y$ hence $Id \in \EndA$. To show $\EndA$ is a subalgebra of $\text{gl}(\Acal)$ it remains to show the sum, scalar multiple and composite of right-$\Acal$-linear maps in once again in $\EndA$. Let $T_1, T_2 \in \EndA$ and $c \in \RN$. Suppose $x,y \in \Acal$ and consider:
\begin{align} 
 (cT_1+T_2)(x \star y) &= c(T_1)(x \star y)+ T_2(x \star y) \\ \notag
 &= cT_1(x )\star y+ T_2(x) \star y \\ \notag
  &= \left(cT_1(x ) + T_2(x) \right)\star y \\ \notag
  &= (cT_1 + T_2)(x) \star y. \notag
\end{align}
Likewise, applying right-linearity of $T_2$ then $T_1$ yields:
\begin{equation} 
 (T_1\comp T_2)(x \star y) = T_1(T_2(x \star y)) = T_1(T_2(x) \star y)) = T_1(T_2(x)) \star y = (T_1 \comp T_2)(x) \star y. 
\end{equation}
Hence $T_1 \comp T_2 \in \EndA$. $\Box$ \\

\noindent
The result above shows that for any algebra $\Acal$ we immediately obtain a related algebra of linear transformations. Given a choice of basis, we also can trade $\Acal$ for a particular set of matrices known as 
the  {\bf regular representation}. 

\begin{de} \label{defn:regrepofA}
Let $\Acal$ have basis $\beta$ then the {\bf regular representation} with respect to $\beta$ is
$$ \MatA(\beta) = \{ [T]_{\beta, \beta} \ | \ T \in \EndA \}. $$ 
In the case $\Acal = \RN^n$ we may forego the $\beta$ notation and write
$$ \MatA = \{ [T] \ | \ T \in \EndA \}$$
for the {\bf regular representation} of $\Acal$.
\end{de}

\noindent
The regular representation of $\text{gl}(V)$ is simply $\RN^{n \times n}$ matrices given $\text{dim}(V)=n$. This is immediate from the fact that any matrix may appear as the matrix of an endomorphism. Just as $\EndA \leq \text{gl}(\Acal)$ we likewise find $\MatA(\beta) \leq \RN^{n \times n}$. In $\RN^{n \times n}$ the algebra multiplication is simply matrix multiplication and $\mathds{1}$ is the $n\times n$ identity matrix.

\begin{thm}
For any choice of $\beta$, the set $\MatA(\beta)$ is a subalgebra of $\RN^{n \times n}$.
\end{thm}

\noindent
{\bf Proof:} since $[Id]_{\beta, \beta} = I$ it follows that $\mathds{1} \in \MatA(\beta)$. It remains to show that $\MatA(\beta)$ is closed under addition, scalar multiplication and matrix multiplication. Assume $A,B \in \MatA( \beta)$ and $c \in \RN$. Observe, by definition, there exist $S,T \in \EndA$ for which $A = [S]_{\beta, \beta}$ and $B = [T]_{\beta, \beta}$. Since the matrix of a linear combination of operators is the linear combination of the matrices of said operators we have:
\begin{equation}
 cA+B =  c[S]_{\beta, \beta}+ [T]_{\beta, \beta} = [cS+T]_{\beta, \beta}
\end{equation}
but, we know $S,T \in \EndA$ hence by Theorem \ref{thm:subalgebraofleftlinearmaps} $cS+T \in \EndA$ which shows $[cS+T]_{\beta, \beta} \in \MatA$ and thus $cA+B \in \MatA$. Finally, recall matrix multiplication was defined precisely so the identity below holds true:
\begin{equation}
 AB=[S]_{\beta, \beta}[T]_{\beta, \beta} = [S \comp T]_{\beta, \beta}. 
 \end{equation}
Observe,  Theorem \ref{thm:subalgebraofleftlinearmaps} indicates that $S,T \in \EndA$ implies $S \comp T \in \EndA$. Therefore, $AB \in \MatA$ and we conclude $\MatA$ forms a subalgebra of $\RN^{n \times n}$. $\Box$ \\

\noindent
Different basis choices give different {\bf matrix} regular representations. 
We use the change of basis theorem from linear algebra to relate the representations:

\begin{thm}
If $\beta,\gamma$ are bases for $\Acal$ then $\MatA(\beta)$ and $\MatA(\gamma)$ are conjugate subalgebras.
\end{thm}

\noindent
{\bf Proof:} Let $\beta, \gamma$ be basis for $\Acal$. If $T: \Acal \rightarrow \Acal$ is a linear transformation then we know from linear algebra that there exists an invertible change of basis matrix $P$ for which $[T]_{\beta, \beta} = P^{-1}[T]_{\gamma, \gamma}P$. Thus, if $A \in \MatA(\beta)$ then $A=[T]_{\beta, \beta} = P^{-1}[T]_{\gamma, \gamma}P \in P^{-1}\MatA(\gamma)P$. In other words, $\MatA(\beta)$ is the image of $\MatA(\gamma)$ under conjugation by $P$. $\Box$ \\

\noindent
For a given algebra $\Acal$ and basis $\beta$ we are free to use transformations in $\EndA$ or matrices in $\MatA(\beta)$ to capture the same structure. These isomorphisms are of fundamental importance to the study of $\Acal$-calculus. Notice, we say $(\Acal, \star)$ and $(\Bcal, \ast)$ are {\bf isomorphic} as real associative algebras and write $\Acal \approx \Bcal$ if there exists an invertible $\RN$-linear transformation $\Psi: \Acal \rightarrow \Bcal$ such that $\Psi (x \star y) = \Psi(x) \ast \Psi(y)$ for all $x,y \in \Acal$.

\begin{thm} \label{thm:isomorphismtrifecta}
If $\beta$ is a basis for $\Acal$ then $\Acal  \approx \EndA \approx \MatA(\beta)$ 
\end{thm}

\noindent
{\bf Proof:} suppose $\beta = \{ v_1, \dots , v_n \}$ is a basis for $\Acal$. Define $\Psi( \alpha) = \ell_{\alpha}$ as was studied in Equation \ref{eqn:rightlinearityofleftmult}. Notice $\Psi(c\alpha+\beta) = \ell_{c\alpha+\beta} = c\ell_{\alpha}+\ell_{\beta} = c \Psi(\alpha)+\Psi(\beta)$ hence $\Psi$ is a linear transformation. Further, as $\Acal$ is unital there exists a multiplicative identity $\mathds{1} \in \Acal$. If $T \in \EndA$ then $T(\mathds{1}) \in \Acal$ and for $x \in \Acal$ we calculate:
\begin{equation} \label{eqn:theformula}
 (\Psi (T(\mathds{1})))(x) =  \ell_{T(\mathds{1})}(x) = T(\mathds{1}) \star x = T(x).
 \end{equation}
Thus $\Psi(T(\mathds{1})) = T$ which shows $\Psi$ is a surjection. Suppose $\Psi(\alpha) =0$ hence $\ell_{\alpha}(x) = 0$ for all $x \in \Acal$. Set $x = \mathds{1}$ to calculate $\ell_{\alpha}(\mathds{1}) = \alpha \star \mathds{1} = \alpha = 0$. Therefore $\text{Ker}(\Psi) = \{ 0 \}$ and we find $\Psi$ is an injection. Hence $\Psi$ is an isomorphism of vector spaces. It remains to show $\Psi$ preserves the algebra multiplication. We need associativity here:
\begin{equation} \label{eqn:assocleftmaps}
 \Psi ( \alpha \star \beta )(x) = \ell_{\alpha \star \beta}(x) = 
(\alpha \star \beta) \star  x = 
\alpha \star (\beta \star  x) = \ell_{\alpha}( \ell_{\beta}(x)) = \ell_{\alpha} \comp \ell_{\beta}(x). 
 \end{equation}
Therefore, $\Psi$ gives the isomorphism $\Acal \approx \EndA$. Next, define $\MM: \Acal \rightarrow \MatA( \beta)$ by
\begin{equation}
 \MM( \alpha ) = [ \ell_{\alpha} ]_{\beta, \beta} = [ [\alpha \star v_1]_{\beta}| \cdots | [\alpha \star v_n]_{\beta}]. 
 \end{equation}
If $A \in \MatA(\beta)$ then there exists $T \in \EndA$ for which $[T]_{\beta, \beta} = A$. Use Equation \ref{eqn:theformula} to see:
 \begin{equation}
 \MM( T(\mathds{1}) ) = [\ell_{T(\mathds{1})}]_{\beta,\beta} = [T]_{\beta,\beta} = A.
\end{equation}
Hence $\MM: \Acal \rightarrow \MatA(\beta)$ is surjective. Suppose $\MM (\alpha) =0$ then $[\ell_{\alpha}]_{\beta, \beta}=0$ from which we find $\ell_{\alpha}=0$. Thus, $\ell_{\alpha}(\mathds{1}) = \alpha \star \mathds{1} = \alpha = 0$. We find $\text{Ker}(\MM) = \{0 \}$ and thus $\MM$ is injective. Finally, we verify $\MM$ preserves the algebra multiplication: we use the middle of Equation \ref{eqn:assocleftmaps} in the second equality:
\begin{equation}
 \MM( x \star y) = [ \ell_{x \star y} ]_{\beta, \beta} = [ \ell_{x} \comp \ell_{y} ]_{\beta, \beta} = [ \ell_{x} ]_{\beta, \beta}[ \ell_{y} ]_{\beta, \beta} = \MM(x) \MM(y). 
 \end{equation}
Thus $\MM$ provides the isomorphism $\Acal \approx \MatA(\beta)$. $\Box$ \\

\noindent
It is convenient to set some notation for the inverse of the $\Psi$ map. For each number $x \in \Acal$ the map $\Psi$ assigns a linear transformation $\Psi(x) \in \EndA$. It seems natural\footnote{this is {\bf not} a hash-tag} to call the inverse of $\Psi$ the {\bf number map}. For convenience of notation, we also use $\#$ for the inverse of $\MM$:

\begin{de} \label{defn:isomorphismOfEndAandA}
The {\bf number map} $\#: \EndA \rightarrow \Acal$ is defined by $\#(T) = T(\mathds{1})$ when the context demands. Likewise, given $\beta$ a basis for $\Acal$ we define $\#: \MatA(\beta) \rightarrow \Acal$ and when the context demands $\#( [T]_{\beta, \beta}) = T(\mathds{1})$.
\end{de}

\noindent
The number map is easiest to understand when $\mathds{1} = v_1 = e_1$ where $\Acal = \RN^n$, however, we have taken care to allow for other possibilities\footnote{Note, $\RN \times \RN$ we have $\mathds{1} = (1,1)$ and in $\RN^{ 2 \times 2}$ we have $\mathds{1} = \left[ \begin{array}{cc} 1 & 0 \\ 0 & 1 \end{array} \right]$ thus the usual standard bases for $\RN^2$ or $\RN^{2 \times 2}$ do not include the identity of the algebra.}. 

\begin{thm} \label{thm:unitybasismatrixugly}
If $\beta = \{ v_1 , \dots , v_n \}$ is a basis for $\Acal$ where $v_1 = \mathds{1}$ then 
$$\MM( x) = \left[ [x]_{\beta}| [x \star v_2]_{\beta}| \cdots | 
[x \star v_n]_{\beta} \right]$$
and $\#(A) = \Phi_{\beta}^{-1}(\mathbf{col}_1(A))$ where $\Phi_{\beta}(x) = [x]_{\beta}$ is the coordinate map.
\end{thm}

\noindent
{\bf Proof:} we define $\MM: \Acal \rightarrow \MatA(\beta)$ as in the proof of Theorem \ref{thm:isomorphismtrifecta}; $ \MM(x) = [ \ell_{x} ]_{\beta, \beta}$. The $j$-th column in $[ \ell_{x} ]_{\beta, \beta}$ is $[ \ell_{x}(v_j) ]_{\beta} = [x \star v_j]_{\beta}$. Therefore, as $v_1 = \mathds{1}$ we find $\MM( x) = \left[ [x]_{\beta}| [x \star v_2]_{\beta}| \cdots | 
[x \star v_n]_{\beta} \right]$. Let $A = \MM(x)$, we wish to solve for $x$. Notice, $\mathbf{col}_1(A)  = [x]_{\beta}$ from which we derive $x = \Phi_{\beta}^{-1}(\mathbf{col}_1(A))$ hence $\#(A) = \Phi_{\beta}^{-1}(\mathbf{col}_1(A))$.  $\Box$ \\

\noindent
In almost all applications of the Theorem \ref{thm:unitybasismatrixugly} we consider the case $\Acal = \RN^n$ with $\beta = \{ e_1, \dots , e_n \}$ the usual standard basis such that $e_1 = \mathds{1}$. Given these special choices we obtain much improved formulae\footnote{ when $\beta = \{ e_1, \dots, e_n \}$ we drop $\beta$ from the notation and simply write $\MatA$ in the place of $\MatA(\beta)$}

\begin{coro}
If $\Acal = \RN^n$ and $\mathds{1} = e_1 = (1,0,\dots , 0)$ then for $x \in \Acal$ and $A \in \MatA$,
$$ \MM(x) = [x| x \star e_2| \cdots | x \star e_n] \qquad \& \qquad \#(A) = \mathbf{col}_1(A), \qquad \# \MM(x) = x. $$
\end{coro}

\noindent
Notice that the first column of $A \in \MatA$ determines the rest through the structure of the multiplication of $\Acal$. Setting aside the special context of the corollary, if $\beta = \{ v_1, \dots , v_n \}$ is a non-standard basis with $v_j = \mathds{1}$ then the $j$-th column of $\MM(x)$ will fix the remaining columns according to the multiplication on $\Acal$. On the other hand, if $\mathds{1} \notin \beta$ then there need not be a single column of each matrix in $\MatA(\beta)$ which fixes the remaining columns. We see this phenomenon explicitly in Examples \ref{Ex:number5} and \ref{Ex:number10}. 

\begin{de} \label{defn:groupofunits}
We say $x \in \Acal$ is a {\bf unit} if there exists $y \in \Acal$ for which $x \star y = y \star x = \mathds{1}$. The set of all units is known as the {\bf group of units} and we denote this by $\Acalx$. We say $a \in \Acal$ is a {\bf zero-divisor} if $a \neq 0$ and there exists $b \neq 0$ for which $a \star b = 0$ or $b \star a =0$. Let $\Acalzd = \{ x \in \Acal \ | \ x=0 \ \text{or $x$ is a zero-divisor} \}$
\end{de}

\noindent
Isomorphisms transfer both units and zero-divisors.

\begin{thm} \label{thm:isomorphismpreservesstructure}
Suppose $(\Acal, \star)$ and $(\Bcal, \ast)$ are real associative algebras and $\Phi: \Acal  \rightarrow \Bcal$ is an isomorphism. Then,
\begin{enumerate}[{\bf (i.)}]
\item $\Phi( \mathds{1}_{\Acal}) = \mathds{1}_{\Bcal}$,
\item for $x \in \Acalx$, $\Phi(x^{-1}) = \Phi(x)^{-1}$;  that is, $\Phi(\Acalx) = \Bcalx$,
\item if $x,y \in \Acal$ and $x \star y = 0$ then $\Phi(x) \ast \Phi(y) = 0$; that is, $\Phi(\Acalzd) = \Bcalzd$.
\end{enumerate}
\end{thm}

\noindent
{\bf Proof:} to prove (i.) simply note for each $y \in \Bcal$ there exists $x\in \Acal$ for which $y=\Phi(x) = \Phi(x \star \mathds{1}_{\Acal}) = \Phi(x ) \ast \Phi( \mathds{1}_{\Acal}) = y \ast \Phi( \mathds{1}_{\Acal})$. Likewise, as $x = \mathds{1}_{\Acal} \star x$ we have $y = \Phi( \mathds{1}_{\Acal}) \ast y$. Thus, $\mathds{1}_{\Bcal} = \Phi( \mathds{1}_{\Acal})$. As $\Phi(0)=0$, the proofs of (ii.) and (iii.) are immediate from the definition of inverse and zero-divisor since $\Phi(x \star y) = \Phi(x) \ast \Phi(y)$. $\Box$ \\

\noindent
Since $\Acal \approx \MatA(\beta) \approx \EndA$ we are free to focus our initial effort where is most convenient. When considering the characterization of $\Acalzd$ the representation $\MatA(\beta)$ is useful due to the theory of determinants.

\begin{thm} \label{thm:MatAclassicalring}
Let $\MatA(\beta)$ be the regular representation of $\Acal$ with respect to basis $\beta$. If $A \in \MatA(\beta)$ then either $A$ is zero, a unit, or a zero-divisor. 
\end{thm}

\noindent
{\bf Proof:} if $A \in \MatA(\beta)$ then by definition there exists $S \in \EndA$ for which $[S]_{\beta,\beta} = A$. Observe that either $\text{det}(A) = 0$ or $\text{det}(A) \neq 0$.  \\

\noindent
In the case $\text{det}(A) \neq 0$ we know $A^{-1} =  \frac{1}{\text{det}(A)} \text{adj}(A)^T$ in $\RN^{n \times n}$. It remains\footnote{in principle, you could worry the inverse exists in $\mathbb{R}^{n \times n}$ yet is not inside the regular representation of $\Acal$, the argument to follow shows this worry is needless.} to show $A^{-1} \in \MatA(\beta)$.  Linear algebra provides the existence of $T: \Acal \rightarrow \Acal$ such that $[T]_{\beta,\beta} = A^{-1}$. Note,
\begin{equation}
  [T]_{\beta,\beta}[S]_{\beta,\beta}=I  \ \ \Rightarrow \ \ \  T \comp S = Id.
 \end{equation}
Suppose $x,y \in \Acal$ and observe by right-$\Acal$-linearity of $S$ we derive:
\begin{equation}
 T(S(x) \star y) = T(S(x \star y)) = x \star y = T(S(x)) \star y. 
 \end{equation}
Since $S$ is a surjection the calculation above shows $T$ is right-$\Acal$-linear hence $[T]_{\beta,\beta}=A^{-1} \in \MatA(\beta)$. Thus $A$ is a unit in $\MatA(\beta)$. \\

\noindent
On the other hand, if $\text{det}(A)=0$ then the constant term in the minimal polynomial $m(t)$ of $A$ is zero. Thus, either $A=0$ or there exists some $k \geq 2$ for which $m(t) = t^k+ c_{k-1}t^{k-1}+ \cdots + c_1 t$ for some $c_1, \cdots , c_{k-1} \in \RN$. From the theory of the minimal polynomial we know $m(A)=0$ hence as $A$ factors either to the left or right:
\begin{equation}
 0 = m(A) = 
A\left(A^{k-1}+ c_{k-1}A^{k-2}+ \cdots + c_1 I\right) = 
\left(A^{k-1}+ c_{k-1}A^{k-2}+ \cdots + c_1I\right) A.
\end{equation}
Thus $B =A^{k-1}+ c_{k-1}A^{k-2}+ \cdots + c_1I$ gives $AB=0=BA$. Furthermore, since $\MatA(\beta)$ is an algebra and $A \in \MatA(\beta)$ it is clear that $B \in \MatA(\beta)$. Thus in the case $\text{det}(A)=0$ either $A=0$ or $A$ is a zero-divisor.  $\Box$ \\

\noindent
The minimal polynomial argument used in the zero-divisor case could also have been used to show the inverse of a unit $A$ in $\MatA(\beta)$ is formed from a suitable polynomial in $A$. 

\begin{coro} \label{thm:classicalring}
Every element of $\Acal$ and $\EndA$ is either zero, a unit, or a zero-divisor.
\end{coro}

\noindent
{\bf Proof:} combine Theorem \ref{thm:isomorphismpreservesstructure} and Theorem \ref{thm:MatAclassicalring}. $\Box$ \\

\noindent
We study finite dimensional real associative algebras in this work. If $\text{dim}(\Acal)=n$ then geometrically $\Acal$ is essentially just $\RN^n$. For $\MatA(\beta)$ we have an $n$-dimensional subspace of $\RN^{n \times n}$. The set $\mathbf{zd}(\MatA(\beta))$ is the solution set of $\text{det}(A)=0$. This is an $n$-th order polynomial equation in the components of $A$ which includes $A=0$ and at most an $(n-1)$-dimensional space of zero-divisors. For example, if we consider the complex numbers $\Acal = \CN$ then $\Acalzd = \{ 0 \}$. On the other hand $\RN^2$ with the standard direct product given by $(a,b) \star (c,d) = (ac, bd)$ has $\mathbf{zd}(\RN^2) = (\{0 \} \times \RN) \cup (\RN \times \{ 0 \})$. We observe that $\Acalzd$ is formed by a union of subspaces of $\Acal$. This is natural given the following:

\begin{thm} \label{thm:zerodivisorsgeometry}
The set $\Acalzd$ is fixed under negation.
\end{thm}

\noindent
{\bf Proof:} suppose $ x \in \Acalzd$ then there exists $y \in \Acal$ for which $x \star y = 0$. Thus $-x \star y = 0$ and we find $-x \in \Acalzd$. $\Box$

\begin{thm} \label{thm:invertiblebasis}
Let $\Acal$ be an $n$-dimensional real associative unital algebra. There exists an {\bf invertible basis} $\beta = \{ v_1, \dots , v_n \} \subset \Acalx$. 
\end{thm}

\noindent
{\bf Proof:} Let $\gamma =\{ w_1, \dots , w_n \}$ be a basis for $\Acal$. Suppose that a basis element $w_i$ is a zero divisor. Recall that the zero divisors in an $n$ dimensional algebra can be at most $n-1$ dimensional. Hence, there exist $c_1, \dots, c_n$ with $c_i \neq 0$ (since this restriction only removes another $n-1$ dimensional subspace of the algebra) such that $\{ w_1, \dots, c_1 w_1 + \dots + c_i w_i + \dots+ c_n w_n, \dots, w_n \}$ is a basis for $\Acal$ where the $i$-th component is now a unit, since the transformation $w_i \mapsto c_1 w_1 + \dots + c_i w_i + \dots+ c_n w_n$ preserves linear independence of the basis so long as $c_i \neq 0$, and we know by linear algebra that the span of the vectors must also be preserved, since we have $n$ linearly independent vectors in an $n$ dimensional vector space. Therefore, applying this procedure iteratively to $w_1, w_2, \dots, w_n$ yields a basis $\beta = \{v_1, v_2, \dots, v_n\}$ for $\Acal$ where $\beta^* \subset \Acalx$  $\Box$ \\


\noindent
It is easy to see the argument above also allows the following result:

\begin{coro} \label{thm:classicalring}
Let $\Acal$ be an $n$-dimensional real associative unital algebra. There exists an  invertible basis of the special form $\beta = \{ \mathds{1},v_2, \dots , v_n \} \subset \Acalx$. 
\end{coro}

\noindent
Finally, we make the following observation that any open ball about a point in $\Acal$ necessarily intersects infinitely many points in $\Acalx$. This observation should be geometrically evident since $\Acalzd$ is a space of smaller dimension than $\Acal$ and $\Acal - \Acalzd = \Acalx$. 

\begin{thm} \label{thm:unitsdense}
Let $\Acal$ be an $n$-dimensional real associative unital algebra. The group of units $\Acalx$ is a dense subset of $\Acal$.  
\end{thm}

\subsection{Examples}

\noindent
To explain the structure of complex numbers it suffices to say $i^2=-1$ and then just add and multiply $a+bi, c+di$ as usual. Of course, we can be more explicit in our construction if the audience knows about field extensions or group algebras, but, as a starting point it is convenient to provide definitions of algebras which are accessible to every level of student.  

\begin{exa} \label{Ex:number0}
The {\bf real numbers} with their usual addition and multiplication is an associative algebra over $\RN$. If $a \in \RN$ then $[a] \in \text{M}_{\RN} = \RN^{ 1 \times 1}$ is its left regular representation. Usually we will not distinguish between $a$ and $[a]$.
\end{exa}

\begin{exa} \label{Ex:number1}
The {\bf complex numbers} are defined by $\CN = \RN \oplus i\RN$ where $i^2=-1$. If $a+ib,c+id \in \CN$ then $(a+ib)(c+id) = ac+iad+ibc+i^2bd = ac-bd+i(ad+bc)$. Note every nonzero complex number $a+ib$ has multiplicative inverse $\frac{a-ib}{a^2+b^2}$ hence $\CN$ is a field. 
Note  $\mathbf{M}(a+ib) = \left[ \begin{array}{cc} a  & -b  \\  b & a \end{array} \right]$ and the set of all such matrices is denoted $\text{M}_{\CN}$.
\end{exa}

\begin{exa} \label{Ex:number2}
The {\bf hyperbolic} numbers are given by $\calH = \RN \oplus j\RN$ where $j^2=1$. If $a+jb,c+jd \in \calH$ then $(a+jb)(c+jd) = ac+adj+jbc+j^2bd = ac+bd+j(ad+bc)$. Observe $ \mathbf{M}(a+bj) = \left[ \begin{array}{cc} a  & b  \\  b & a \end{array} \right] \in \text{M}_{\calH}$ and $\mathbf{zd}( \calH ) = \{ a+bj \ | \ a^2=b^2 \}$ whereas $\calH^{\times} = \{ a+bj \ | \ a^2 \neq b^2 \}$. The reciprocal of an element in $\calH^{ \times}$ is simply
\begin{equation} 
\frac{1}{a+bj} = \frac{a-bj}{a^2-b^2} 
\end{equation}
this follows from the identity $(a+bj)(a-bj) = a^2-b^2$ given $a^2-b^2 \neq 0$. Let $\Bcal = \RN \times \RN$ with $(a,b)(c,d) = (ac, bd)$ for all $(a,b),(c,d) \in \Bcal$. We can show that
\begin{equation} \Psi(a,b) = a\left(\frac{1+j}{2}\right)+b\left(\frac{1-j}{2}\right) \qquad \& \qquad
\Psi^{-1}(x+jy) = (x+y,x-y) 
\end{equation} 
provide an isomorphism of $\calH$ and $\RN \times \RN$. We can use this isomorphism to transfer problems from $\calH$ to $\Bcal$ and vice-versa. For example, to solve $z^2+Bz+C=0$ in the hyperbolic numbers we note
\begin{equation} 
z^2+Bz+C= 0 \ \ \Rightarrow \ \ \Psi^{-1}(z)^2+
\Psi^{-1}(B)\Psi^{-1}(z)+ \Psi^{-1}(C) = 0 
\end{equation}
Setting $\Psi^{-1}(B)=(b_1,b_2)$ and $\Psi^{-1}(C) = (c_1,c_2)$ and $\Psi^{-1}(z) = (x,y)$ we arrive at 
\begin{equation}
 (x,y)^2+(b_1,b_2)(x,y)+(c_1,c_2) = 0 
 \end{equation}
which reduces to
\begin{equation}
 (x^2+b_1x+c_1, y^2+b_2y+c_2) = (0,0). 
 \end{equation}
Of course, these are just quadratic equations in $\RN$ so we can solve them and transfer back the result to the general solution of $z^2+Bz+C=0$ in $\calH$. Given this correspondence, we deduce there are either zero, two or four solutions to the quadratic hyperbolic equation.
\end{exa}

\begin{exa} \label{Ex:number3}
The {\bf dual numbers} are given by $\calN = \RN \oplus \eps \RN$ where $\eps^2=0$. If $a+\eps b,c+\eps d \in \calN$ then 
\begin{equation}
 (a+\eps b)(c+\eps d) = ac+ad\eps+bc\eps+\eps^2bd = ac+(ad+bc)\eps. 
\end{equation}
Observe $\mathbf{M}(a+b \eps) = \left[ \begin{array}{cc} a  & 0  \\  b & a \end{array} \right] \in \text{M}_{\calN}$ and $\mathbf{zd}( \calN) = \{ a+b\eps \ | \ a^2=0 \} = \eps \RN$. The units in the dual numbers are of the form $a+b\eps$ where $a \neq 0$. Note $(a+b \eps)(a-b \eps) = a^2$ hence $\frac{1}{a+b \eps} = \frac{a-b \eps}{a^2}$ provided $a \neq 0$.
\end{exa}

\noindent
For higher dimensional algebras the multiplicative inverse of a general element can be calculated by computing the inverse of the element's regular representation.

\begin{exa} \label{Ex:number3n}
The {\bf $n$-th order dual numbers} are given by $ \Ncal_n = \RN \oplus \eta\RN \oplus \cdots \oplus \epsilon^{n-1}\RN$ where $\eps^n=0$ and $\eps^k \neq 0$ for $1 \leq k \leq n-1$. The regular representation is formed by lower triangular matrices of a particular type:
\begin{equation}
 \mathbf{M}(a_1+a_2\eps+ \cdots + a_n \eps^{n-1}) = \left[ 
\begin{array}{llllll} 
a_1  & 0 & \cdots & 0 & 0  \\  
a_2  & a_1 & \cdots & 0 & 0  \\
\vdots & \vdots & \ddots  & \vdots & \vdots \\
a_{n-1} & a_{n-2} & \cdots & a_1 & 0  \\
a_n & a_{n-1} & \cdots & a_2 & a_1  
\end{array} \right] \in \text{M}_{\calN_n} 
\end{equation}
Notice $a_1+a_2\eps + \cdots + a_n \eps^{n-1} \in\mathbf{zd}( \calN_n )$ only if $a_1 \neq 0$.
\end{exa}

\begin{exa} \label{Ex:number4}
Let $\mathcal{A} = \RN \oplus j\RN \oplus j^2 \RN$ where $j^3=1$. The matrix representatives of these numbers have an interesting shape; note: $A \in \text{M}_{\calA}$ implies $A = \left[ \begin{array}{ccc} a  & c & b  \\  b & a & c \\ c & b & a \end{array} \right]$. We note an isomorphism $\calA \approx \RN \times \CN$ is given by mapping $j$ to $(1, \omega)$ where $\omega$ is a third root of unity.
\end{exa}

\begin{exa} \label{Ex:number5}
Let $\mathcal{A} = \RN \times \calH$ where $\mathds{1} = (1,1+0j)$. Let $\beta = \{ (1,0), (0,1), (0,j) \}$ gives block-diagonal $A \in \text{M}_{\calA}(\beta)$; 
\begin{equation} 
\mathbf{M}_{\beta}((a,b+cj)) = \left[ \begin{array}{ccc} a  & 0 & 0  \\  0 & b & c \\ 0& c & b \end{array} \right]. 
\end{equation} 
This algebra is isomorphic to $\RN \times \RN \times  \RN$ with $(a_1,a_2,a_3)\star (b_1,b_2,b_3)=(a_1b_1,a_2b_2,a_3b_3)$.
\end{exa}

\begin{exa} \label{Ex:number6}
Let $\mathcal{A} = \RN \oplus j\RN \oplus j^2 \RN\oplus j^3\RN$ where $j^4=1$. Observe,
\begin{equation}
 \mathbf{M}(a+bj+cj^2+dj^3) = \left[ \begin{array}{cccc} a  &  d  & c & b \\  b & a & d & c \\ c & b & a & d \\ d & c & b & a \end{array} \right]. \end{equation}
 This algebra is naturally isomorphic to $\CN \oplus \calH$ which is clearly isomorphic to $\CN \times  \RN \times  \RN$.
\end{exa}

\begin{exa} \label{Ex:number7}
Let $\mathcal{A} = \calH \times \calH$ where $\mathds{1} = (1+0j,1+0j)$. This means $(1,1)$ is naturally represented by the identity matrix. Set $\beta = \{ (1,0), (j,0), (0,1), (0,j) \}$ and observe
\begin{equation}
 \mathbf{M}_{\beta}((a+bj,c+dj)) = \left[ \begin{array}{cc|cc} a  &  b  & 0 & 0 \\  b & a & 0 & 0 \\ \hline 0 & 0 & c & d \\ 0 & 0 & d & c \end{array} \right]. 
 \end{equation}
This algebra is isomorphic to $\RN \times \RN \times  \RN\times \RN$ with the Hadamard product $(a_1,a_2,a_3,a_4)*(b_1,b_2,b_3,b_4)=(a_1b_1,a_2b_2,a_3b_3,a_4b_4)$.
\end{exa}

\begin{exa} \label{Ex:number8}
Let $\mathcal{A} = \CN \times \CN$ where $\mathds{1} = (1+0i,1+0i)$. Here we study the problem of two complex variables. In this algebra $(1+0i, 1+0i)$ corresponds to the identity and hence $(1,1)$ is naturally represented by the identity matrix. In total we have once more a block-diagonal representation: $A \in \text{M}_{\calA}$ implies $A = \left[ \begin{array}{cc|cc} a  &  -b  & 0 & 0 \\  b & a & 0 & 0 \\ \hline 0 & 0 & c & -d \\ 0 & 0 & d & c \end{array} \right]$ and this matrix represents $(a+bi,c+di)$.
\end{exa}

\begin{exa} \label{Ex:number9}
Let $\mathbb{H} = \RN \oplus i\RN \oplus j \RN\oplus k\RN$ where $i^2=j^2=k^2=-1$ and $ij=k$.  These are Hamilton's famed {\bf quaternions}. We can show $ij=-ji$ hence these are not commutative. With respect to the natural basis $e_1=1, e_2=i, e_3=j,e_4=k$ we find the matrix representative of $a+ib+cj+dk$ is as follows: 
\begin{equation}
A= \left[ \begin{array}{rrrr} a & -b & -c & -d \\ b & a & -d & c \\ c & d & a & -b \\ d & -c & b & a  \end{array} \right] \in \text{M}_{\mathbb{H}}.
\end{equation}
\end{exa}

\begin{exa} \label{Ex:number10}
Let $\mathcal{A} = \RN_2$ with the multiplication $\star$ induced from the multiplication of $2 \times 2$ matrices.  This again forms a noncommutative algebra. In particular, this multiplication is induced in the natural manner:  
\begin{equation}
 \left[\begin{array}{cc} a & b \\ c & d \end{array} \right]\left[\begin{array}{cc} t & x \\ y & z \end{array} \right] = \left[\begin{array}{cc} at+by & ax+bz \\ ct+dy & cx+dz \end{array} \right]. 
 \end{equation}
It follows that $ (a,b,c,d) \star (t,x,y,z) = ( at+by, ax+bz,ct+dy,cx+dz)$. We can read from this multiplication that the representative of $(a,b,c,d) \in \RN_2$ is given by
\begin{equation}
 A= \left[\begin{array}{cc|cc} a & 0 & b & 0 \\ 0 & a & 0 & b \\ \hline c & 0 & d & 0 \\ 0 & c & 0 & d \end{array} \right] = \left[\begin{array}{c|c} aI & bI \\ \hline cI & dI \end{array} \right] \in \Amat. 
\end{equation}
Note, the basis $\beta = \{ E_{11}, E_{12}, E_{21}, E_{22} \}$ does not contain the multiplicative identity $I = E_{11}+E_{22}$. We define $\Bcal = \RN^4$ by
\begin{equation}
 (a,b,c,d) \star (w,x,y,z) = (aw+by,  ax+bz,  cw+dy, cx+dz) 
 \end{equation}
Of course, $(1,0,1,0) = \mathds{1}_{\Bcal}$ and $\Bcal$ is really just another notation for the $2 \times 2$ matrix algebra. In fact, $\Psi \left[ \begin{array}{cc} a & b \\ c & d \end{array} \right] = (a,b,c,d)$ defines an isomorphism of $\Acal$ and $\Bcal$. The reader will verify that $\Psi (AB) = \Psi(A) \star \Psi(B)$.

\end{exa}

\begin{exa}
 $Z$ in $\HN^{2 \times 2}$ is a $\Acal$-number. There is a natural injection $ \Psi: \HN^{ 2 \times 2} \rightarrow \RN^{ 8 \times 8}$ induced from $\mathbf{M}: \mathbb{H} \rightarrow \RN^{  4\time 4}$ from Example \ref{Ex:number9}, 
\begin{equation}
\Psi \left( \left[\begin{array}{cc} x & y \\ 
z & w \end{array} \right] \right) = 
\left[\begin{array}{cc} \mathbf{M}(x) & \mathbf{M}(y) \\ 
\mathbf{M}(z) & 
\mathbf{M}(w) \end{array} \right]
\end{equation}
for all $x,y,z,w \in \mathbb{H}$. The matrices in $\Psi( \HN^{ 2 \times 2})$ are isomorphic to $\HN^{ 2 \times 2}.$ 
\end{exa}

\noindent

\begin{exa}
Let $G$ be a finite multiplicative group; $G = \{ g_1, \dots , g_n \}$ then we define
\begin{equation} 
\Acal_G = g_1\RN \oplus \cdots \oplus g_n \RN 
\end{equation}
with natural multiplication inherited from $G$. For example, for all $a,b,c,d \in \RN$,
\begin{equation}
 (ag_1+bg_2)(cg_3+dg_4) = acg_1g_3+adg_1g_4 +bcg_2g_3+bdg_2g_4. 
 \end{equation}
The group algebra allows us to multiply $\RN$-linear combinations of group elements by extending the group multiplication linearly. By construction, $\{g_1, \dots, g_n \}$ serves as a basis for $\Acal_G$. As $G$ is a group we know for each $i,j \in \{ 1, \dots , n \}$ there exists $k \in \{ 1, \dots , n \}$ for which $g_ig_j=g_k$. If we define structure constants $C_{ijk}$ by $g_ig_j = \sum_l C_{ijl}g_l$ then $g_ig_j=g_k$ implies $C_{ijl} = \delta_{kl}$. 
\end{exa}

\begin{exa}
The cyclic group of order $n$ in multiplicative notation has the form $G = \{ e, g,g^2, \dots, g^{n-1} \}$. The group algebra $\Acal_G = e\RN \oplus g \RN \oplus \cdots \oplus g^{n-1} \RN$. We usually call this algebra the $n$-hyperbolic numbers. 
\end{exa}

\section{$\Acal$-differentiable functions} \label{sec:adifffunctions}
\noindent
Given $\Acal$ with basis $\beta = \{ \mathds{1}, v_2, \dots , v_n \}$ we may define an inner-product on $\Acal$ by bilinearly extending $g(v_i,v_j) = \delta_{ij}$. The $g$-induced norm $|| x|| = \sqrt{ g(x,x)}$ has $||v_i||=1$ for $i = 1,2, \dots , n$. In this construction we have $\beta$ is $g$-orthonormal and
\begin{equation} 
|| x_1v_1+x_2v_2+ \cdots + x_nv_n ||^2 = x_1^2+x_2^2+  \cdots + x_n^2.
\end{equation}
hence $|| \zeta || = || \overline{\zeta}_j ||$ for $j=2,\dots , n$. \\

\noindent
We should caution, there are examples where the length of $\mathds{1}$ is not $1$ in the natural norm for the example. For example, $\RN^n$ with $\mathds{1} = (1,\dots , 1)$ has length $|| \mathds{1} || = \sqrt{n}$ in the usual Euclidean metric. \\

\noindent
The multiplicative structure of the absolute value on $\RN$, the modulus on $\CN$ or $\HN$ are very special. It is not typically the case that we have $|| x \star y|| = ||x|| \, ||y||$, however, we can always find a norm which is submultiplicative over $\Acal$.

\begin{thm} \label{thm:submultiplicative}
If $\Acal$ is an associative $n$-dimensional algebra over $\RN$ then there exists a norm $|| \cdot ||$ for $\Acal$ and $\Abound >0$ for which $|| x \star y|| \leq \Abound ||x|| ||y||$ for all $x,y \in \Acal$. Moreover, for this norm we find $\Abound = \mathbf{C}(n^2-n+1) \sqrt{n}$ where $\mathbf{C} = \text{max} \{ C_{ijk} \ | \ 1 \leq i,j,k \leq n \}$.
\end{thm}

\noindent
{\bf Proof:} suppose $\beta = \{ v_1, \dots , v_n \}$ is a basis for $\Acal$ for which $|| x_1v_1+ \cdots +x_nv_n ||^2 = x_1^2+ \cdots +x_n^2$ for each $x_1v_1+ \cdots +x_nv_n \in \Acal$. It is always possible to construct such a norm as we described at the outset of this section. Next, suppose $v_i \star v_j = \sum_{i,j,k}C_{ijk}v_k$ and define $\mathbf{C} = \text{max} \{ C_{ijk} \ | \ 1 \leq i,j,k \leq n \}$. Let $x,y \in \Acal$ where $x = \sum_i x_i v_i$ and $y = \sum_j y_jv_j$. Calculate,
\begin{align} \label{eqn:submultiplicative}
|| x \star y ||^2 &= \big{|}\big{|} \sum_{i} x_iv_i  \star  \sum_{j} y_jv_j  \big{|}\big{|}^2 \\ \notag
&= \big{|}\big{|} \sum_{i,j,k} x_i v_j C_{ijk} v_k \big{|}\big{|}^2 \\ \notag
&= \sum_k \biggl( \sum_{i,j} C_{ijk} x_i y_j \biggr)^2 \\ \notag
&\leq \sum_k  \biggl( \mathbf{C} \sum_{i,j} | x_i| |y_j| \biggr)^2 \\ \notag
&= n \mathbf{C}^2\biggl(\sum_{i,j}| x_i| |y_j| \biggr)^2 \\ \notag
&= n \mathbf{C}^2\biggl(\sum_{i}| x_i|\biggr)^2  \biggl(\sum_{j}|y_j| \biggr)^2 \\ \notag
&= n \mathbf{C}^2\biggl(\sum_{i}| x_i|^2+ \sum_{k \neq l} |x_k||x_l| \biggr) \biggl(\sum_{j}| y_j|^2+\sum_{k \neq l}|y_k||y_l| \biggr) \\ \notag
&\leq n \mathbf{C}^2\biggl(||x||^2+ \sum_{k \neq l} ||x||^2 \biggr) \biggl(||y||^2+\sum_{k \neq l}||y||^2 \biggr) \\ \notag
&= n(n^2-n+1)^2 \mathbf{C}^2||x||^2||y||^2.
\end{align}
Thus, set $\Abound = \mathbf{C}(n^2-n+1) \sqrt{n}$ as to obtain $||x \star y|| \leq \Abound ||x|| \, ||y||$ for all $x,y \in \Acal$. $\Box$ \\

\noindent
Typical applications have $\mathbf{C}=1$ and we find the bound found in the proof above can usually be sharpened for a particular algebra.

\begin{exa}
In the hyperbolic numbers $\Hcal = \RN \oplus j \RN$ it can be shown geometrically $|| z \star w|| \leq \sqrt{2} ||z|| \, ||w||$ for all $z,w \in \Hcal$ where $||x+jy|| = \sqrt{x^2+y^2}$ defines the Euclidean norm. In fact, this bound cannot be made smaller. Consider $(1+j)(1+j) = 2+2j$ and $||1+j|| = \sqrt{2}$ whereas $||2+2j|| = \sqrt{8} = \sqrt{2}\sqrt{2}\sqrt{2}$ hence $||(1+j)(1+j)|| = \sqrt{2}||1+j|| \cdot ||1+j||$. In this case the structure constants take values of $0$ or $1$ in all cases thus $\mathbf{C}=1$ and as $n=2$ we find $\Abound = (4-2+1)\sqrt{2} = 3 \sqrt{2}$. Thus $\Abound$ is not sharp.
\end{exa}

\noindent
Also, beware that $||a/b|| = ||a||/||b||$ in general\footnote{we assume $\Acal$ is commutative and thus denote both $a \star b^{-1}$ and $b^{-1} \star a$ by $a/b$. }. However, we can say something productive:

\begin{coro} \label{thm:quotientinequality}
Suppose $\Abound >0$ is a real constant such that $|| x \star y|| \leq \Abound ||x|| ||y||$ for all $x,y \in \Acal$. If $b \in \Acalx$ and $a \in \Acal$ then 
\begin{equation}
 \frac{||a||}{||b||} \leq \Abound \, \bigg{|}\bigg{|} \frac{a}{b} \bigg{|}\bigg{|}. 
 \end{equation}
\end{coro}

\noindent
{\bf Proof:} Theorem \ref{thm:submultiplicative} provides $\Abound>0$ for which $|| x \star y|| \leq \Abound ||x|| ||y||$ for all $x,y \in \Acal$. Consider, if $b \in \Acalx$ and $a \in \Acal$ then $a/b \in \Acal$ and $a = b \star (a/b)$ thus \\ $ ||a|| = ||b \star (a/b)|| \leq \Abound \, ||b|| \, \big{|}\big{|} \frac{a}{b}\big{|}\big{|}$ and we deduce $\frac{||a||}{||b||} \leq \Abound \, \big{|}\big{|} \frac{a}{b} \big{|}\big{|}$.  $\Box$ \\

\noindent
In what follows we assume $\Acal$ has a norm denoted $|| \cdot ||$. For a given $\Acal$ there may be many choices for the norm, however, these all produce the same topology  as $\Acal$ is a finite dimensional real vector space.

\begin{de} \label{defn:Adiff}
Let $U \subseteq \Acal$ be an open set containing $p$. If $f: U \rightarrow \Acal$ is a function then we say $f$ is {\bf $\Acal$-differentiable at $p$} if there exists a linear function $d_pf \in \EndA$ such that
\begin{equation} \label{eqn:frechetquotAcal}
 \lim_{h \rightarrow 0}\frac{f(p+h)-f(p)-d_pf(h)}{||h||} = 0. 
\end{equation}
 \end{de}

\noindent
Recall, $d_pf \in \EndA$ implies $d_pf: \Acal \rightarrow \Acal$ is $\RN$-linear mapping on $\Acal$ and $d_pf( v \star w) = d_pf(v) \star w$ for all $v,w \in \Acal$.

\begin{thm} \label{thm:AdiffisFrechetDiff}
If $f$ is $\Acal$ differentiable at $p$ then $f$ is $\RN$-differentiable at $p$.
\end{thm}

\noindent
{\bf Proof:} if $f$ is $\Acal$-differentiable at $p$ then we know $d_pf \in \EndA$ satisfies the Frechet limit given in Equation \ref{eqn:frechetquotAcal}. Moreover, $d_pf$ is $\RN$-linear. $\Box$ \\ 

\noindent
There are $\RN$-differentiable functions which are not $\Acal$-differentiable. The condition $d_pf \in \EndA$ is not met by all functions on $\Acal$. 

\begin{exa}
Let $\Acal$ be an algebra of dimension $n \geq 2$ and $ \{v_1, \dots , v_n \}$ an invertible basis with $v_1 = \mathds{1}$ and coordinates $x_1, \dots , x_n$ then we have $\zeta = x_1+ \cdots +x_nv_n$ and $\overline{\zeta}_2 = x_1-v_2x_2+ \cdots +x_nv_n$. The function $f(\zeta) = \overline{\zeta}_2$ is everywhere real differentiable and nowhere $\Acal$-differentiable. These observations are most easily verified using the tools of Section \ref{sec:CReqns} which culminate in Theorem \ref{thm:CReqnsII}.
\end{exa}

\noindent
If $f$ is $\Acal$-differentiable at each $p \in V$ then $f$ is {\bf $\Acal$-differentiable on $V$} and we write $f \in \Fun(V)$. If there exists an open set containing $p$ on which $f$ is $\Acal$-differentiable then we say $f$ is {\bf $\Acal$-differentiable near $p$} and write $f \in \Fun(p)$. There are several ways to characterize $\Acal$-differentiability at a point. These all follow from the isomorphism $\Acal \approx \EndA \approx \MatA(\beta)$.

\begin{thm} \label{thm:AcalCRequations}
Let $U \subseteq \Acal$ and $p \in U$.  Let $f: U \rightarrow \Acal$ be a $\RN$-differentiable function at $p$. The following are equivalent
\begin{enumerate}[{\bf (i.)}]
\item $d_pf( v \star w) = d_pf(v) \star w$ for all $v,w \in \Acal$,
\item there exists $\lambda \in \Acal$ for which $d_pf(v) = \lambda \star v$ for each $v \in \Acal$,
\item for any basis $\beta$ of $\Acal$, $[d_pf]_{\beta, \beta} \in \MatA(\beta)$.
\end{enumerate}
\end{thm}

\noindent
{\bf Proof:} If $f$ is a function on $\Acal$ which is $\RN$-differentiable at $p$ then there exists an $\RN$-linear function $d_pf: \Acal \rightarrow \Acal$ which satisfies the Frechet quotient condition set-forth in Equation \ref{eqn:diffwrtR}. Notice condition (i.) is true iff $d_pf \in \EndA$. Therefore, (i.) is equivalent to (iii.) in view of the definition of the regular representation built from $\beta$ (see Definition \ref{defn:regrepofA}). Suppose (i.) is true. Since $\Acal$ is unital we have $v = \mathds{1} \star v$,
\begin{equation}
d_pf(v) = d_pf(\mathds{1} \star v) = d_pf(\mathds{1}) \star v 
\end{equation}
hence (ii.) follows with $d_pf(\mathds{1} ) = \lambda$. Suppose (ii.) is true and let $\beta$ be a basis for $\Acal$. Let $v,w \in \Acal$ and consider
\begin{equation}
 d_pf(v \star w) = \lambda \star (v \star w) = (\lambda \star v) \star w = d_pf(v) \star w.
\end{equation}
thus (i.) is true. Since (i.) $\Leftrightarrow$ (ii.) and (i.) $\Leftrightarrow$ (iii.) we have (i.) $\Leftrightarrow$ (ii.) $\Leftrightarrow$ (iii.).  $\Box$ \\
                   
\noindent
The $\Acal$-number $\lambda$ which appears above in (ii.) is known as the {\bf derivative} of $f$ at $p$. Notice: if $d_pf(v) = \lambda \star v$ for each $v \in \Acal$ then $\#(d_pf) = \lambda$. We should appreciate the importance of the isomorphism of $\EndA$ and $\Acal$ as it allows derivatives of functions on $\Acal$ to be viewed once more as functions on $\Acal$. This is a great simplification as derivatives of $\RN$-differentiable maps are not usually objects of the same type. The following quote\footnote{from the introduction to Dieudonn\a' e's chapter on differentiation in Modern Analysis Chapter VIII} is from Dieudonn\a' e in \cite{Dmaster}
\begin{quote}
...on a one-dimensional vector space, there is a one-to-one correspondence between linear forms and numbers, and therefore the derivative at a point is defined as a number instead of a linear form. 
\end{quote}
Dieudonn\a' e says this to encourage students to place linear transformations at the center stage of their analysis. In contrast, we find the correspondence of right-$\Acal$-linear transformations, or later $k$-linear transformations on $\Acal$\footnote{in particular, see Theorem \ref{thm:multilinearisomorphictoA}}, with $\Acal$ itself allows us to perform calculations in $\Acal$-calculus in nearly the same fashion as introductory real calculus. In other words, since there is also a natural correspondence between $\Acal$-linear transformations and $\Acal$ we escape the sophistication which Dieudonn\a' e could not avoid.

\begin{de} \label{defn:derivative}
Let $U\subseteq \Acal$ be an open set and $f: U  \rightarrow \Acal$ an $\Acal$-differentiable function on $U$ then we define $f': U \rightarrow \Acal$ by $f'(p) = \# (d_pf) $ for each $p \in U$.
\end{de}

\noindent
Equivalently, we could write $f'(p) = d_pf( \mathds{1})$ since $\#(T) = T( \mathds{1})$ for each $T \in \EndA$.
Many theorems of calculus hold for $\Acal$-differentiable functions. 

\begin{thm} \label{thm:Acalculusbasicprops}
If $f,g \in \Fun(p)$ and $c \in \Acal$ and we define $f+g$ by the usual rule and  $(c\star f)(x) = c\star f(x)$ for each $x \in \text{dom}(f)$. Then
\begin{enumerate}[{\bf (i.)}]
\item $f+g \in \Fun(p)$ and $ \ds (f+g)'(p) = f'(p)+g'(p)$,
\item $c \star f \in \Fun(p)$ and $ \ds (c \star f)'(p) = c \star f'(p)$.
\end{enumerate}
\end{thm}

\noindent
{\bf Proof (i.):} suppose $f,g \in \Fun(p)$ and $c \in \Acal$. Recall from the advanced calculus of normed linear spaces that real differentiability of $f,g$ at $p$ implies real-differentiablitly at $p$ of $f+g$ and $d_p(f+g) = d_pf+d_pg$. But, we assume $f,g \in \Fun(p)$ hence $d_pf,d_pg \in \EndA$ and by Theorem \ref{thm:subalgebraofleftlinearmaps} we deduce $d_p(f+g) \in \EndA$ and thus $f+g \in \Fun(p)$. The number map is linear hence
$d_p(f+g) = d_pf+d_pg$ implies $\#(d_p(f+g)) = \#(d_pf)+\#(d_pg)$ which gives $(f+g)'(p) = f'(p)+g'(p)$ which proves (i.). \\

\noindent
{\bf Proof(ii.):} If $f \in \Fun(p)$ then $d_pf \in \EndA$ which means $d_pf(v \star w) = d_pf(v) \star w$. Let $c \in \Acal$ and define $g(p) = c \star f(p)$. Let $L(h) = c \star d_pf(h)$ for each $h \in \Acal$.  If $v, w \in \Acal$ then
\begin{equation}
 L(v \star w) = c \star d_pf(v \star w) = c \star d_pf(v) \star w = L(v) \star w. 
\end{equation}
thus $L \in \EndA$. It remains to show $g$ is differentiable with $d_pg=L$. If $h \neq 0$ let $\mathcal{F}_f, \mathcal{F}_g$ denote the Frechet quotients of $f,g$ respective. Since $f$ is differentiable at $p$ means $\lim_{h \rightarrow 0} \mathcal{F}_f =0$. Calculate:
\begin{align} 
\mathcal{F}_g =\frac{g(p+h)-g(p)-L(h)}{||h||} &= \frac{c\star f(p+h)-c\star f(p)-c\star d_pf(h)}{||h||} \\ \notag
&= c \star \frac{ f(p+h)- f(p)-d_pf(h)}{||h||} \\ \notag
&= c \star \mathcal{F}_f.
\end{align}
Thus, by Theorem \ref{thm:submultiplicative} we find $ ||\mathcal{F}_g || = || c \star \mathcal{F}_f|| \leq \Abound ||c|| \, || \mathcal{F}_f||$. Since $|| \mathcal{F}_f|| \rightarrow 0$ as $h \rightarrow 0$ it follows $\lim_{h \rightarrow 0} \mathcal{F}_g =0$. Hence, $g$ is $\RN$-differentiable with $d_pg = L$. Thus $c \star f \in \Fun(p)$ with $d_p(c \star f) = c \star d_pf$. Note 
\begin{equation}
d_p(c \star f)(\mathds{1}) = c \star d_pf(\mathds{1}) = c \star f'(p)
\end{equation}
Consequently, $\#(d_p(c \star f)) = c \star f'(p)$ and we conclude $(c \star f)'(p) = c \star f'(p)$. $\Box$ \\

\noindent
The product of two $\Acal$-differentiable functions is not necessarily $\Acal$-differentiable in the case that $\Acal$ is a noncommutative algebra. However, the product of two $\Acal$-differentiable functions is always $\RN$-differentiable and we have the following result:

\begin{thm} \label{thm:Aproductrule}
Suppose $f,g$ are $\Acal$-differentiable at $p$ then 
$$ d_p( f \star g)(v) = d_pf(v) \star g(p)+ f(p) \star d_pg(v) $$
for each $v \in \Acal$. Furthermore, if $\Acal$ is commutative then $f \star g$ is $\Acal$-differentiable at $p$ and
$$ (f \star g)'(p) = f'(p) \star g(p)+ f(p) \star g'(p). $$
\end{thm}

\noindent
{\bf Proof:} suppose $f, g$ are $\Acal$ differentiable at $p$ then $f,g$ are $\RN$-differentiable at $p$ with differentials $d_pf, d_pg \in \EndA$. Furthermore, $f(p+h) = f(p)+d_pf(h) + \eta_f$ where the Frechet quotient $\mathcal{F}_f = \eta_f /||h|| \rightarrow 0$ as $h \rightarrow 0$. Likewise, $g(p+h) = g(p)+d_pg(h) + \eta_g$ where the Frechet quotient $\mathcal{F}_g = \eta_g /||h|| \rightarrow 0$ as $h \rightarrow 0$.  Calculate,
\begin{align}
 f(p+h)\star g(p+h) &= \bigl(f(p)+d_pf(h) + \eta_f \bigr)\star\bigl(g(p)+d_pg(h) + \eta_g \bigr) \\ \notag
 &= f(p) \star g(p) + \underbrace{f(p) \star d_pg(h) + d_pf(h) \star g(p)}_{L(h)} + \eta_{f \star g}
 \end{align}
where $ \eta_{f \star g} = f(p) \star \eta_g+ \eta_f \star g(p)+ \eta_f \star \eta_g $. For $h \neq 0$ we define, then simplify
\begin{equation}
 \mathcal{F}_{f \star g} = \frac{f(p+h)\star g(p+h)-f(p)\star g(p) - L(h)}{||h||}  = \frac{\eta_{f \star g}}{||h||}
\end{equation}
Thus, by the triangle inequality and submultiplicativity of $|| \cdot ||$
\begin{equation}
 ||  \mathcal{F}_{f \star g}|| \leq \Abound\biggl( ||f(p)|| \, \bigg{|}\bigg{|} \frac{\eta_g}{ ||h||} \bigg{|}\bigg{|}+\bigg{|}\bigg{|} \frac{\eta_f}{ ||h||} \bigg{|}\bigg{|}\, ||g(p)|| + \frac{||\eta_f||||\eta_g||}{ ||h||} \biggr) 
\end{equation}
Note, for $h \neq 0$ we have $||h|| \neq 0$ and $\frac{||\eta_f||||\eta_g||}{ ||h||} = ||h|| \frac{||\eta_f||}{||h||}\frac{||\eta_g||}{||h||} = ||h|| \, || \mathcal{F}_f|| \, || \mathcal{F}_g||$ where $\mathcal{F}_f,\mathcal{F}_g$ denote the Frechet quotients of $f,g$ respectively. In summary,
\begin{equation}
 ||  \mathcal{F}_{f \star g}|| \leq \Abound \bigl( ||f(p)|| \, ||\mathcal{F}_g || + ||\mathcal{F}_f || \, ||g(p)||  + ||h|| \, ||\mathcal{F}_f || \, ||\mathcal{F}_g || \bigr) 
 \end{equation}
thus $||  \mathcal{F}_{f \star g}|| \rightarrow 0$ as $h \rightarrow 0$. We defined $L = f(p) \star d_pg(h) + d_pf(h) \star g(p)$ hence $L(cv+w) = cL(v)+L(w)$ for $c \in \RN$ and $v,w \in \Acal$. The real linearity of $L$ follows from linearity of $d_pf$ and $d_pg$ as well as the structure of the $\star$-multiplication. Therefore, $f \star g$ is real-differentiable at $p$ with $d_p(f \star g)(v) = d_pf(v) \star g(p)+ f(p) \star d_pg(v)$. \\

\noindent
Next, suppose $\Acal$ is commutative. If $v, w \in \Acal$ then
\begin{align} 
 d_p(f \star g)(v \star w) &= d_pf(v \star w) \star g(p)+ f(p) \star d_pg(v \star w) \\ \notag
 &= d_pf(v) \star w \star g(p)+ f(p) \star d_pg(v) \star w \\ \notag
 &= \bigl(d_pf(v) \star g(p)+ f(p) \star d_pg(v) \bigr) \star w \\ \notag
 &= d_p(f \star g)(v) \star w
 \end{align}
Thus $d_p(f \star g) \in \EndA$ and we have shown $f \star g$ is $\Acal$-differentiable at $p$. Moreover,
\begin{align} 
(f \star g)'(p) &= d_p(f \star g)(\mathds{1}) \\ \notag
&= d_pf(\mathds{1}) \star g(p)+ f(p) \star d_pg(\mathds{1})\\ \notag
&= f'(p) \star g(p) + f(p) \star g'(p). \ \ \Box
\end{align}

\noindent
If $\Acal$ is not commutative then it is possible
\begin{equation} 
 d_pf(v) \star w \star g(p)  \neq d_pf(v)  \star g(p) \star w 
 \end{equation}
If $f \star g$ is to be $\Acal$ differentiable at $p$ then we need that
\begin{equation} \label{eqn:noncommuteprobI}
 d_pf(v) \star w \star g(p)  - d_pf(v)  \star g(p) \star w = 0 
 \end{equation}
for all $v, w \in \Acal$. Equivalently,
\begin{equation} \label{eqn:noncommuteprobII}
 d_pf(v) \star \bigl( w \star g(p)  -  g(p) \star w \bigr) = 0 
 \end{equation}
For example, if $f(p)=c$ is the constant function then $f \star g$ is differentiable as we already saw in Theorem \ref{thm:Acalculusbasicprops} part (ii.). For a less trivial example, we could seek a function $g$ for which 
\begin{equation} \label{eqn:criterionnilp}
w \star g(p)  -  g(p) \star w =0
\end{equation}
for all $w \in \Acal$ and some $p$. The {\bf center} of $\Acal$ is $\mathbf{Z}( \Acal) = \{ x \in \Acal \ | \ x \star y = y \star x \ \ \text{for all $y \in \Acal$} \}$. The center forms an ideal of $\Acal$. We say $\Acal$ is a simple algebra if it has no ideals except $\{0 \}$ and $\Acal$. If $g(p) \in \mathbf{Z}( \Acal)$ and $g(p) \neq 0$ then we find $\Acal$ is not a simple algebra.  Another aspect to this result is the nonexistence of higher than first-order $\Acal$-polynomials. In particular, Rosenfeld shows in \cite{Rosenfeld} that there are only linear $\Acal$-differentiable functions over simple associative or alternative algebras. Simple algebras aside, there are algebras with nontrivial centers which in turn support nontrivial functions which meet the criteria of Equation \ref{eqn:criterionnilp}.

\begin{exa} \label{exa:uppertriangular}
Let $\Acal = \RN^6$ with the following noncommutative multiplication:
\begin{equation}
 (a,b,c,d,e,f) \star (x,y,z,u,v,w) = (ax, by, cz, au+dy, bv+ez, aw+dv+fz) \end{equation}
The regular representation of $\Acal$ has typical element
\begin{equation}
 \mathbf{M}(a,b,c,d,e,f) = 
\left[
\begin{array}{cccccc}
a & 0 & 0 & 0 & 0 & 0 \\
0 & b & 0 & 0 & 0 & 0 \\
0 & 0 & c & 0 & 0 & 0 \\
0 & d & 0 & a & 0 & 0 \\
0 & 0 & e & 0 & b & 0 \\
0 & 0 & f & 0 & d & a 
\end{array}
\right]
\end{equation}
Suppose $\Acal$ has variables $\zeta = (x_1, \dots, x_6)$ and define $f(\zeta) = (1,1,1,1,1,x_3^2)$ and define $g(\zeta) = (0,0,0,x_2,0,x_5)$. Calculate $(f \star g)(\zeta) = (0,0,0,x_2,0,x_5) $. We calculate,
\begin{equation}
 \left[ \frac{\partial f}{\partial x_i} \right] = \left[
\begin{array}{cccccc}
0 & 0 & 0 & 0 & 0 & 0 \\
0 & 0 & 0 & 0 & 0 & 0 \\
0 & 0 & 0 & 0 & 0 & 0 \\
0 & 0 & 0 & 0 & 0 & 0 \\
0 & 0 & 0 & 0 & 0 & 0 \\
0 & 0 & 2x_3 & 0 & 0 & 0 
\end{array}
\right] \qquad \& \qquad
\left[ \frac{\partial g}{\partial x_i} \right] = \left[
\begin{array}{cccccc}
0 & 0 & 0 & 0 & 0 & 0 \\
0 & 0 & 0 & 0 & 0 & 0 \\
0 & 0 & 0 & 0 & 0 & 0 \\
0 & 1 & 0 & 0 & 0 & 0 \\
0 & 0 & 0 & 0 & 0 & 0 \\
0 & 0 & 0 & 0 & 1 & 0 
\end{array}
\right]
\end{equation}
Observe $f$ and $g$ are $\Acal$-differentiable and $f \star g = g$ is likewise $\Acal$-differentiable. In contrast, $(g \star f)(\zeta) = (0,0,0,x_2,0,x_2+x_5)$ is not $\Acal$-differentiable as its Jacobian matrix is  nonzero in the $(2,6)$-entry and hence is not\footnote{The algebra in Example \ref{exa:uppertriangular} is isomorphic to the algebra formed by upper triangular matrices in $\RN^{ 3 \times 3}$. The center of the upper triangular matrices is formed by the strictly upper triangular matrices. The element $\small A=\left[ \begin{array}{ccc} 0 & 0 & a \\ 0 & 0 & 0 \\ 0 & 0 & 0 \end{array}\right]$ annihilates everything in the center of the triangular matrices. This is the reason that $f \star g$ is differentiable whereas $g \star f$ is not; $d_pf$ corresponds to $A$ whereas $d_pg$ does not annihilate the center of the algebra. 
} in $\MatA$. Recall Equation \ref{eqn:noncommuteprobII} showed we need $d_pf$ to annihilate the center of the algebra in order that $f \star g$ be $\Acal$-differentiable at $p$. Likewise, to have $g \star f$ differentiable over $\Acal$ at $p$ we need $d_pg$ to annihilate the center of $\Acal$. This is the distinction between $f$ and $g$ in this example, only $f$ has $d_pf$ annihilating the center of $\Acal$.
\end{exa}

\begin{thm} \label{thm:Achainrule}
Suppose $U,V \subseteq \Acal$ are open sets and $g: U \rightarrow V$ and $f: V \rightarrow \Acal$ are $\Acal$-differentiable functions. If $p \in U$ then
$$ (f \comp g)'(p) = f'(g(p)) \star g'(p). $$
\end{thm}

\noindent
{\bf Proof:} if $f \in \Fun (U)$ and $g \in \Fun (V)$ then $f$ and $g$ are $\RN$-differentiable on $U,V$ respective. Moreover, by the usual real calculus of a normed linear space, if the composite $f \comp g$ is defined at $p$ we have an elegant chain rule in terms of differentials:
$ d_p(f \comp g) = d_{g(p)}f \comp d_p g $. Let $v,w \in \Acal$ and consider
\begin{align} 
 d_p(f \comp g)(v \star w) &=  (d_{g(p)}f \comp d_p g) (v \star w)  &\text{: real chain rule} \\ \notag
 &= d_{g(p)}f(d_pg(v \star w))  &\text{: def. of composite} \\ \notag
 &= d_{g(p)}f(d_pg(v) \star w))  &\text{: as $g \in \Fun(p)$} \\ \notag
 &= d_{g(p)}f(d_pg(v)) \star w  &\text{: as $f \in \Fun(g(p))$} \\ \notag
&= d_p (f \comp g)(v) \star w  &\text{: real chain rule}
 \end{align}
 Thus $d_p(f \comp g) \in \EndA$ which shows $f \comp g$ is $\Acal$-differentiable at $p$. Moreover, as $f \in \Fun(g(p))$ implies $d_{g(p)}f(w) = f'(g(p)) \star w$ and $g \in \Fun(p)$ implies $d_pg(v) = g'(p) \star v$ we derive
\begin{equation}
 d_p(f \comp g)( v) = d_{g(p)}f(d_pg( v)) = f'(g(p)) \star d_pg(v) = f'(g(p)) \star g'(p) \star v 
\end{equation}
for each $v \in \Acal$. Therefore, $(f \comp g)'(p) = f'(g(p)) \star g'(p)$. $\Box$

\noindent
\begin{exa}
{\bf Claim:} Let $f(\zeta) = \zeta^n$ for some $n \in \NN$ then $f'(\zeta) = n \zeta^{n-1}$. \\ 

\noindent
We proceed by induction on $n$. If $n=1$ then $f(\zeta) = \zeta$ which is to say $f = Id$ and hence $d_pf = Id$ for each $p \in \Acal$. Moreover, 
\begin{equation}
d_pf(x \star y) = Id(x \star y) = x \star y = Id(x) \star y = d_p f(x) \star y
\end{equation}
which shows $f$ is $\Acal$-differentiable on $\Acal$. We calculate $\# (d_pf) = \#(Id) = Id( \mathds{1}) = \mathds{1} = 1 \zeta^{0}$ hence the claim is true for $n=1$. Suppose the claim holds for some $n \in \NN$. Define $f(\zeta) = \zeta^n$ and $g(\zeta) = \zeta$. We have $f'(\zeta) = n \zeta^{n-1}$ by the induction hypothesis and we already argued $g'(\zeta) = 1$. Thus, Theorem \ref{thm:Aproductrule} applies to calculate $\zeta^{n+1} = f(\zeta) \star g(\zeta)$:  
\begin{equation}
  (f \star g)'(\zeta)  = n \zeta^{n-1} \star \zeta + \zeta^{n} \star \mathds{1} = (n+1)\zeta^{(n+1)-1}
 \end{equation}
thus the claim is true for $n+1$ and we conclude $\frac{d}{d\zeta}(\zeta^n) = n \zeta^{n-1}$ for all $n \in \NN$.
\end{exa}

\noindent
Admittedly, we just introduced a new notation; if $f$ is an $\Acal$-differentiable function and $\zeta$ denotes an $\Acal$ variable then we write
\begin{equation}
 f'(\zeta) = \frac{df}{d\zeta}(\zeta) = \frac{d}{d\zeta}( f( \zeta)) \qquad \& \qquad f' = \frac{df}{d\zeta}
\end{equation}
If the $\Acal$-differentiability of $f: \Acal \rightarrow \Acal$ is not certain then we may still meaningfully calculate $\frac{\partial}{\partial \zeta}$ as a particular $\Acal$-linear combination of real partial derivatives. In other words, we are able to find an $\Acal$-generalization of Wirtinger's calculus. Details are given in the next section. \\

\noindent
If $\Acal \approx \Bcal$ then $\Acal$ and $\Bcal$ differentiable functions are related through the isomorphism. 

\begin{thm} \label{thm:difftransfer}
Let $\Psi: \Acal \rightarrow \Bcal$ be an isomorphism of unital, associative finite dimensional algebras over $\RN$. If $f$ is $\Acal$ differentiable at $p$ then $g = \Psi \comp f \comp \Psi^{-1}$ is $\Bcal$-differentiable at $\Psi(p)$. Moreover, $g'(p) = (\Psi \comp f' \comp \Psi^{-1})(p)$.
\end{thm}

\noindent
{\bf Proof:} Let $(\Acal, \star)$ and $(\Bcal, \ast)$ be finite dimensional isomorphic unital associative algebras via the isomorpism $\Psi: \Acal \rightarrow \Bcal$. In particular, $\Psi$ is a linear bijection and $\Psi( v \star w) = \Psi(v) \ast \Psi(w)$ for all $v,w \in \Acal$. Since $\Psi$ and $\Psi^{-1}$ are linear maps on normed linear spaces of finite dimension we know these are smooth real maps with $d_p \Psi = \Psi$ for each $p \in \Acal$ and $d_q\Psi^{-1} = \Psi^{-1}$ for each $q \in \Bcal$. If $f$ is $\Acal$ differentiable at $p$ then $d_pf \in \EndA$. Define $g = \Psi \comp f \comp \Psi^{-1}$ and notice $d_qg$ exists as $g$ is formed from the composite of differentiable maps. The chain rule\footnote{explicitly $ d_qg =  d_{f(\Psi^{-1}(q))}\Psi \comp d_{\Psi^{-1}(q)}f \comp d_q\Psi^{-1} $} provides,
\begin{equation}
 dg =  d\Psi \comp df \comp d\Psi^{-1}   \ \ \Rightarrow \ \ dg = \Psi \comp df \comp \Psi^{-1}
 \end{equation}
as $\Psi, \Psi^{-1}$ are linear maps. We seek to show $d_qg$ is right-$\Bcal$-linear at $q = \Psi(p)$. Calculate,
\begin{align} 
d_qg (v \ast w) & = \Psi ( d_pf ( \Psi^{-1}(v \ast w))) \\ \notag
& = \Psi ( d_pf ( \Psi^{-1}(v) \star \Psi^{-1}(w))) \\ \notag
& = \Psi ( d_pf ( \Psi^{-1}(v)) \star \Psi^{-1}(w)) \\ \notag
& = \Psi ( d_pf ( \Psi^{-1}(v))) \ast \Psi(\Psi^{-1}(w)) \\ \notag
& = d_qg(v) \ast w.
\end{align}
Thus $d_qg \in \EndB$ and we find $g$ is $\Bcal$-differentiable at $q = \Psi(p)$ as claimed. $\Box$ \\

\noindent
Isomorphic algebras have algebra-differentiable functions which naturally correspond:

\begin{coro}
If $\Psi: \Acal \rightarrow \Bcal$ be an isomorphism of unital, associative finite dimensional algebras over $\RN$ and $U \subseteq \Acal$ an open set then each $f \in \text{C}_{\Acal}(U)$ can be written as a composite $f = \Psi^{-1} \comp g \comp \Psi$ for some $g \in  \text{C}_{\Bcal}(\Psi(U))$
\end{coro}

\noindent
{\bf Proof:} Observe $g = \Psi \comp f \comp \Psi^{-1}$ satisfies $f = \Psi^{-1} \comp g \comp \Psi$. Moreover, $\Bcal$-differentiability of $g$ at $\Psi(p) \in \Psi(U)$ is naturally derived from the given $\Acal$-differentiability of $f$ at $\Psi(p) \in \Psi(U)$ via the result of Theorem \ref{thm:difftransfer}. $\Box$

\section{$\Acal$-Cauchy Riemann equations} \label{sec:CReqns}
Our first goal in this section is to describe a Wirtinger calculus for $\Acal$. In particular, we define the partial derivatives of an algebra variable and its conjugate variables. We should caution, the term {\it conjugate} ought not be taken too literally. These conjugates generally do not form automorphisms of the algebra. Their utility is made manifest that the play much the same role as the usual conjugate in complex analysis. We follow the construction of Alvarez-Parrilla, Fr{\' i}as-Armenta, L{\' o}pez-Gonz{\' a}lez and Yee-Romero directly for the definition below (this is Equation 4.3 of \cite{pagr2012}):

\begin{de} \label{defn:conjugates}
Suppose $\Acal$ has an invertible basis which begins with the multiplicative identity of the algebra. In particular, $\beta = \{ v_1, v_2, \dots , v_n \}$ is a basis for $\Acal$ with $v_1 = \mathds{1}$. If $\zeta = x_1v_1+x_2v_2+ \cdots + x_nv_n$ then we define the {\bf $j$-th conjugate of $\zeta$} as follows:
$$ \overline{\zeta}_j  = \zeta - 2x_j v_j  = x_1\mathds{1}+ \cdots + x_{j-1}v_{j-1}-x_jv_j+ x_{j+1}v_{j+1}+ \cdots + x_n v_n$$
for $j=2,3 \dots , n$.
\end{de}

\noindent
In some sense, the variables $\zeta, \overline{\zeta}_2, \dots , \overline{\zeta}_n$ are simply an algebra notation for $n$ real variables. Given a function of $x_1, \dots , x_n$ we are free to express the function in terms of the algebra variables $\zeta, \overline{\zeta}_2, \dots , \overline{\zeta}_n$. 

\begin{thm} \label{thm:conjugatevariableformulas}
Suppose $\Acal$ has invertible basis $\beta = \{ \mathds{1}, v_2, \dots , v_n \}$ and $\zeta = \sum_{i=1}^n x_i v_i$ and $\overline{\zeta}_j  = \zeta - 2x_j v_j$ for $j=2, \dots , n$ then using $\frac{1}{v_j}$ to denote $v_j^{-1}$ and omit $\mathds{1}$ we find:
\begin{enumerate}[{\bf (i.)}]
\item $ \ds x_j = \frac{1}{2v_j} \bigl( \zeta - \overline{\zeta}_j \bigr)$ for $j=2, \dots , n$.
\item $ \ds x_1 = \frac{1}{2}\biggl((3-n)\zeta + \sum_{j=2}^n \overline{\zeta}_j \biggr)$.
\end{enumerate}
\end{thm}

\noindent
{\bf Proof:} to obtain (i.) simply solve $\overline{\zeta}_j  = \zeta - 2x_j v_j$ for $x_j$. Then to derive (ii.) we solve $\zeta = x_1\mathds{1}+x_2v_2+ \cdots +x_nv_n$ for $x_1\mathds{1}$:
\begin{equation}
 x_1\mathds{1} = \zeta -\sum_{j=2}^n x_jv_j \\ 
                        = \zeta - \frac{1}{2}\sum_{j=2}^n (\zeta - \overline{\zeta}_j) 
                        = \frac{1}{2}\biggl[(3-n)\zeta+\sum_{j=2}^n \overline{\zeta}_j \biggr]
\end{equation}
finally, omit $\mathds{1}$ to obtain (ii.) $\Box$ \\

\noindent
It may be helpful to review the usual results of the Wirtinger's \cite{wirtinger} calculus for complex analysis. Consider $\Acal = \CN$ where $z = x+iy$ and $\bar{z} = x-iy$ hence $x = \frac{1}{2}(z+ \bar{z})$ and $y = \frac{i}{2}(\bar{z}-z)$. Hence, formally,
\begin{equation}
 \frac{\partial}{\partial z} = 
\frac{\partial x}{\partial z}\frac{\partial}{\partial x}
+  
\frac{\partial y}{\partial z}\frac{\partial}{\partial y} = \frac{1}{2} \left( \frac{\partial}{\partial x}
 - i\frac{\partial}{\partial y}
\right) \ \& \ \  
 \frac{\partial}{\partial \bar{z}} = 
\frac{\partial x}{\partial \bar{z}}\frac{\partial}{\partial x}
+  
\frac{\partial y}{\partial \bar{z}}\frac{\partial}{\partial y} = \frac{1}{2} \left( \frac{\partial}{\partial x}
 + i\frac{\partial}{\partial y}
\right).
\end{equation}
Furthermore, $\partial_x = \partial_z + \partial_{\bar{z}}$ and  $\partial_y = i(\partial_z - \partial_{\bar{z}})$. We seek similar formulas for $\Acal$. Hence, consider Theorem \ref{thm:conjugatevariableformulas} shows that if $\Acal$ has invertible basis $\beta = \{ \mathds{1}, v_2, \dots , v_n \}$ and $\zeta = \sum_{i=1}^n x_i v_i$ and $\overline{\zeta}_j  = \zeta - 2x_j v_j$ for $j=2, \dots , n$ then
\begin{equation}
  x_j = \frac{1}{2v_j} \bigl( \zeta - \overline{\zeta}_j \bigr) \ \ \ \ \& \ \ \ \
   x_1 = \frac{1}{2}\biggl((3-n)\zeta + \sum_{j=2}^n \overline{\zeta}_j \biggr). 
\end{equation}
Formally, $\frac{\partial x_j}{\partial \zeta} = \frac{1}{2v_j}$ for $j=2,\dots , n$ and $\frac{\partial x_1}{\partial \zeta} = \frac{3-n}{2}$ hence we suspect
\begin{equation}
 \frac{\partial}{\partial \zeta} = \sum_{j=1}^n \frac{\partial x_j}{\partial \zeta}\frac{\partial}{\partial x_j} = \frac{1}{2}\left((3-n)\frac{\partial}{\partial x_1}+ \frac{1}{v_2}\frac{\partial}{\partial x_2}+ \cdots + \frac{1}{ v_n}\frac{\partial}{\partial x_n} \right)
\end{equation}
whereas $\frac{\partial x_j}{\partial \overline{\zeta}_k} = \frac{-1}{2v_j}\delta_{jk}$ and $\frac{\partial x_1}{\partial \overline{\zeta}_k} = \frac{1}{2}$ thus we speculate:
\begin{equation}
 \frac{\partial }{\partial \overline{\zeta}_k} = 
\sum_{j=1}^n \frac{\partial x_j}{\partial \overline{\zeta}_k}\frac{\partial}{\partial x_j} = \frac{1}{2}\left(\frac{\partial}{\partial x_1}-\frac{1}{v_k} \frac{\partial}{\partial x_k}\right). 
\end{equation}
The calculations above convinced us to make the definition below (which is slightly different than the defintion offered in \cite{pagr2012} where $\partial/ \partial \zeta$ is defined differently in their Equation 4.7). 

\begin{de} \label{defn:partialAderivatives}
Suppose $f: \Acal \rightarrow \Acal$ is $\RN$-differentiable. Furthermore, suppose $\beta = \{ \mathds{1},v_2, \dots , v_n \}$ is an invertible basis and $\zeta = x_1 \mathds{1}+x_2v_2 + \cdots + x_nv_n$. We define
$$ 
\frac{\partial }{\partial \zeta} = \frac{1}{2}\left((3-n)\frac{\partial}{\partial x_1}+ \frac{1}{v_2}\frac{\partial}{\partial x_2}+ \cdots + \frac{1}{ v_n}\frac{\partial}{\partial x_n} \right)
 \qquad \& \qquad
\frac{\partial }{\partial \overline{\zeta}_k} = \frac{1}{2} \left( \frac{\partial}{\partial x_1} - \frac{1}{v_k}\frac{\partial }{\partial x_k}\right)$$
for $j=2,3,\dots , n$.
\end{de}

\noindent
The merit of Definition \ref{defn:partialAderivatives} is seen in the theorem below\footnote{compare with Equation 4.9 in \cite{pagr2012}}:

\begin{thm} \label{thm:wirtingerworks}
Given the notation of Definition \ref{defn:partialAderivatives},
$$ \frac{\partial \zeta}{\partial \zeta}= 1, \ \ \ \ \frac{\partial \overline{\zeta}_j}{\partial \zeta} = 0, \ \ \ \
\frac{\partial \overline{\zeta}_j}{\partial \overline{\zeta}_j}= 1, \ \ \ \
\frac{\partial \overline{\zeta}_j}{\partial \overline{\zeta}_k}= 0, \ \ \ \ 
\frac{\partial \zeta}{\partial \overline{\zeta}_j}= 0  $$
for all $j=2, \dots , n$ and $k \neq j$.
\end{thm}

\noindent
{\bf Proof:} simple calculation. Consider:
\begin{align} 
 \frac{\partial \zeta}{\partial \zeta} 
&= \frac{1}{2}\left((3-n)\frac{\partial}{\partial x_1}+ \frac{1}{v_2}\frac{\partial}{\partial x_2}+ \cdots + \frac{1}{ v_n}\frac{\partial}{\partial x_n} \right)\left( x_1+ \cdots + x_nv_n \right) \\ \notag
 &= \frac{(3-n)+n-1}{2} \\ \notag
 &= 1.
\end{align}
For $j=2,\dots , n$ we calculate:
\begin{align} 
 \frac{\partial \overline{\zeta}_j}{\partial \zeta} 
&= \frac{1}{2}\left((3-n)\frac{\partial}{\partial x_1}+ \frac{1}{v_2}\frac{\partial}{\partial x_2}+ \cdots + \frac{1}{ v_n}\frac{\partial}{\partial x_n} \right)\left( x_1+ \cdots  -x_jv_j+ \cdots + x_nv_n \right) \\ \notag
&= \frac{3-n+n-3}{2} \\ \notag
&= 0,
\end{align}
and
\begin{equation}
 \frac{\partial \overline{\zeta}_j}{\partial \overline{\zeta}_j} 
= \frac{1}{2}\left(\frac{\partial}{\partial x_1}- \frac{1}{v_j}\frac{\partial}{\partial x_j}\right)\left( x_1+ \cdots  -x_jv_j+ \cdots + x_nv_n \right) 
= \frac{1}{2}\left( 1+ \frac{1}{v_j}v_j \right) = 1, 
\end{equation}
and for $k \neq j$,
\begin{equation}
 \frac{\partial \overline{\zeta}_j}{\partial \overline{\zeta}_k} 
= \frac{1}{2}\left(\frac{\partial}{\partial x_1}- \frac{1}{v_k}\frac{\partial}{\partial x_k}\right)\left( x_1+ \cdots  -x_jv_j+ \cdots + x_nv_n \right) 
= \frac{1}{2}\left( 1 - \frac{1}{v_k}v_k \right) = 0 . 
\end{equation}
and finally
\begin{equation}
 \frac{\partial \zeta}{\partial \overline{\zeta}_j} 
= \frac{1}{2}\left(\frac{\partial}{\partial x_1}- \frac{1}{v_j}\frac{\partial}{\partial x_j}\right)\left( x_1+ \cdots  +x_jv_j+ \cdots + x_nv_n \right)
= \frac{1}{2}\left( 1 - \frac{1}{v_j}v_j \right) = 0 .  
\end{equation}
In summary, the derivatives above show we may think of $\zeta$ and $\overline{\zeta}_j$ as independent variables. $\Box$ \\

\noindent
We should connect the formal derivatives above with $\Acal$-differentiability of a function. In the context of complex analysis, you may recall that complex differentiability is also characterized by the equation $\frac{\partial f}{\partial \bar{z}}=0$. For $\Acal$, we expect additional conditions which make use of the $(n-1)$-conjugates. Notice $\frac{\partial f}{\partial \overline{\zeta}_j}=0$ for $j=2, \dots , n$ is a likely generalization. Count $n-1$ algebra equations which decomposes into $n^2-n$ real equations which are necessarily solved by the components of an $\Acal$-differentiable function.

\begin{thm} \label{thm:CReqnsI}
Suppose $\Acal$ is an associative unital algebra of finite dimension over $\RN$ with basis $\beta = \{ v_1, \dots , v_n \}$. If $f: \Acal \rightarrow \Acal$ is $\Acal$-differentiable at $p$ then 
\begin{enumerate}[{\bf (i.)}]
\item given $C_{ij}^k \in \RN$ for which $v_i \star v_j = \sum_k C_{ij}^k v_k$ we have $\ds \frac{\partial f}{\partial x_i} \star v_j =  \sum_kC_{ij}^k \frac{\partial f}{\partial x_k}$ 
\item given basis $\beta$ has $v_1 =\mathds{1}$ we find $\ds \frac{\partial f}{\partial x_j} =  \frac{\partial f}{\partial x_1} \star v_j $ for $j=2,\dots , n$.
\item if $\Acal$ is commutative then $\frac{\partial f}{\partial x_i} \star v_j = v_i \star \frac{\partial f}{\partial x_j} $ for all $i,j = 1, 2, \dots , n$.
\end{enumerate}
\end{thm}

\noindent
{\bf Proof:} suppose $f$ is $\Acal$-differentiable at $p$ then $d_pf \in \EndA$.  Partial derivatives with respect to the basis $\beta$ in $\Acal$ were given by Equation \ref{eqn:partialderNLS}; $d_pf(v_k) = \frac{\partial f}{\partial x_k}$. To derive (i.) suppose $v_i \star v_j = \sum_k C_{ij}^k v_k$ and calculate:
\begin{equation}
 \frac{\partial f}{\partial x_i} \star v_j = d_pf(v_i) \star v_j = d_pf (v_i \star v_j) = d_pf \left(  \sum_kC_{ij}^k v_k \right) = \sum_kC_{ij}^kd_pf(v_k) =  \sum_kC_{ij}^k \frac{\partial f}{\partial x_k}.
\end{equation}
Likewise, (ii.) follows as 
\begin{equation}
 d_pf(v_j) = d_pf(\mathds{1} \star v_j) = d_pf( \mathds{1}) \star v_j \ \ \Rightarrow \ \ \frac{\partial f}{\partial x_j} = \frac{\partial f}{\partial x_1} \star v_j 
\end{equation} 
for $j=2, \dots , n$. Finally, in the case $\Acal$ is commutative we derive (iii.) as follows:
\begin{equation}
 \frac{\partial f}{\partial x_i} \star v_j = d_pf(v_i) \star v_j = d_pf(v_i \star v_j) = d_pf(v_j \star v_i) = d_pf(v_j) \star v_i = \frac{\partial f}{\partial x_j} \star v_i.  
 \end{equation}
Thus, once more using commutativity of $\Acal$, $ \frac{\partial f}{\partial x_i} \star v_j = v_i \star \frac{\partial f}{\partial x_j}$.  $\Box$ \\

\noindent
The (i.), (ii.) or (iii.) equations above are known as {\it generalized Cauchy Riemann Equations} by many authors. We prefer to call them the $\Acal$-CR-equations as to be specific. Both (i.) and (ii.) are suitable sets of equations when $\Acal$ is noncommutative. If $\mathds{1}$ is conveniently presented in a basis for $\Acal$ then (ii.) is the convenient description of $\Acal$-differentiable functions.

\begin{thm} \label{thm:CReqnsII}
Let $\beta = \{ \mathds{1}, v_2, \dots , v_n \}$ be an invertible basis for the commutative algebra $\Acal$. If $f: \Acal \rightarrow \Acal$ is $\Acal$-differentiable at $p$ then $\ds \frac{\partial f}{\partial \overline{\zeta}_j} = 0$ for $j=2, \dots , n$.
\end{thm}

\noindent
{\bf Proof:} following Definition \ref{defn:partialAderivatives}
\begin{equation} 
\frac{\partial f}{\partial \overline{\zeta}_k} = 
\frac{1}{2} \left( \frac{\partial f}{\partial x_1} - \frac{1}{v_k}\star\frac{\partial f}{\partial x_k}\right) 
=\frac{1}{2} \left( \frac{\partial f}{\partial x_1} - \frac{1}{v_k}\star \frac{\partial f}{\partial x_1} \star v_k \right) 
= 0 
\end{equation} 
as $\Acal$ is assumed commutative and $\frac{1}{v_k} \star v_k= 1$.  $\Box$ \\

\noindent
The usual additive and product rules hold for $\partial/ \partial \zeta$ and $\partial / \partial \overline{\zeta}_j$. 

\begin{thm} \label{thm:rulesofcalculus}
Using the notation of Definition \ref{defn:partialAderivatives}, if $f, g: \Acal \rightarrow \Acal$ are differentiable then
$$ \frac{\partial}{\partial \zeta}(f+g) = \frac{\partial f}{\partial \zeta} + \frac{\partial g}{\partial \zeta} \qquad \& \qquad \frac{\partial}{\partial \overline{\zeta}_j}(f+g) = \frac{\partial f}{\partial \zeta} + \frac{\partial g}{\partial \overline{\zeta}_j} $$
for $j=2,\dots , n$. Likewise,
$$ \frac{\partial}{\partial \zeta}(f \star g) = \frac{\partial f}{\partial \zeta} \star g+f \star \frac{\partial g}{\partial \zeta} \qquad \& \qquad \frac{\partial}{\partial \overline{\zeta}_j}(f \star g) = \frac{\partial f}{\partial \overline{\zeta}_j}\star g+f \star \frac{\partial g}{\partial \overline{\zeta}_j} $$ 
\end{thm}

\noindent
{\bf Proof:} suppose $f$ and $g$ are real differentiable functions on $\Acal$, 
\begin{equation}
 \frac{1}{2} \left( \frac{\partial}{\partial x_1} - \frac{1}{v_k}\frac{\partial }{\partial x_k}\right)(f+g) = \frac{1}{2} \left( \frac{\partial f}{\partial x_1} - \frac{1}{v_k}\frac{\partial f }{\partial x_k}\right)
+ \frac{1}{2} \left( \frac{\partial g}{\partial x_1} - \frac{1}{v_k}\frac{\partial g }{\partial x_k}\right) 
\end{equation}
thus $\frac{\partial}{\partial \overline{\zeta}_j}(f+g) = \frac{\partial f}{\partial \zeta} + \frac{\partial g}{\partial  \overline{\zeta}_j}$. Using the structure constants $C_{ijk}$ we express $f \star g = \sum_{ijk} C_{ijk}f_ig_jv_k$ and it follows $\partial_j(f \star g) = \partial_j f \star g+f \star \partial_jg$. Consequently,
\begin{align} 
 \frac{1}{2} \left( \frac{\partial}{\partial x_1} - \frac{1}{v_k}\frac{\partial }{\partial x_k}\right)(f \star g) &= \frac{1}{2} \left( \frac{\partial f}{\partial x_1} \star g + f \star \frac{\partial g}{\partial x_1}- \frac{1}{v_k}\left[\frac{\partial f }{\partial x_k} \star g+f \star \frac{\partial g }{\partial x_k} \right]\right) \\ \notag
 &= \frac{1}{2} \left( \frac{\partial f}{\partial x_1}  + - \frac{1}{v_k}\frac{\partial f }{\partial x_k} \right) \star g + f \star \frac{1}{2} \left( \frac{\partial g}{\partial x_1}  + - \frac{1}{v_k}\frac{\partial g }{\partial x_k} \right).
\end{align}
Hence $\frac{\partial}{\partial \overline{\zeta}_j}(f \star g) = \frac{\partial f}{\partial \overline{\zeta}_j}\star g+f \star \frac{\partial g}{\partial \overline{\zeta}_j}$. The identities for $\partial/ \partial \zeta$ follow from similar calculations.
 $\Box$ \\

\noindent
The following example can be constructed in nearly every $\Acal$.

\begin{exa}
Suppose $\text{dim}(\Acal) \geq 2$. Let $f( \zeta) = \zeta \, \overline{\zeta}_2$ where $f: \Acal \rightarrow \Acal$ then 
$$\frac{\partial f}{\partial \zeta} = \frac{\partial \zeta}{\partial \zeta}\overline{\zeta}_2 + \zeta \frac{\partial\overline{\zeta}_2}{\partial \zeta} = \overline{\zeta}_2 \qquad \& \qquad 
\frac{\partial f}{\partial \overline{\zeta}_2} = \frac{\partial \zeta}{\partial \overline{\zeta}_2}\overline{\zeta}_2 + \zeta \frac{\partial\overline{\zeta}_2}{\partial \overline{\zeta}_2} = \zeta$$ 
This function is only $\Acal$-differentiable at the origin. In the usual complex analysis it is simply the square of the modulus; $f(z) = z\overline{z} = x^2+y^2$ where $z=x+iy$ has $\overline{z}_2 = x-iy$.
\end{exa}

\noindent
Inverting Definition \ref{defn:partialAderivatives} for $\partial/ \partial x_1, \dots, \partial / \partial x_n$ in terms of $\partial/ \partial \zeta, \dots , \partial / \partial \overline{\zeta}_n$ yields:

\begin{thm} \label{thm:inverserelations}
Using the notation of Definition \ref{defn:partialAderivatives}, 
$$ \frac{\partial}{\partial x_1} = \frac{\partial }{\partial \zeta} + \frac{\partial }{\overline{\zeta}_2} + \cdots + \frac{\partial }{\overline{\zeta}_n} \qquad \& \qquad  
\frac{\partial}{\partial x_k} = v_k \left(\frac{\partial }{\partial \zeta} + \frac{\partial }{\overline{\zeta}_2} + \cdots + \frac{\partial }{\overline{\zeta}_n}- 2\frac{\partial }{\overline{\zeta}_k} \right)
$$ 
\end{thm}

\noindent
{\bf Proof:} Begin with Definition \ref{defn:partialAderivatives} and note the identity for $\partial/\partial x_1$ follows immediately from summing $\partial/ \partial \zeta$ with the $n-1$ conjugate derivatives $\frac{\partial }{\partial \overline{\zeta}_2}, \dots , \frac{\partial }{\partial \overline{\zeta}_n}$:
\begin{equation} \label{eqn:derivative1}
 \frac{\partial }{\partial \zeta} + \frac{\partial }{\overline{\zeta}_2} + \cdots + \frac{\partial }{\overline{\zeta}_n} = \frac{3-n}{2}\frac{\partial}{\partial x_1}+ \frac{n-1}{2}\frac{\partial}{\partial x_1} = \frac{\partial}{\partial x_1}. 
\end{equation}
Following Definition \ref{defn:partialAderivatives} we substitute Equation \ref{eqn:derivative1} into the definition of  $\frac{\partial }{\partial \overline{\zeta}_k}$  to obtain:
\begin{equation}
 \frac{\partial }{\partial \overline{\zeta}_k} = \frac{1}{2} \left( \frac{\partial }{\partial \zeta} + \frac{\partial }{\overline{\zeta}_2} + \cdots + \frac{\partial }{\overline{\zeta}_n} - \frac{1}{v_k}\frac{\partial }{\partial x_k}\right) 
\end{equation}
It is now clear we can solve for $\frac{\partial}{\partial x_k}$ to obtain the desired result. $\Box$ \\ 

\noindent
In principle we can take a given PDE in $x_1, \dots, x_n$ and convert it to an $\Acal$ differential equation in $\zeta, \overline{\zeta}_1, \dots , \overline{\zeta}_n$. If we assume a solution for which all the conjugate derivative vanish then the PDE simplifies to an ordinary $\Acal$-differential equation. For some PDEs the corresponding $\Acal$-ODE may be solvable using elementary calculus. See Example \ref{exa:AODEbackwards} for a demonstration. 

\begin{exa} \label{exa:trihyperbolicidentities}
Consider $\Acal = \RN \oplus j \RN \oplus j^2\RN$ where $j^3=1$. We consider the algebra variable $\zeta = x +jy + z j^2$ and conjugate variables 
\begin{equation} 
\overline{\zeta}_2 = x-jy+j^2z \qquad \& \qquad \overline{\zeta}_3 = x+jy-j^2z 
\end{equation} 
In our current notation $\{ 1, j , j^2 \}$ forms an invertible basis with $v_2=j$ and $v_3=j^2$. Note $1/v_2 = j^2$ and $1/v_3 = j$.
It follows we have derivatives
\begin{equation}
\frac{\partial}{\partial \zeta} = \frac{1}{2} \left[j \frac{\partial}{\partial y} + j^2 \frac{\partial}{\partial z} \right] 
\ \ \& \ \ 
\frac{\partial}{\partial \overline{\zeta}_2} = \frac{1}{2} \left[ \frac{\partial}{\partial x} - j^2 \frac{\partial}{\partial y} \right]
\ \ \& \ \
\frac{\partial}{\partial \overline{\zeta}_3} = \frac{1}{2}\left[ \frac{\partial}{\partial x} - j \frac{\partial}{\partial z}\right].
\end{equation}
Thus, by Theorem \ref{thm:inverserelations}, or direct calculation, we find:
\begin{equation}
 \frac{\partial}{\partial x} = \frac{\partial }{\partial \zeta} + \frac{\partial }{\partial\overline{\zeta}_2} + \frac{\partial }{\partial\overline{\zeta}_3}, \ \  
\frac{\partial}{\partial y} = j \left(\frac{\partial }{\partial \zeta} - \frac{\partial }{\partial\overline{\zeta}_2} + \frac{\partial }{\partial\overline{\zeta}_3} \right), \ \   
\frac{\partial}{\partial z} = j^2 \left(\frac{\partial }{\partial \zeta} + \frac{\partial }{\partial\overline{\zeta}_2} - \frac{\partial }{\partial\overline{\zeta}_3} \right)
\end{equation}
From the formulas above we can derive the following differential identities:
\begin{align} 
\frac{\partial^2}{\partial x^2} - 
\frac{\partial}{\partial y} \frac{\partial}{\partial z} &= 2\left(\frac{\partial }{\partial \zeta}\frac{\partial }{\partial\overline{\zeta}_2} + \frac{\partial }{\partial \zeta}\frac{\partial }{\partial\overline{\zeta}_3} + \frac{\partial^2 }{\partial\overline{\zeta}_2^2} \right)\\ \notag 
\frac{\partial^2}{\partial y^2} - 
\frac{\partial}{\partial z} \frac{\partial}{\partial x} &= 2j^2\left(-2\frac{\partial }{\partial \zeta}\frac{\partial }{\partial\overline{\zeta}_2} + \frac{\partial }{\partial \zeta}\frac{\partial }{\partial\overline{\zeta}_3} -\frac{\partial }{\partial\overline{\zeta}_2}\frac{\partial }{\partial\overline{\zeta}_3} + \frac{\partial^2 }{\partial\overline{\zeta}_3^2} \right)\\ \notag 
\frac{\partial^2}{\partial z^2} - 
\frac{\partial}{\partial x} \frac{\partial}{\partial y} &= 2j\left(-2\frac{\partial }{\partial \zeta}\frac{\partial }{\partial\overline{\zeta}_3} + \frac{\partial }{\partial \zeta}\frac{\partial }{\partial\overline{\zeta}_2} -\frac{\partial }{\partial\overline{\zeta}_2}\frac{\partial }{\partial\overline{\zeta}_3} + \frac{\partial^2 }{\partial\overline{\zeta}_2^2} \right)
\end{align}
If $f = u+vj+j^2w$ is an $\Acal$-differentiable function then $\frac{\partial f }{\partial\overline{\zeta}_2} = 0$ and $\frac{\partial f }{\partial\overline{\zeta}_3}=0$. Therefore, $f$ is annihilated by the operators $\partial_x^2- \partial_y\partial_z$, $\partial_y^2- \partial_z\partial_x$ and $\partial_z^2- \partial_x\partial_y$. It follows that the component functions of $f$ must solve the corresponding PDEs:
$$ \Phi_{xx}-\Phi_{yz}=0, \ \ \Phi_{yy}-\Phi_{zx}=0, \ \ 
\Phi_{zz}-\Phi_{xy}=0. $$
These are known as the {\bf generalized Laplace Equations} for the $3$-hyperbolic numbers.
\end{exa}

\noindent
Generalized $\Acal$-Laplace Equations are differential consequences of the $\Acal$-CR equations. When system of PDEs happens to be the $\Acal$-Laplace equations we find any $\Acal$-differentiable function provides solutions to system. One may also wonder when a given system is consistent with the $\Acal$-Laplace equations. In the event a given system of PDEs was consistent then we may impose the $\Acal$-CR equations and their differential consequences on the given system of PDEs as to find a special subclass of $\Acal$-differentiable solutions. Computationally this section provides a roadmap for this procedure:
\begin{quote}
\begin{enumerate}[{\bf (1.)}]
\item given a PDE in real independent variables $x_1,x_2, \dots, x_n$ choose an algebra $\Acal$ of dimension $n$ to study in conjunction with the system.
\item convert the derivatives in the PDE with respect to $x_1, x_2, \dots , x_n$ to derivatives with respect to the algebra variables $\zeta, \overline{\zeta}_2, \dots , \overline{\zeta}_n$
\item impose that the derivatives with respect to $\overline{\zeta}_2, \dots , \overline{\zeta}_n$ vanish, study the resulting ordinary differential equation in $\zeta$. If possible, solve the $\Acal$-ODE which results. 
\end{enumerate}
\end{quote}

\noindent
The possibility that the technique above may produce novel solutions to particular systems of PDEs is one of the major motivations of this work.


\section{Deleted difference quotients} \label{sec:deleteddiffquotients}
There are several popular definitions of differentiability with respect to an algebra. Either we can follow the path of first semester calculus and use a difference quotient\footnote{suitably modified to avoid zero-divisors} or we can follow something involving a Frechet quotient\footnote{where $\Acal$-differentiability is imposed by an algebraic condition on the differential}. To see a rather detailed exposition of how these are related in the particular context of bicomplex or multi-complex numbers see \cite{price}.   \\

\noindent
The general concept this section is an adaptation and generalization of the arguments given in \cite{gadeaD1vsD2} for the context of the hyperbolic numbers. We show how some introductory results in \cite{gadeaD1vsD2} generalize to any commutative semisimple algebra of finite dimension over $\RN$. Ultimately the section demonstrates why we prefer  the definition of $\Acal$-differentiability given in Definition \ref{defn:Adiff} as opposed to the deleted-difference quotient definition. It is helpful to have a precise and abbreviated terminology for the discussion which follows:

\begin{de} \label{defn:terms}
Let $f: \text{dom}(f) \rightarrow \Acal$ be a function where $\text{dom}(f)$ is open and $p \in \text{dom}(f)$.
\begin{enumerate}[{\bf(1.)}]
\item If $f$ is $\Acal$-differentiable at $p$ then $f$ is $D_1$ at $p$. . 
\item If $\ds \lim_{ \Acalx \ni \zeta \rightarrow p} \frac{f(\zeta)-f(p)}{\zeta-p}$ exists then $f$ is $D_2$ at $p$.  
\end{enumerate}
If $f$ is $D_1 (D_2)$ for each $p \in  U$ then $f$ is $D_1 (D_2)$ on $U$. 
\end{de}

\noindent
If we fix our attention to a point then the class of $D_1$ and $D_2$ functions at $p$ are inequivalent. 

\begin{exa} \label{exa:dirichletmod}
In the spirit of Dirichlet we define $f(z) = \begin{cases} 0 & \text{if} \ \zeta \in \Acalx \cup \{ 0 \}\\ 1 & \text{if} \ \zeta \in \Acalzd - \{ 0 \} \end{cases}$. Since $f(\zeta) = 0$ for all $\zeta \in \Acalx$ we find $\frac{f (\zeta) - f(0)}{ \zeta} = 0 $ for all $\zeta \in \Acalx$. Thus $f$ is $D_2$-differentiable over $\Acal$ at $p=0$. Notice, $\Acal$-differentiability implies real differentiability and thus continuity. Clearly $f$ is not continuous at $p=0$ thus $f$ is not $D_1$ differentiable at $p=0$. 
\end{exa}

\noindent
In fact, $D_1$ on $\Acal$ at a point need not imply $D_2$ at a  the given point in $\Acal$.  For an explicit demonstration of this in the case of $\Acal = \calH$ see Example 2.2 part (2) of \cite{gadeaD1vsD2}.

\begin{thm} \label{thm:D2impliesD1prime}
Let $f$ be a function on $\Acal$ with $\zeta_o \in \text{dom}(f)$ where $\text{dom}(f)$ is open in $\Acal$. If $f$ is $D_2$ differentiable at $\zeta_o$ then 
\begin{equation} \label{eqn:D1prime} 
\lim_{ \overset{\zeta \rightarrow \zeta_o}{\zeta-\zeta_o \, \in \, \Acalx}} \frac{||f(\zeta)-f(\zeta_o)- \lambda \star (\zeta -\zeta_o)||}{||\zeta -\zeta_o||} = 0 
\end{equation}
\end{thm} 

\noindent
{\bf Proof:} since $f$ is $D_2$ differentiable at $\zeta_o$ there exists $\lambda \in \Acal$ to which the deleted-difference quotient of $f$ converges at $\zeta_o$. In particular, for each $\epsilon >0$ there exists $\delta>0$ for which $\zeta -\zeta_o \in \Acalx$ with $0< || \zeta - \zeta_o || < \delta$ implies
\begin{equation}
 \bigg{|}\bigg{|}\frac{f(\zeta)-f(\zeta_o)}{\zeta -\zeta_o} - \lambda \bigg{|}\bigg{|} < \epsilon \ \ \Rightarrow \ \ \bigg{|}\bigg{|}\frac{f(\zeta)-f(\zeta_o)- \lambda \star (\zeta -\zeta_o)}{\zeta -\zeta_o}  \bigg{|}\bigg{|} < \epsilon 
\end{equation}
Apply Corollary \ref{thm:quotientinequality}, 
\begin{equation}
  \frac{||f(\zeta)-f(\zeta_o)- \lambda \star (\zeta -\zeta_o)||}{||\zeta -\zeta_o||} \leq \Abound \bigg{|}\bigg{|}\frac{f(\zeta)-f(\zeta_o)- \lambda \star (\zeta -\zeta_o)}{\zeta -\zeta_o}  \bigg{|}\bigg{|} < \Abound \epsilon
 \end{equation}
Thus $\ds  \frac{||f(\zeta)-f(\zeta_o)- \lambda \star (\zeta -\zeta_o)||}{||\zeta -\zeta_o||} \rightarrow 0$ as $\zeta \rightarrow \zeta_o$ for $\zeta -\zeta_o \in \Acalx$. $\Box$ \\

\noindent
The deleted limit in Equation \ref{eqn:D1prime} almost provides $D_1$-differentiability at $\zeta_o$. To overcome the difficulty of Example \ref{exa:dirichletmod} it suffices to assume continuity of $f$ near $\zeta_o$.

\begin{thm} \label{thm:D2andcontinuousimpliesD1atpoint}
Let $f$ be a function on $\Acal$ which is continuous in some open set containing $\zeta_o$. If $f$ is $D_2$ differentiable at $\zeta_o$ then $f$ is $D_1$ differentiable at $\zeta_o$.
\end{thm}

\noindent
{\bf Proof:} let $\epsilon >0$ and use Theorem \ref{thm:D2impliesD1prime} to choose $\delta >0$ such that $\zeta- \zeta_o \in \Acalx$ and $||\zeta- \zeta_o|| < \delta$ implies $ ||f(\zeta)-f(\zeta_o)- \lambda \star (\zeta -\zeta_o)|| < \epsilon ||\zeta -\zeta_o||$. It remains to show $||f(\zeta)-f(\zeta_o)- \lambda \star (\zeta -\zeta_o)|| < \epsilon ||\zeta -\zeta_o||$ for $\zeta- \zeta_o \notin \Acalx$. We begin by making $\delta$ smaller (if necessary) such that $f$ is continuous on $U=\{ \zeta \in \text{dom}(f) \ | \ ||\zeta- \zeta_o|| < \delta \}$. Theorem \ref{thm:unitsdense} implies $U \cap \Acalx$ is dense in $U$. Consequently, if $\zeta_1- \zeta_o \in \Acalzd \cap U$ then $\zeta_1 - \zeta_o$ is a limit point of $U \cap \Acalx$.  Hence there exists a sequence of points $\zeta_n- \zeta_o \in U \cap \Acalx$ for which $\zeta_n- \zeta_o \rightarrow \zeta_1- \zeta_o$. Hence, $\zeta_n \rightarrow \zeta_1$ and by continuity of $f$ near $\zeta_o$ we find $f(\zeta_n) \rightarrow f(\zeta_1)$. Observe, as $\zeta_n - \zeta_o \in U \cap \Acalx$ we have the estimate $ ||f(\zeta_n)-f(\zeta_o)- \lambda \star (\zeta_n -\zeta_o)|| < \epsilon ||\zeta_n -\zeta_o||$. Hence, as $n \rightarrow \infty$ we find $||f(\zeta_1)-f(\zeta_o)- \lambda \star (\zeta_1 -\zeta_o)|| < \epsilon ||\zeta_1 -\zeta_o||$. But, as $\zeta_1$ was an arbitrary zero-divisor near $\zeta_o$ we find $ \lim_{ \zeta \rightarrow \zeta_o} \frac{||f(\zeta)-f(\zeta_o)- \lambda \star (\zeta -\zeta_o)||}{||\zeta -\zeta_o||} = 0 $. Thus the Frechet derivative of $f$ at $\zeta_o$ exists and the differential $d_{\zeta_o}f \in \EndA$ since $d_{\zeta_o}f(h) = \lambda \star h$. We conclude $f$ is $D_1$ at $\zeta_o$. $\Box$ \\

\begin{thm} \label{thm:D2holomorphicisalsoD1WHAT}
Let $U \subset \Acal$ be open. If $f$ is $D_2$ at each point in $U$ then $f$ is continuous on $U$
\end{thm}

\noindent
{\bf Proof:} Suppose $f$ is $D_2$ at each point of the open set $U$. Let $\zeta_o \in U$. By Theorem \ref{thm:D2impliesD1prime}
\begin{equation} \label{eqn:mostlycontinuous}
\lim_{ \overset{\zeta \rightarrow \zeta_o}{\zeta-\zeta_o \, \in \, \Acalx}} \frac{||f(\zeta)-f(\zeta_o)- \lambda \star (\zeta -\zeta_o)||}{||\zeta -\zeta_o||} = 0 
\end{equation}
But, $ \lim_{ \overset{\zeta \rightarrow \zeta_o}{\zeta-\zeta_o \, \in \, \Acalx}}  ||\zeta -\zeta_o|| = 0$ and  $ \lim_{ \overset{\zeta \rightarrow \zeta_o}{\zeta-\zeta_o \, \in \, \Acalx}}  ||\lambda \star (\zeta -\zeta_o)|| = 0$ hence we deduce
\begin{equation} \label{eqn:almostzero}
 \lim_{ \overset{\zeta \rightarrow \zeta_o}{\zeta-\zeta_o \, \in \, \Acalx}}||f(\zeta)-f(\zeta_o)|| = 0.
\end{equation}
It remains to show $\zeta$ for which $\zeta \rightarrow \zeta_o$ with $\zeta- \zeta_o \in \Acalzd$ also have $|| f(z) - f(\zeta_o) || \rightarrow 0$. Let $\eps >0$ and choose $\delta>0$ with $\{ \zeta \ | \ || \zeta - \zeta_o || \} \subset U$ and for which $\zeta- \zeta_o \in \Acalx $ and 
$||\zeta - \zeta_o|| < \delta$ imply $|| f( \zeta) - f(\zeta_o) || < \eps/2$. Selection of such $\delta>0$ is possible by Equation \ref{eqn:almostzero}. Form a triangle with vertices $\zeta_o, \zeta_2, \zeta_1$ where $|| \zeta_2-\zeta_o||< || \zeta_1- \zeta_o|| < \delta$. By construction $\zeta_1 \in U$ thus $f$ is $D_2$ at $\zeta_1$. Hence, following the thought behind Equation \ref{eqn:almostzero} once more, we find there exists $\delta'>0$ for which $|| \zeta_2 - \zeta_1 || < \delta'$ and $\zeta_2 - \zeta_1 \in \Acalx$ imply $|| f(\zeta_2) - f(\zeta_1) || < \eps/2$. Since $\Acalx$ is dense in $\Acal$ we are free to move $\zeta_2$ as close as we wish to $\zeta_1$ while maintaining $\zeta_2 - \zeta_1 \in \Acalx$ and $\zeta_2 - \zeta_o \in \Acalx$. Hence, for $\zeta_1$ such that $\zeta_1- \zeta_o \in \Acalzd$ with $|| \zeta_1- \zeta_o ||< \delta$ we find
\begin{equation}
 ||f(\zeta_1)-f(\zeta_o)|| \leq ||f(\zeta_1)-f(\zeta_2)||+||f(\zeta_2)-f(\zeta_o)|| < \epsilon/2+ \epsilon/2 = \epsilon.
\end{equation}
Therefore, $f$ is continuous at $\zeta_o$ and hence $f$ is continuous on $U$. $\Box$ \\

\noindent
Given Theorem \ref{thm:D2holomorphicisalsoD1WHAT} and Theorem \ref{thm:D2andcontinuousimpliesD1atpoint} we obtain the main result of this section:

\begin{thm} \label{thm:D2holomorphicisalsoD1}
Let $U$ be an open set in $\Acal$. If $f$ is $D_2$ at each point in $U$ then $f$ is $D_1$ on $U$.
\end{thm}
 
\noindent
In other words, functions which are $D_2$-holomorphic are necessarily $D_1$-holomorphic. We will see the converse need not be true. There are functions which are $D_1$ on an open set and yet fail to be $D_2$ at even a single point in the set. 

\begin{exa} \label{exa:nowhereD1}
Let $\Acal = \RN \oplus \epsilon \RN$ where $\epsilon^2=0$. In this algebra, $(a+b\epsilon) \epsilon = a\epsilon$ hence a typical matrix in $\MatA$ has the form $\left[ \begin{array}{cc} a & 0 \\ b & a \end{array} \right]$. If $f=u+\epsilon v: \Acal \rightarrow \Acal$ is $\Acal$-differentiable in the $D_1$ sense then $u_x = v_y$ and $u_y=0$. Conversely, if $u,v$ are continuously differentiable on $\Acal$ and satisfy $u_x = v_y$ and $u_y=0$ then $f=u+\epsilon v$ is $\Acal$-differentiable in the $D_1$ sense on $\Acal$. Observe $u = c_1$ and $v = c_2+y\frac{dc_1}{dx}$ where both $c_1$ and $c_2$ are real-valued functions of $x$ alone describe the general form of a $\Acal$-differentiable function in the $D_1$ sense. \\

\noindent
For example, setting $c_1=x$ and $c_2=0$ provides the function $f(x+\epsilon y) = x+y \epsilon$. The $D_1$ derivative is simply the constant function $f' = 1$ on $\Acal$. Let us study the $D_2$ differentiability of $f$ at $z_o=x_o+y_o \epsilon$. First, note $(a-b\epsilon/a)(a+b\epsilon) = a^2$ hence for $a \neq 0$
\begin{equation} \label{eqn:zerodivisorsdualnumbers}
 \frac{1}{a+b\epsilon} = \frac{a-b\epsilon/a}{a^2}. 
\end{equation}
We use this identity to begin the calculation below: for $x \neq x_o$,
\begin{align} 
 \frac{f(x+y\epsilon)-f(x_o+y_o\epsilon)}{(x-x_o)+\epsilon(y-y_o)} &= \frac{ \left[x-x_o+(y-y_o)\epsilon\right]\left[x-x_o-(y-y_o)\epsilon /(x-x_o) \right]}{(x-x_o)^2}  \\ \notag 
 &= 1+\epsilon\left[ \frac{y-y_o}{x-x_o} - \frac{y-y_o}{(x-x_o)^2}  \right] 
\end{align}
Notice Equation \ref{eqn:zerodivisorsdualnumbers} shows  $(\RN \oplus \epsilon \RN)^{ \times} = \{ x+y\epsilon \ | \ x \neq 0 \}$. Thus, we study how the difference quotient of $f$ behaves as $x+y\epsilon \rightarrow x_o+\epsilon y_o$ for $x \neq x_o$. Observe the $1$ agrees with the $D_1$ derivative. However, the remaining terms do not converge in the deleted limit hence $f$ is not $D_2$ at $x_o+\epsilon y_o$. But, $x_o+y_o \epsilon$ is arbitrary so we have shown $f$ is {\bf nowhere} $D_2$.
\end{exa}

\noindent
It seems for general finite dimensional commutative unital algebras over $\RN$ it may be difficult or even impossible to obtain nontrivial functions of $D_2$ type. Fortunately, we are free to study $D_1$-differentiability as it includes $D_2$-functions when they exist. \\

\noindent
Much of the literature on hypercomplex variables is largely centered on semisimple algebras. Upto isomorphism in the commutative case we face $\Acal = \RN^n \times \CN^m$. In such a context, it can be shown the set of $D_1$ and $D_2$ differentiable functions on an open set coincide. 

\begin{thm} \label{thm:D2sameasD1forsemisimple}
Let $U$ be an open set in a commutative semisimple finite dimensional real algebra $\Acal$. The set of $D_1$ functions on $U$ coincides with the set of $D_2$ functions on $U$.
\end{thm}

\noindent
{\bf Proof:} Suppose $\Acal$ is a commutative semisimple finite dimensional real algebra and $U \subseteq \Acal$ is open. In Theorem \ref{thm:D2holomorphicisalsoD1} we showed that the set of $D_2$-differentiable functions on $U$ are a subset of the $D_1$-differentiable functions on $U$.  It remains to show $D_1$ functions on $U$ are necessarily $D_2$ functions on $U$. Our proof involves several steps. First, we show the results hold for the direct product algebra $\RN^n$. Second, we show the result holds for the direct product $\CN^m$. Third, Wedderburn's Theorem tells us $\Acal \approx \RN^n \times \CN^m$ and we show how our result filters naturally through the isomorphism to complete the proof. \\

\noindent
Consider $\Acal_R = \RN^n$ with $U_R$ open in $\Acal_R$. Suppose $f = (f_1,\dots , f_n)$ is $D_1$-differentiable on $\Acal_R$. Then $f$ is differentiable on $U_R$ and as the regular representation of $\Acal_R$ is formed by diagonal matrix we find the Cauchy Riemann equations simply indicate that $f_j$ is a function of $x_j$ alone\footnote{as is often notated $f_j = f(x_j)$.}. In total,
\begin{equation}
 f(x) = (f_1(x_1), f_2(x_2), \dots , f_n(x_n)). 
\end{equation}
Moreover, differentiability on $U_R$ implies the partial derivatives of $f$ likewise exist on $U_R$ hence $f_j$ is a real-differentiable function of $x_j$ for $j=1,2,\dots, n$. Notice, $\mathds{1} \in \RN^n$ has the explicit form $\mathds{1} = (1,1,\dots, 1)$ and it follows for $h = (h_1,h_2, \dots , h_n) \neq 0$
\begin{equation}
 \frac{1}{h} = \left( \frac{1}{h_1},\frac{1}{h_2}, \dots, \frac{1}{h_n}\right). 
\end{equation}
Consider the difference quotient at $p = (p_1,p_2, \dots , p_n)$. Define $\triangle_j f = f_j(p_j+h_j)-f_j(p_j)$
\begin{align} 
\frac{f(p+h)-f(p)}{h} &= \left( \frac{1}{h_1},\frac{1}{h_2}, \dots, \frac{1}{h_n}\right)\left(
\triangle_1 f,\triangle_2 f, \dots, \triangle_n f \right) \\ \notag
&= \left( \frac{\triangle_1 f}{h_1},\frac{\triangle_2 f}{h_2}, \dots, \frac{\triangle_n f}{h_n}\right)
\end{align}
To prove $f$ is $D_2$ at $p$ we must show the limit of the  difference quotient exists as $h \rightarrow 0$ for $h \in \Acalx_R$. The condition $h \in \Acalx_R$ simply requires $h_j \neq 0$ for all $j = 1, 2, \dots , n$. Differentiability of $f_j$ at $p_j$ gives
\begin{equation}
 \lim_{h_j \rightarrow 0} \frac{\triangle_j f}{h_j} = \lim_{h_j \rightarrow 0} \frac{f_j(p_j+h_j)-f_j(p_j)}{h_j} = f_j'(p_j). 
\end{equation}
Notice, the condition that $h_i \neq 0$ for $i=1,2, \dots , n$ has no bearing on the limit above. Hence,
\begin{equation}
 \lim_{ \overset{h \rightarrow 0}{h \, \in \, \Acalx_R}} \frac{\triangle_j f}{h_j} = f_j'(p_j).
\end{equation}
Since this holds for each component of $\triangle f/h$ we find
\begin{equation}
 \lim_{ \overset{h \rightarrow 0}{h \, \in \, \Acalx_R}}\frac{f(p+h)-f(p)}{h} = 
\left( f_1'(p_1),f_2'(p_2), \dots, f_n'(p_n) \right)
\end{equation}
Therefore, $f$ is $D_2$ at $p \in U_R$. But, $p$ was arbitrary hence $f$ is $D_2$ on $U_R$. \\

\noindent
Next, if $\Acal_C = \CN^m$ and $U_C$ is open in $\Acal_C$ we consider $g$ which is $D_1$ on $U_C$. Extending the result seen in Example \ref{Ex:number8} we find the Jacobian matrix of $g$ will be block-diagonal with $m$-blocks of the form $\left[ \begin{array}{cc}a_j & -b_j \\ b_j & a_j \end{array}\right]$ for $j=1,2, \dots , m$. The $j$-th diagonal block serves to give the ordinary Cauchy Riemann equations for $g_j$. The zero blocks in for $g_j$ serve to indicate $g_j$ is a function of $z_j$ alone. Here we use $(z_1,z_2, \dots , z_m)$ as the variable on $\Acal_C$. In summary, $g = (g_1, g_2, \dots , g_m)$ where $g_j$ is a complex differentiable function of $z_j$ alone. Moreover,  we may follow the arguments for $\Acal_R$ simply replacing real with complex limits. We find,
\begin{equation}
 \lim_{ \overset{h \rightarrow 0}{h \, \in \, \Acalx_C}}\frac{g(p+h)-g(p)}{h} = 
\left( g_1'(p_1),g_2'(p_2), \dots, g_m'(p_m) \right)
\end{equation}
where $g_j' = \frac{dg_j}{dz_j}$ are complex derivatives.  \\

\noindent
If $\Bcal = \RN^n \times \CN^m$ then we can fit together our result for $\RN^n$ and $\CN^m$ if we make the usual identification that $x \in \RN^n$ and $z \in \CN^m$ gives $(x,z) \in \RN^n \times \CN^m$. Notice $(x,z) \in \Bcalx$  only if $x \in (\RN^n)^{ \times}$ and $z \in (\CN^m)^{\times}$. Moreover,
\begin{equation}
 \frac{1}{(x,z)} = \left( \frac{1}{x}, \frac{1}{z} \right) = \left( \frac{1}{x_1}, \dots , \frac{1}{x_n}, \frac{1}{z_1}, \dots , \frac{1}{z_m} \right).
\end{equation}
It follows that if $(f,g)$ is $D_1$ differentiable on $U$ open in $\Bcal$ then $(f,g)$ is $D_2$ differentiable with $(f,g)' = (f_1', \dots, f_n', g_1', \dots , g_m')$ on $U$. \\

\noindent
Finally, if $\Acal$ is commutative and semisimple associative algebra over $\RN$ then Wedderburn's Theorem\footnote{see Dummit and Foote page 854-855, Theorem 4 part (5.) in \cite{DF}} provides an isomorphism of $\Acal$ and $\Bcal = \RN^n \times \CN^m$ for some $n,m \in \NN$. Suppose $\Psi: \Acal \rightarrow \Bcal$ provides the isomorphism. If $U$ is open in $\Acal$ then $\Psi(U) = U'$ is open in $\Bcal$. Furthermore, suppose $F$ is $D_1$ with respect to $\Acal$ on $U$. Apply Theorem \ref{thm:difftransfer} at each point in $U$ to find that $G = \Psi \comp f \comp \Psi^{-1}$ is $D_1$ with respect to $\Bcal$ at each point in $U'$. Therefore, $G$ is $D_2$ differentiable on $U'$ as we have already shown $D_1$ implies $D_2$ for an open subset of $\RN^n \times \CN^m$. Theorem \ref{thm:isomorphismpreservesstructure} provides that $\Psi$ preserves differences and multiplicative inverses and $  \Psi \comp f =  G \comp \Psi $ thus:
\begin{equation}
 \Psi \left( \frac{f(p+h)-f(p)}{h} \right)  = \frac{\Psi(f(p+h))- \Psi(f(p))}{\Psi(h)}  = \frac{G(\Psi(p+h)) - G(\Psi(p))}{\Psi(h)}.
\end{equation}
Consequently,
\begin{equation}
 \frac{f(p+h)-f(p)}{h}  = \Psi^{-1} \left( 
\frac{G(\Psi(p)+\Psi(h)) - G(\Psi(p))}{\Psi(h)}\right) 
\end{equation}
If $p \in U$ then $\Psi(p) \in U'$ where $G$ is $D_2$ differentiable. Note $h \rightarrow 0$ with $h \in \Acalx$ implies $\Psi(h) \rightarrow 0$ with $\Psi(h) \in \Bcalx$. Using continuity of $\Psi^{-1}$ and that $G$ is $D_2$ at $\Psi(p)$ we find
\begin{equation}
 \lim_{ \overset{h \rightarrow 0}{h \, \in \, \Acalx}} \frac{f(p+h)-f(p)}{h} = \Psi^{-1} \left( 
\lim_{ \overset{\Psi(h) \rightarrow 0}{\Psi(h) \, \in \, \Bcalx}}\frac{G(\Psi(p)+\Psi(h)) - G(\Psi(p))}{\Psi(h)} \right)  = \Psi^{-1}(G'(\Psi(p))). 
\end{equation}
Therefore, $f$ is $D_2$ at $p$ with $f'(p) = \Psi^{-1}(G'(\Psi(p)))$. $\Box$ \\ 

\noindent
In conclusion, the distinction between $D_1$ and $D_2$ differentiability is lost in the commutative semisimple case. This is reflected in Definition 2.11 of \cite{gadeaD1vsD2}. However, if we drop the semisimple condition and just consider general real associative algebras then we argue from Theorem \ref{thm:D2holomorphicisalsoD1} and Example  \ref{exa:nowhereD1} that $D_1$-differentiability provides a more general concept of differentiation over an algebra. For these reasons we take Definition \ref{defn:Adiff} as primary. 

\section{Higher $\Acal$-derivatives} \label{sec:higherderivatives}
The calculus of higher derivatives for functions on $\RN^n$ requires the study of symmetric multilinear maps\footnote{see, for example,  Zorich, Mathematical Analysis II, see Section 10.5 pages 80-87. The results we claim without proof from advanced calculus can all be found in \cite{zorich}.}. However, in $\Acal$-calculus this is avoided due to a fortunate isomorphism between $\Acal$ and symmetric multi-$\Acal$-linear mappings of $\Acal$. Let us begin by generalizing $\EndA$ to its multilinear analog:

\begin{de} \label{defn:multilinearmaps}
We say $T: \underbrace{\Acal \times \cdots \times \Acal}_{k} \rightarrow \Acal$ is a {\bf $k$-linear map on $\Acal$} if $T$ is right-$\Acal$-linear in each of its arguments. That is, $T$ is additive in each entry and
$ T(z_1,\dots, z_j \star w, \dots, z_n) = T(z_1,\dots, z_j, \dots, z_n) \star w. $ for all $z_1, \dots, z_n, w \in \Acal$.
\end{de}

\noindent
If $T(v_{\sigma(1)}, \dots , v_{\sigma(k)}) = T(v_1, \dots , v_k)$ for all permutations $\sigma$ then $T$ is {\bf symmetric}. We continue to assume $\Acal$ is a unital, associative and finite-dimensional algebra over $\RN$.

\begin{thm} \label{thm:multilinearisomorphictoA}
The set of symmetric $k$-linear maps on $\Acal$ is isomorphic to $\Acal$.\end{thm}

\noindent
{\bf Proof:} the sum and scalar multiple of symmetric $k$-linear map is once more $k$-linear and symmetric. Since $v_j = \mathds{1} \star v_j$ we find:
\begin{equation}
 T(v_1, \dots, v_k) = T(\mathds{1},\dots , \mathds{1}) \star v_1 \star \cdots \star v_k. 
\end{equation}
thus $T$ is uniquely fixed by $k$-linearity on $\Acal$ together with its value on $(\mathds{1},\dots , \mathds{1})$. $\Box$ \\

\noindent
We already saw this argument in Equation \ref{eqn:theformula} where we proved the $k=1$ case $\Acal \approx \EndA$. 

\begin{de} \label{defn:higherderivative}
Suppose $f$ is a function on $\Acal$ for which the derivative function $f'$ is $\Acal$-differentiable at $p$ then we define
$ f''(p) = (f')'(p)$. Furthermore, supposing the derivatives exist, we define $f^{(k)}(p) = (f^{(k-1)})'(p)$ for $k =2,3, \dots$. 
\end{de}

\noindent
Naturally we define functions $f'', f''', \dots, f^{(k)}$ in the natural pointwise fashion for as many points as the derivatives exist. Furthermore, with respect to $\beta = \{ v_1, \dots , v_n \}$ where $v_1 = \mathds{1}$, we have $f'(p) = d_pf(\mathds{1}) = \frac{\partial f}{\partial x_1}(p)$.  Thus, $f' = \frac{\partial f}{\partial x_1}$. Suppose $f''(p)$ exists. Note,
\begin{equation}
 f''(p) = (f')'(p) = \#( d_p f'(\mathds{1}) )= \frac{\partial f'}{\partial x_1}(p) = \frac{\partial^2f}{\partial x_1^2}(p). 
\end{equation}
Thus, $f'' = \frac{\partial^2 f}{\partial x_1^2}$. By induction, we find the following theorem:

\begin{thm} \label{thm:highderivativeallone}
If $f: \Acal \rightarrow \Acal$,  $\beta = \{ \mathds{1}, \dots , v_n\}$ a basis, and $f^{(k)}$ exists then $f^{(k)} = \frac{\partial^k f}{\partial x_1^k}$.
\end{thm}

\noindent
The algebra derivatives naturally dovetail with the iterated-symmetric-Frechet differentials which are used to describe higher derivatives of a map on normed linear spaces\footnote{The iterated-differentials are developed in many advanced calculus texts. See \cite{zorich} where the theory of real higher derivatives  is developed in Section 10.5 pages 80-87.}.

\begin{thm} \label{thm:rightAconspriacy}
Suppose $f: \Acal \rightarrow \Acal$ is a function for which $f^{(k)}(p)$ exists. Then the iterated $k$-th Frechet differential exists and is related to the $k$-th $\Acal$ derivative as follows: 
$$ d_p^kf(v_1, \dots, v_k) = f^{(k)}(p) \star v_1 \star \cdots \star v_k. $$
for all $v_1, \dots, v_k \in \Acal$.
\end{thm}

\noindent
{\bf Proof:} Suppose $f: \Acal \rightarrow \Acal$ is a function for which $f^{(k)}(p)$ exists. The existence of the iterated $\Acal$-derivatives implies that $f$ is also $k$-fold $\RN$-differentiable and thus $d^k_pf: \Acal \times \cdots \times \Acal \rightarrow \Acal$ exists and is a real symmetric $k$-linear map. Let  $\beta = \{ v_1, \dots , v_n \}$ with $v_1 = \mathds{1}$ be a basis for $\Acal$ with coordinates $x_1, \dots , x_n$. The iterated $k$-th Frechet differential and iterated partial derivatives are related by:
\begin{equation} \label{eqn:iteratedFrechetdiffandpartial}
d^k_pf(v_{i_1},v_{i_2}, \dots , v_{i_k}) = \frac{\partial^k f}{\partial x_{i_1}\partial x_{i_2} \cdots \partial x_{i_k}}. 
\end{equation}
Differentiating the $\Acal$-CR equations
$\frac{\partial f}{\partial x_j} = \frac{\partial f}{\partial x_1} \star v_j$ with respect to $x_i$ yields:
\begin{align} \notag
\frac{\partial^2 f}{\partial x_i\partial x_j}  = \frac{\partial}{\partial x_i} \left[ \frac{\partial f}{\partial x_j} \right] &= \frac{\partial}{\partial x_i} \left[ \frac{\partial f}{\partial x_1} \star v_j \right] \\ \notag
 &= \frac{\partial}{\partial x_1} \left[ \frac{\partial f}{\partial x_i} \right] \star v_j \\ \label{eqn:proofsymmetric}
 &= \frac{\partial^2 f}{\partial x_1^2} v_i \star v_j  
 \end{align}
Apply Equation \ref{eqn:proofsymmetric} repeatedly as to exchange partial derivatives with respect to $x_{i_j}$ for partial derivatives with respect to $x_1$ and multiplication by $v_{i_j}$ obtain:
\begin{equation} \label{eqn:partialstomultiples}
 \frac{\partial^k f}{\partial x_{i_1}\partial x_{i_2} \cdots \partial x_{i_k}} = \frac{\partial^k f}{\partial x_1^k} \star v_{i_{1}} \star v_{i_{2}} \star \cdots \star v_{i_{k}}  = f^{(k)} \star v_{i_{1}} \star v_{i_{2}} \star \cdots \star v_{i_{k}} 
\end{equation}
We used Theorem \ref{thm:highderivativeallone} in the last step. Compare Equations \ref{eqn:iteratedFrechetdiffandpartial} and \ref{eqn:partialstomultiples} to conclude the proof.   $\Box$ \\

\noindent
In fact, the first equality in Equation \ref{eqn:partialstomultiples} should be emphasized:

\begin{thm} \label{thm:partialsproducts}
If $f: \Acal \rightarrow \Acal$ is $k$-times $\Acal$-differentiable then 
$$ \frac{\partial^k f}{\partial x_{i_1}\partial x_{i_2} \cdots \partial x_{i_k}} = \frac{\partial^k f}{\partial x_1^k} \star v_{i_{1}} \star v_{i_{2}} \star \cdots \star v_{i_{k}}.$$
\end{thm}

\noindent
Theorem \ref{thm:rightAconspriacy} and  \ref{thm:partialsproducts} provide the basis for both the formulation of an $\Acal$-variable Taylor Theorem as well as a lucid derivation of generalized Laplace Equations.

\subsection{The $\Acal$-harmonic equations}
\noindent
In the case of complex analysis the second order differential consequences of the Cauchy Riemann equations include the Laplace equations. It is interesting to determine what equations form the analog to Laplace's Equation for $\Acal$. 
In 1948 Wagner derived {\bf generalized Laplace Equations} in \cite{wagner1948} via calculations performed through the lens of the {\it paraisotropic matrix}. Then, in 1992, Waterhouse derived the same results by using the trace on a Frobenius algebra \cite{waterhouseII}. In both cases, the argument is essentially a pairing of the commutativity of mixed real partial derivatives and the generalized Cauchy Riemann equations. 

\begin{thm} \label{thm:equationtoDEqn}
Let $U$ be open in $\Acal$ and suppose $f: U \rightarrow \Acal$ is twice $\Acal$-differentiable on $U$. If there exist $B_{ij} \in \RN$ for which $\sum_{i,j} B_{ij}v_i \star v_j = 0$ then $\sum_{i,j}B_{ij}\frac{\partial^2 f}{\partial x_i \partial x_j}=0$.
\end{thm}

\noindent
{\bf Proof:} suppose $f$ is twice continuously $\Acal$-differentiable on $U \subset \Acal$ and suppose there exist $B_{ij} \in \RN$ for which $\sum_{i,j} B_{ij}v_i \star v_j = 0$. Multiply the given equation by $\frac{\partial^2 f}{\partial x_1^2}$ to obtain:
\begin{equation}
  \sum_{i,j} B_{ij}\frac{\partial^2 f}{\partial x_1^2} \star v_i \star v_j = 0. 
\end{equation}
Then, by Equation \ref{eqn:partialstomultiples} we deduce $\frac{\partial^2 f}{\partial x_1^2} \star v_i \star v_j = \frac{\partial^2 f}{\partial x_i \partial x_j}$. Therefore,
 $ \ds \sum_{i,j} B_{ij}\frac{\partial^2 f}{\partial x_i\partial x_j}=0$.  $\Box$

\noindent
Theorem \ref{thm:equationtoDEqn} essentially says that a quadratic equation in $\Acal$ yields a corresponding Laplace-type equation for $\Acal$-differentiable functions.  Hence we find:

\begin{coro} \label{thm:laplacepatternmatch}
Generalized Laplace equations can be assembled by mimicking patterns in the multiplication table for $\Acal$ to matching patterns in the Hessian matrix. Moreover, each component of an $\Acal$-differentiable function is a solution to the generalized Laplace equations.
\end{coro}

\noindent
This result was given by Wagner in \cite{wagner1948}.

\begin{exa}
Consider $\Acal = \RN \oplus j \RN\oplus j^2 \RN$ where $j^3=1$. Notice, we have multiplication table and Hessian matrix
\begin{equation}
  \begin{array}{ c || c | c |c|}
     & 1 & j & j^2 \\ \hline \hline
    1 & 1 & j & j^2\\ \hline
    j & j& j^2 & 1 \\ \hline
    j^2 & j^2 & 1 & j \\ \hline
  \end{array} \qquad \& \qquad 
  \begin{array}{ c || c | c |c|}
     & x & y & z \\ \hline \hline
    x & f_{xx} & f_{xy} & f_{xz}\\ \hline
    y& f_{yx}& f_{yy} & f_{yz} \\ \hline
    z & f_{zx} & f_{yz} & f_{zz} \\ \hline
  \end{array}  
 \end{equation}
Theorem \ref{thm:laplacepatternmatch} allows us to find the following generalized Laplace equations by inspection of the tables above:
\begin{equation}
 f_{xx}= f_{yz}, \ \ f_{xy}= f_{zz}, \ \ f_{xz}= f_{yy} 
\end{equation}
You might recognize these from Example \ref{exa:trihyperbolicidentities}.
\end{exa}

\begin{exa}
Consider $\Acal = \RN \oplus i \RN$ where $i^2 = -1$. Notice, we have multiplication table and Hessian matrix
\begin{equation}  \begin{array}{ c || c | c |}
     & 1 & i  \\ \hline \hline
    1 & 1 & i \\ \hline
    i & i& -1  \\ \hline
  \end{array} \qquad \& \qquad 
  \begin{array}{ c || c | c |}
     & x & y  \\ \hline \hline
    x & f_{xx} & f_{xy} \\ \hline
    y& f_{yx}& f_{yy}  \\ \hline
  \end{array}  
 \end{equation}
Theorem \ref{thm:laplacepatternmatch} allows us to find the Laplace equations by inspection of the tables above: if $f=u+iv$ then
\begin{equation}
 f_{xx}= -f_{yy}, \ \ \Rightarrow \ \ u_{xx}+u_{yy}=0 \ \ \& \ \ v_{xx}+v_{yy}=0.  
\end{equation}
\end{exa}

\noindent
If we consider the component formulation of the Cauchy Riemann equations then differentiation of these equations will produce second order homogeneous PDEs which include the generalized Laplace Equations and other less elegant equations coupling distinct components. The elegance of the generalized Laplace equations is seen in the fact that {\bf every} component of an $\Acal$-differentiable function is what we may call $\Acal$-harmonic. In other words, the algebra $\Acal$ provides a natural function theory to study $\Acal$-harmonic functions. Notice, the concept of $\Acal$-harmonicity involves solving a system of PDEs. When is it possible to find an $\Acal$ for which a given system of real PDEs for the $\Acal$-harmonic equations for $\Acal$? Theorem \ref{thm:laplacepatternmatch} gives at least a partial answer. If we replicate patterns imposed on the Hessian matrix to produce a multiplication table then we can test if the table is a possible multiplication table for an algebra. It is interesting to note that Ward already solved the corresponding problem for generalized Cauchy Riemann equations in 1952. In particular, Ward showed in \cite{ward1952} how to construct an algebra $\Acal$ which takes a given set of $n^2-n$ independent PDEs as its generalized $\Acal$-CR equations. \\

\noindent
The next example was inspired by Example 4.6 in \cite{waterhouseII}. It demonstrates how algebraic insight can be wielded to produce solutions to PDEs.

\begin{exa} \label{exa:waveqn}
Consider the wave equation $c^2u_{xx} = u_{tt}$ where $c$ is a positive constant which characterizes the speed of the transverse waves modelled by this PDE. Let us find an algebra $\mathcal{W}_c$ which takes the speed-$c$ wave equation as its generalized Laplace Equation. Let $(x,t) = x+kt$ form a typical point in the algebra. What rule should we give to $k$? Following Corollary \ref{thm:laplacepatternmatch} we should consider the correspondence:
\begin{equation}
c^2u_{xx} = u_{tt} \ \ \ \leftrightarrow \ \ \ c^2 = k^2
\end{equation}
thus set $k^2=c^2$. The algebra $\mathcal{W}_c = \RN \oplus k \RN$ with $k^2=c^2$ has $\mathcal{W}_c$-differentiable functions $f = u+kv$ for which $c^2u_{xx} = u_{tt}$. Observe $\Gamma: \mathcal{W}_c \rightarrow \Hcal$ defined by $\Gamma( x+kt) = x+cjt$ serves as an isomorphism of $\mathcal{W}_c$ and the hyperbolic numbers of Example \ref{Ex:number2}. Combine $\Psi^{-1}(x+jy) = (x+y,x-y)$ of Example \ref{Ex:number2} with $\Gamma$ to construct the isomorphism $\Phi = \Psi^{-1} \comp \Gamma$ from $\mathcal{W}_c$ to $\RN \times \RN$. In particular,
\begin{equation}
\Phi(x+kt) = \Psi^{-1}(  \Gamma (x+kt)) = \Psi^{-1}(x+cjt) = (x+ct, x-ct)   
 \end{equation}
Following the insight of Theorem \ref{thm:difftransfer} we associate to each $\mathcal{W}_c$-differentiable function $f: \mathcal{W}_c \rightarrow \mathcal{W}_c$ a corresponding $\RN \times \RN$ differentiable function $F$ as follows:
\begin{equation}
 f = \Phi^{-1} \comp F \comp \Phi  
 \end{equation}
where $F: \RN \times \RN \rightarrow \RN \times \RN$. The structure of $\RN \times \RN$-differentiable functions is rather simple; $F(a,b) = (F_1(a), F_2(b))$ where $F_1,F_2$ are differentiable functions on $\RN$. Thus,
\begin{align}
 f(x+kt) &= \Phi^{-1}(  F ( \Phi(x+kt))) \\ \notag  
&= \Phi^{-1}(  F ( (x+ct, x-ct))) \\ \notag  
&= \Phi^{-1}(  F_1 (x+ct),F_2(x-ct)) \\ \notag 
&= \frac{1}{2}(F_1 (x+ct)+F_2(x-ct))+ \frac{k}{2c}(F_1 (x+ct)-F_2(x-ct))
 \end{align}
 We have shown that $\Acal$-differentiable functions $f= u+kv$ have (using $x+kt = (x,t)$ to make the formulas more recognizable)
\begin{equation}
u(x,t) = \frac{1}{2}(F_1 (x+ct)+F_2(x-ct)) \qquad \& \qquad v(x,t) = \frac{k}{2}(F_1 (x+ct)-F_2(x-ct)). 
\end{equation}
We've shown how d'Alembert's solution to the wave-equation appears naturally in the function theory of $\mathcal{W}_c$. 
\end{exa}

\noindent
Naturally, there are higher order versions of the Laplace Equations.

\begin{thm}
Let $U$ be open in $\Acal$ and suppose $f: U \rightarrow \Acal$ is $k$-times $\Acal$-differentiable. If there exist $B_{i_1i_2\dots i_k} \in \RN$ for which $\sum_{i_1i_2\dots i_k} B_{i_1i_2\dots i_k}v_{i_1} \star v_{i_2} \star \cdots \star v_{i_k} = 0$ then 
$$\sum_{i_1i_2\dots i_k} B_{i_1i_2\dots i_k}\frac{\partial^k f}{\partial x_{i_1} \partial x_{i_2} \cdots  \partial x_{i_k} }=0. $$
\end{thm}

\noindent
{\bf Proof:} multiply the assumed relation by $\frac{\partial^kf}{\partial x_1^k}$ and apply Theorem \ref{thm:partialsproducts}. $\Box$ \\

\subsection{$\Acal$-variate Taylor's Theorem}

\noindent
If we are given that $f: \Acal \rightarrow \Acal$ is smooth in the sense of real analysis then it is simple to show that the existence of the first $\Acal$-derivative implies the existence of all higher $\Acal$-derivatives.

\begin{thm} \label{thm:fromonemany}
Let $\Acal$ be a commutative unital finite dimensional algebra over $\RN$. Suppose $f: \Acal \rightarrow \Acal$ has arbitrarily many continuous real derivatives at $p$ and suppose $f$ is once $\Acal$-differentiable at $p$ then $f^{(k)}(p)$ exists for all $k \in \NN$. \end{thm}

\noindent
{\bf Proof:} suppose $f$ is smooth and once $\Acal$-differentiable at $p$. We assume $\Acal$ is a commutative unital algebra over $\RN$ with basis $\beta = \{ \mathds{1}, \dots, v_n \}$. Assume inductively that $f^{(k)}(p)$ exists hence Theorem \ref{thm:highderivativeallone} provides $f^{(k)}(p) = \frac{\partial^k f(p)}{\partial x_1^k}$. Consider, omitting $p$ to reduce clutter,
\begin{equation}
 \frac{\partial f^{(k)}}{\partial x_j} = \frac{\partial}{\partial x_j} \left[ \frac{\partial^k f}{\partial x_1^k} \right] =
\frac{\partial^k}{\partial x_1^k} \left[ \frac{\partial f}{\partial x_j} \right]  = \frac{\partial^k}{\partial x_1^k} \left[ \frac{\partial f}{\partial x_1}  \right] \star v_j = \frac{\partial f^{(k)}}{\partial x_1} \star v_j. 
\end{equation}
Thus $f^{(k)}$ is $\Acal$-differentiable at $p$ which proves $f^{(k+1)}(p)$ exists.  $\Box$ \\

\noindent
In \cite{zorich} a multivariate Taylor's Theorem over a finite dimensional normed linear space is given. In particular, if $f: V \rightarrow V$ is {\bf real analytic} then $f$ is represented by its multivariate Taylor series on some open set containing $p$. The multivariate Taylor series of $f$ centered at $p$ is given, for $h$ sufficiently small, by the convergent series\footnote{equivalently, $ f(p+h) = f(p)+ \sum_{k=1}^{\infty}\sum_{i_1, \dots, i_k} \frac{1}{k!} \frac{\partial^k f (p)}{\partial x_{i_1}\partial x_{i_2} \cdots \partial x_{i_k}} h_{i_1}h_{i_2} \cdots h_{i_k}. $}:  
\begin{equation} \label{eqn:multivariateTaylorSeries}
f(p+h) = f(p)+ d_pf(h)+\frac{1}{2}d^2_pf(h,h)+ \frac{1}{3!}d^3_pf(h,h,h) + \cdots
\end{equation}
Notice, for $(h,\dots , h) \in \Acal^k$ we may expand $h = \sum_{i_j} h_{i_j}v_{i_j}$ for $j=1,2, \dots , k$,
\begin{equation} \label{eqn:kthtermmultiTaylor}
 d^k_pf(h, \dots, h) = \sum_{i_1, \dots, i_k}  h_{i_1}h_{i_2} \cdots h_{i_k}d^k_pf(v_{i_1},v_{i_2}, \dots , v_{i_k})
 \end{equation}

\noindent
Compare the formulas above to the form of the $k$-term in the Taylor expansion for an $\Acal$-differentiable function given below. This simplification has Theorem \ref{thm:multilinearisomorphictoA} at its root.

\begin{thm} \label{thm:AcalTaylor}
Let $\Acal$ be a commutative, unital, associative algebra over $\RN$. If $f$ is real analytic at $p \in \Acal$ then
$$  f(p+h) = f(p)+f'(p) \star h + \frac{1}{2}f''(p) \star h^2+ \cdots + \frac{1}{k!}f^{(k)}(p) \star h^k+ \cdots $$
where $h^2 = h \star h$ and $h^{k+1} = h^{k} \star h$ for $k \in \NN$.
\end{thm}

\noindent
{\bf Proof:} Suppose $f$ is real analytic and $\Acal$-differentiable. Since real analytic implies $f$ is smooth over $\RN$ we apply Theorem \ref{thm:fromonemany} to see $f$ is smooth over $\Acal$. Therefore, we may follow the proof of Theorem \ref{thm:rightAconspriacy} and obtain: 
\begin{equation} \label{eqn:taylorproofI}
 d^k_pf(v_{i_1},v_{i_2}, \dots , v_{i_k}) = f^{(k)}(p) \star v_{i_1} \star v_{i_2} \star  \dots \star  v_{i_k}. 
 \end{equation}
Observe the $k$-th power of $h \in \Acal$ is given by
\begin{equation} \label{eqn:taylorproofII}
 h^k = \sum_{i_1, \dots , i_k} h_{i_1} h_{i_2}\cdots h_{i_k}v_{i_1} \star v_{i_2} \star  \dots \star  v_{i_k}. 
\end{equation}
Therefore, combining Equations \ref{eqn:kthtermmultiTaylor}, \ref{eqn:taylorproofI}, and  \ref{eqn:taylorproofII} we find
\begin{align} 
 d^k_pf(h, \dots, h) &= \sum_{i_1, \dots, i_k}  h_{i_1}h_{i_2} \cdots h_{i_k}d^k_pf(v_{i_1},v_{i_2}, \dots , v_{i_k}) \\ \notag
&= \sum_{i_1, \dots, i_k}  h_{i_1}h_{i_2} \cdots h_{i_k}f^{(k)}(p) \star v_{i_1} \star v_{i_2} \star  \dots \star  v_{i_k}. \\ \notag
&=  f^{(k)}(p) \star \sum_{i_1, \dots, i_k}  h_{i_1}h_{i_2} \cdots h_{i_k} v_{i_1} \star v_{i_2} \star  \dots \star  v_{i_k} \\ \notag
&= f^{(k)}(p)\star h^k. 
\end{align}
We conclude, $ f(p+h) = f(p)+f'(p) \star h + \frac{1}{2}f''(p) \star h^2+ \cdots + \frac{1}{k!}f^{(k)}(p) \star h^k+ \cdots $. $\Box$ \\

\noindent
We study the theory of convergence series in $\Acal$ in the sequel to this paper which is a joint work with Daniel Freese \cite{cookfreese}.

\section{An approach to the inverse problem} \label{sec:inverseproblem}
The inverse problem of $\Acal$-calculus is roughly this:
\begin{quote}
When can we translate a problem of real calculus to a corresponding problem of $\Acal$-calculus ?
\end{quote}
Naturally, this raises a host of questions. What kind of real calculus problems? How do we choose $\Acal$? In  \cite{pagr2012} and \cite{pagr2015} the authors study how to modify certain ODEs in terms of $\Acal$-calculus. In contrast, our study on the inverse problem has centered around systems of PDEs. In particular, we seek to answer the following question:
\begin{quote}
When can we find solutions to a system of PDEs which are simultaneously solutions to the generalized Laplace equations of some algebra $\Acal$ ?
\end{quote}

\noindent
Our hope is that if the answer to the question above is affirmative then it may be possible to rewrite the system of PDEs in an $\Acal$-based notation where the PDE in many variables simply becomes an $\Acal$-ODE in a single algebra variable. Of course, this is just an initial conjecture, there are many directions we could explore at the level of algebra-based differential equations\footnote{In fact at the time this paper is prepared the author has already shown how to solve many $\Acal$-ODEs. The joint work \cite{cookbedell} with Nathan BeDell is currently under preparation. Nathan BeDell has three other papers \cite{bedellI},\cite{bedellII},\cite{bedellIII} in preparation which discuss zero-divisors and basic algebra, the construction of logarithms, and identities for generalized trigonmetric functions and many other algebraic preliminaries.}. Let us examine a simple example of how an $\Acal$-ODE can replace a system of PDEs.

\begin{exa} \label{exa:AODEbackwards}
Let $z = x+jy$ denote an independent hyperbolic variable and $w = u+jv$ the solution of $\frac{dw}{dz} = w^2$. Separating variables gives $\frac{dw}{w^2} = dz$ hence $\frac{-1}{w} = z+c$ where $c=c_1+jc_2$ is a hyperbolic constant. Thus, $w = \frac{-1}{z+c}$ is the solution. What does this mean at the level of real calculus? Note, \begin{equation}
 \frac{dw}{dz} = w^2 \ \ \Rightarrow \ \ u_x+jv_x = (u+jv)^2 = u^2+v^2+2juv
\end{equation}
In other words, the $\Acal$-ODE is the nonlinear system of PDEs
\begin{equation} \label{eqn:systemofPDEs}
 u_x = u^2+v^2, \ \ \& \ \ v_x = 2uv 
\end{equation}
paired with the $\Acal$-CR equations $u_x=v_y$ and $u_y=v_x$. We have the solution already from direct calculus on $\Acal$,
\begin{equation}
 w = \frac{-1}{z+c}  \ \ \Rightarrow \ \ u+jv = -\frac{x+c_1-j(y+c_2)}{(x+c_1)^2-(y+c_2)^2}
\end{equation}
Thus, $u = -\frac{x+c_1}{(x+c_1)^2-(y+c_2)^2}$ and $v = \frac{y+c_2}{(x+c_1)^2-(y+c_2)^2}$ are the real solutions to \ref{eqn:systemofPDEs} and you can check that $u_x=v_y$ and $u_y=v_x$ as well.
\end{exa}

\noindent
Our first goal in understanding the inverse problem was to decide when it is possible to pair a system like \ref{eqn:systemofPDEs} with the $\Acal$-CR equations of an appropriate algebra $\Acal$. Essentially, our first concern is whether there is at least an algebra whose $\Acal$-CR equations and their {\it differential consequences} are not inconsistent with a given system of PDEs. \\

\noindent
We construct the {\bf generic tableau} $\Tcal = \Acal \oplus \Tcal_1 \oplus \Tcal_2 \oplus \cdots \oplus \Tcal_k \oplus+ \cdots $ 
\begin{equation}
 \Tcal_k = \text{span} \{ dx^{i_1} \otimes dx^{i_2}\otimes \cdots \otimes dx^{i_k} \otimes  v_j \ | \ 1 \leq i_1 \leq i_2 \leq \dots \leq i_k \leq n, 1 \leq j \leq n \} 
\end{equation}
We use these spaces to account for dependencies amongst variables and their derivatives. In particular, we focus our attention on PDEs which are formed from $n$-dependent and $n$-independent variables. Notice, we only need increasing indices since these symbols represent partial derivatives which we can commute to be in increasing order\footnote{this construction was inspired by a more sophisticated, but similar, construction in Chapter 4 of \cite{cartan4beginners}}. The {\bf Gauss map} of a function $f: \Acal \rightarrow \Acal$ into the generic tableau is formed as follows: if $f = \sum u^jv_j$
\begin{equation}
 \gamma(f) = f \oplus \left(\sum_{i,j}\frac{\partial u^j}{\partial x_i} v_j \otimes dx^i  \right) \oplus \left(\sum_{i_1\leq i_2} \sum_j\frac{\partial^2 u^j}{\partial x_{i_1}\partial x_{i_2}} v_j \otimes dx^{i_1} \otimes dx^{i_2}  \right) \oplus \cdots  
\end{equation}
For the $k$-th term,
\begin{equation}
 \gamma(f) = \cdots \oplus \left( \sum_{i_1\leq i_2 \leq \dots \leq i_k} \sum_j\frac{\partial^k u^j}{\partial x_{i_1}\partial x_{i_2} \cdots \partial x_{i_k}} v_j \otimes dx^{i_1} \otimes dx^{i_2}\otimes \cdots \otimes dx^{i_k}\right) \oplus \cdots 
 \end{equation}
Or, more concisely, 
\begin{align} \label{eqn:gaussmap}
\gamma(f) &= f+ \sum_i (\partial_i f) \, dx^i+ \sum_{i_1 \leq i_2} (\partial_{i_1}\partial_{i_2}f ) \, dx^{i_1} \otimes dx^{i_2}+\\ \notag
& \qquad \cdots + \sum_{i_1 \leq \dots \leq i_k}  (\partial_{i_1}\cdots \partial_{i_k}f) dx^{i_1} \otimes \cdots \otimes dx^{i_k} + \cdots
\end{align}
We seek to represent systems of PDEs as particular subspaces inside $\Tcal$. The PDEs we study have finite order and hence the calculation ultimately amounts to a question of finite dimensional linear algebra. We should note the subspace of $\Tcal$ which is given by the $\Acal$-CR equations is particularly simple.

\begin{thm} \label{thm:gaussAmap}
Let $\Tcal_{\Acal} \leq \Tcal$ denote the subspace of $\Tcal$ generated by infinitely $\Acal$-differentiable functions then $\gamma \in \Tcal_{\Acal}$ has the form
$$ \gamma = \alpha_o + \alpha_1 \star \sum_i  v_i \otimes dx^i  +  \alpha_2 \star  \sum_{i_1 \leq i_2}  v_{i_1} \star v_{i_2} \otimes dx^{i_1} \otimes dx^{i_2}  + \cdots + \alpha_k \star  \sum_{|I|} v_I \otimes dx^I + \cdots $$
where $|I|$ indicates the sum over increasing $k$-tuples of indices taken from $\{ 1,2,\dots, n \}$ and $v_I = v_{i_1} \star v_{i_2} \star \cdots \star v_{i_k}$ and $dx^I = dx^{i_1} \otimes dx^{i_2}\otimes \cdots \otimes dx^{i_k}$.
\end{thm}

\noindent
{\bf Proof:} If $f$ has arbitrarily many $\Acal$-derivatives then we find
\begin{align} 
 \partial_i f &= (\partial_1 f) \star v_i, \\ \notag
 \partial_{i_1} \partial_{i_2} f &= (\partial^2_1 f) \star v_{i_1} \star v_{i_1} \\ \notag
 \partial_{i_1} \partial_{i_2} \cdots \partial_{i_k} f &= (\partial^k_1 f) \star v_{i_1} \star v_{i_1} \star \cdots \star v_{i_k}.
 \end{align}
Then by Equation \ref{eqn:gaussmap} we find $\gamma_f$ has the form given in the Theorem. $\Box$ \\

\section{Integration} \label{sec:integration}

\noindent
Let $U$ be a connected open subset of $\CN$. If $f$ is complex differentiable on $U$ then Goursat's Theorem shows that the derivative function $f'$ is continuous. Furthermore, Cauchy's integral formula shows that $f^{(k)}$ exists for any $k \in \NN$. In summary, if $f$ is once complex differentiable on $U$ then $f$ is infinitely many times complex differentiable. This well-known result stands in marked contrast to the structure of real differentiable functions. Naturally, we wish to understand to what extent this story may be replicated for an associative, unital, finite dimensional algebra over $\RN$. \\

\noindent
If $\Acal = \RN$ then know a function can be once $\Acal$-differentiable and not $\Acal$-smooth. In addition, even when $\Acal \neq \RN$, Cauchy's integral formula does not hold for an arbitrary $\Acal$. For example, both \cite{khrennikov} and \cite{MotterRosa} show that Cauchy's integral formula holds for hyperbolic numbers. \\

\noindent
In this Section we assume $\Acal$ is an associative, unital, commutative algebra of finite dimension over $\RN$ with basis $\beta = \{ v_1, \dots , v_n \}$ where $v_1 = \mathds{1}$.

\begin{de} \label{defn:Avaluedintegral}
Suppose $f: [a,b] \rightarrow \Acal$ is continuous function with $f = u_1v_1+\cdots + u_nv_n$ where $u_j : [a,b] \rightarrow \RN$ for $j=1,2, \dots, n$. We define
$$ \int_a^b f\, dt  = 
 \left(\int_a^b u_1\, dt\right)v_1+ \cdots + \left(\int_a^b u_n \, dt \right)v_n.$$
\end{de}

\noindent
This integral has the expected linearity properties:

\begin{thm} \label{thm:realAintegralprops}
Let $f,g : [a,b] \rightarrow \Acal$ and be continuous and $c \in \Acal$ then
\begin{enumerate}[{\bf (i.)}]
\item $\int_a^b (f+g)\, dt = \int_a^b f\, dt+ \int_a^b g\, dt$
\item $\int_a^b \alpha \star f \, dt = \alpha \star \int_a^b f \, dt$ 
\item for $F$ differentiable, $\int_a^b \frac{dF}{dt} \, dt = F(b)-F(a)$
\end{enumerate}
\end{thm}

\noindent
{\bf Proof:} the proof of (i.) is straightforward. To see (ii.) notice if $c_o \in \RN$ then 
\begin{equation}
 \int_a^b c_o f \,dt = \sum_i \left( \int_{a}^{b}  c_o f_i \, dt \right) v_i = c_o\sum_i\left( \int_{a}^{b} f_i \, dt \right) v_i = c_o \int_{a}^{b} f \, dt. 
\end{equation}
Recall the structure constants $C_{ij}^k \in \RN$ for which $v_i \star v_j = \sum_k C_{ij}^k v_k$ allow us to express the product of $\alpha = \sum_i \alpha_i v_i$ and $f = \sum_j f_j v_j$ as $\alpha \star f = \sum_{i,j,k} \alpha_i f_j C_{ij}^k v_k$. Hence,
\begin{align}
 \int_a^b (\alpha \star f) \, dt  &= 
\sum_k \left(\int_a^b \sum_{i,j} \alpha_i f_j C_{ij}^k \, dt \right) v_k \\ \notag
&= \sum_{k}\left(\sum_{i,j} C_{ij}^k  \alpha_i\int_a^b f_j \, dt \right) v_k  \\ \notag
&= \alpha \star \int_a^b f \, dt  
\end{align}
which proves (ii.). Suppose $F = \sum_i F_i v_i$ is differentiable on $[a,b]$ then by applying the Fundamental Theorem of Calculus on each component we find:
\begin{equation}
 \int_a^b \frac{dF}{dt} \, dt = \sum_i \left( \int_a^b \frac{dF_i}{dt} \, dt \right) v_i = \sum_i (F_i(b)-F_i(a)) v_i = F(b)-F(a). 
\end{equation}
This proves (iii.). $\Box$ \\

\noindent
There are at least two natural ways to define the integral of a function on $\Acal$ along a curve. We focus our attention to smooth curve. If faced with a piecewise smooth curve then we agree to form the integral by taking the sum over smooth pieces. 

\begin{de} \label{defn:Aintegral}
Suppose $f$ is a function on $\Acal$ which is continuous near a rectifiable curve $C$ which begins at $P$ and terminates at $Q$. Let $\zeta_o = P$ and $\zeta_m = Q$. Denote the line-segment from $\zeta_i$ to $\zeta_j$ by $[\zeta_i, \zeta_j]$. If we assume each $\zeta_i \in C$ then the concatenation of $[\zeta_o,\zeta_1], [\zeta_1,\zeta_2], \dots , [\zeta_{m-1},\zeta_m]$ forms a broken line segment which tends to $C$ as $m \rightarrow \infty$. We define,
$$ \int_C f(\zeta) \star d\zeta = \lim_{m \rightarrow \infty} \sum_{i=1}^{m} f(\zeta_i) \star \triangle \zeta_i \ \ \\ \text{where $\triangle \zeta_i = \zeta_i - \zeta_{i-1}$ for $i=1,\dots, m$.} $$
\end{de}

\noindent
Pragmatically, this definition is not of much direct use. However, we begin here as this is a natural generalization of the usual Riemann integral. In addition, it allows an economical proof of the theorem below.

\begin{thm} \label{thm:subML}
Let $C$ be a rectifiable curve with arclength $L$. Suppose $||f(\zeta) || \leq M$ for each $\zeta \in C$ and suppose $f$ is continuous near $C$. Then 
$$ \bigg{|}\bigg{|} \int_C f( \zeta) \star d\zeta \bigg{|}\bigg{|} \leq \Abound ML $$
where $\Abound$ is a constant such that $|| z \star w || \leq \Abound ||z|| \, ||w||$ for all $z,w \in \Acal$.
\end{thm}

\noindent
{\bf Proof:} assume the notation and conditions of the theorem and Definition \ref{defn:Aintegral}. Observe,
\begin{equation}
\bigg{|}\bigg{|}\sum_{i=1}^{m} f(\zeta_i)\star\triangle \zeta_i \bigg{|}\bigg{|} \leq  \sum_{i=1}^{m} || f(\zeta_i) \star \triangle \zeta_i || \leq 
\sum_{i=1}^{m} \Abound ||f(\zeta_i)|| \, || \triangle \zeta_i || \leq \Abound M \sum_{i=1}^{m} || \triangle \zeta_i ||.  
\end{equation}
Notice, as $m \rightarrow \infty$ we find $ \sum_{i=1}^{m} || \triangle \zeta_i ||$ tends to the arclength $L$ of $C$. Since the norm $|| \cdot ||$ is continuous, we may pass the limit inside the norm and the theorem follows. $\Box$ \\

\noindent
If $\zeta = x_1v_1+ \cdots +x_nv_n$ is  differentiable for $t$ in an interval $I \subseteq \RN$ then each component function $x_i$ is differentiable. Furthermore, by the mean value theorem on $[t_{k-1},t_k] \subseteq I$ we have $x_i(t_k)-x_i(t_{k-1}) = \frac{dx_i}{dt}(t_{k_i}^*)(t_k-t_{k-1})$ for some $t_{k_i}^* \in [t_{k-1},t_k]$.  If $\zeta_j = \zeta (t_j)$ for $j=0,1,\dots , m$ then notice
\begin{equation}
 \triangle \zeta_j = \sum_{i} (x_i(t_j)-x_i(t_{j-1}) v_i = \sum_{i}\frac{dx_i}{dt}(t_j^*)(t_j-t_{j-1}) v_i 
= \triangle t_j \sum_{i}\frac{dx_i}{dt}(t_{j_i}^*) v_i 
\end{equation}
where $\triangle t_j = t_j - t_{j-1}$ and $t_{j_i}^* \in [t_{j-1},t_j]$ for $i=1,\dots , n$. Therefore,
\begin{equation}
 \sum_{j=1}^{m} f(\zeta_j) \star \triangle \zeta_j = \sum_{j=1}^{m} f(\zeta_j) \star \left( \sum_{i}\frac{dx_i}{dt}(t_{j_i}^*)v_i \right)\,  \triangle t_j. 
\end{equation}
As $m \rightarrow \infty$ we derive the following: 

\begin{thm} \label{thm:Aintegralformula}
If $\zeta: [t_o,t_1] \rightarrow \Acal$ is differentiable parametrization of a curve $C$ and $f$ is continuous near $C$ then $ \ds \int_C f ( \zeta) \star d\zeta = \int_{t_o}^{t_f} f( \zeta (t)) \star \frac{d\zeta}{dt} \, dt. $
\end{thm}

\noindent
We could have used the formula in Theorem \ref{thm:Aintegralformula} in the place of Definition \ref{defn:Aintegral}. 

\begin{exa}
For $\Acal = \RN \oplus j \RN$ with $j^2=1$ and $\zeta = x+jy$ and $f = u+jv$ then
\begin{equation}
f \,\frac{d\zeta}{dt} = (u+jv)\left( \frac{dx}{dt}+j\frac{dy}{dt}\right) = u\frac{dx}{dt}+v\frac{dy}{dt}+ j \left( v\frac{dx}{dt}+u\frac{dy}{dt}\right)
\end{equation}
Thus $ \ds \int_C f \, d\zeta = \int_C udx+vdy+ j \int_C vdy+udx.$
\end{exa}

\noindent
Formally, this integral is the natural generalization of the complex integral. We need to interconnect the real calculus of paths and the $\Acal$-calculus in a natural fashion to make further progress.

\begin{thm} \label{thm:chainruleApaths}
Let $f = \sum_i f_i v_i$ be $\Acal$-differentiable on a curve $C$ which is parametrized by $\zeta: I \rightarrow C$ where $\zeta = \sum_i x_i v_i$ for some interval $I \subseteq \RN$ (hence $x_i: I \rightarrow \RN$). Then,
$$ \frac{d}{dt} f (\zeta (t)) = f'( \zeta (t)) \star \frac{d\zeta}{dt}$$
where $\frac{d\zeta}{dt} = \sum_i \frac{dx_i}{dt} v_i$.
\end{thm}

\noindent
{\bf Proof:} given basis $\{ v_1, \dots , v_n \}$ where $v_1 = \mathds{1}$ we know the $\Acal$-CR equations read $\frac{\partial f}{\partial x_i} = \frac{\partial f}{\partial x_1} \star v_i$ for each $i = 2, \dots, n$.  Begin by applying the chain rule from real multivariate calculus,
\begin{align} 
\frac{d}{dt} f (\zeta (t)) &= \sum_i \frac{\partial f}{\partial x_i}(\zeta (t))\frac{d x_i}{dt} \\ \notag
&= \sum_i \frac{\partial f}{\partial x_1}(\zeta (t))\star v_i \frac{d x_i}{dt} \\ \notag
&= \frac{\partial f}{\partial x_1}(\zeta (t)) \star \left( \sum_i \frac{d x_i}{dt}v_i \right) \\ \notag
&= f'( \zeta(t)) \star \frac{d \zeta}{dt}. \ \ \Box
\end{align}

\noindent
There is a fundamental theorem of calculus for $\Acal$-integrals.

\begin{thm} \label{thm:FTCforalg}
Suppose $f = \frac{dF}{d\zeta}$ near a curve $C$ which begins at $P$ and ends at $Q$ then
$$ \int_C f( \zeta) \star d \zeta = F( Q) - F(P). $$
\end{thm}

\noindent
{\bf Proof:} Let $\gamma:[ t_o, t_1] \rightarrow C$ parametrize $C$. If $f = \frac{dF}{d\zeta}$ then
$ \frac{d}{dt} F( \gamma(t)) = F'(\gamma(t)) \star \frac{d\gamma}{dt} = f( \gamma(t)) \star \frac{d\gamma}{dt}$ by Theorem \ref{thm:chainruleApaths}. Thus, applying Theorems \ref{thm:Aintegralformula} and part (iii.) of \ref{thm:realAintegralprops} we derive:
\begin{align} 
 \int_C f( \zeta) \star d \zeta &= \int_{t_o}^{t_1} f( \gamma(t)) \star \frac{d\gamma}{dt} \, dt  \\ \notag
 &= \int_{t_o}^{t_1} \frac{d}{dt} \left[ F( \gamma(t)) \right] \, dt \\ \notag
 &= F(\gamma(t_1))-F(\gamma(t_o))  \\ \notag
 &= F(Q)-F(P). \ \ \Box
 \end{align}

\noindent
In other words, there is a {\it Fundamental Theorem of Calculus} II (FTC II) for $\Acal$-calculus. We conclude this Section with the analog of {\it FTC} I for $\Acal$-calculus.

\begin{thm} \label{thm:topologicalAintegral}
Let $f: U \rightarrow \Acal$ be a function where $U$ is a connected subset of $\Acal$ then the following are equivalent:
\begin{enumerate}[{\bf (i.)}]
\item $\int_{C_1} f \star d\zeta = \int_{C_2} f \star d\zeta$ for all curves $C_1,C_2$ in $U$ beginning and ending at the same points,
\item $\int_C f \star d\zeta =0$ for all loops in $U$,
\item $f$ has an antiderivative $F$ for which $\frac{dF}{d\zeta} = f$ on $U$.
\end{enumerate}
\end{thm}

\noindent
{\bf Proof:} to prove (i.) equivalent to (ii.) we need only observe $\int_{-C} f  \star d\zeta = -\int_{C} f  \star d\zeta$ and note the clear geometric connection between a loop and a pair of curves with matching end points. To prove (iii.) implies (i.) we simply assume the existence of a primitive $F$ for $f$ and make use of Theorem \ref{thm:FTCforalg}. In particular, if $f(\zeta) = F'( \zeta)$ and if $C_1$ and $C_2$ both begin at $P$ and terminate at $Q$ then 
\begin{equation}
 \int_{C_1} f  \star d\zeta = F(Q)-F(P) = \int_{C_2} f  \star d\zeta. 
\end{equation}
To prove (i.) implies (iii.), assume path-independence of the algebra integral.  Fix $\zeta_o \in U$ and define, using path independence to avoid ambiguity,
\begin{equation}
 F(\zeta) = \int_{\zeta_o}^{\zeta} f \star d\eta. 
\end{equation}
Formally, if $f = u_1v_1+ \cdots + u_nv_n$ and $d\eta = dx_1+ dx_2v_2+ \cdots + dx_nv_n$ we have
\begin{equation}
 f \star d\eta = f dx_1+ (f \star v_2)  dx_2+ \cdots + (f \star v_n) dx_n.\end{equation}
Following the usual argument, consider a path from $\zeta_o$ to $\zeta$ which is along the $v_j$ - direction near $\zeta$. Using this path we derive
\begin{equation}
 \frac{\partial F}{\partial x_j} = f \star v_j
 \end{equation}
Therefore, as we assume $v_1 = \mathds{1}$,
\begin{equation}
 \frac{\partial F}{\partial x_1} = f \star v_1 = f \ \ \Rightarrow \ \ \frac{\partial F}{\partial x_j} = f \star v_j = \frac{ \partial F}{\partial x_1} \star v_j. 
\end{equation}
Hence\footnote{we are using part (ii.) of Theorem \ref{thm:CReqnsI} to characterize the $\Acal$-CR equations}, $F$ satisfies the $\Acal$-CR equations  with $\frac{dF}{d\zeta} = f$. $\Box$  \\

\noindent
{\bf Remark:} If we added the assumption that $f$ was smooth in the real sense to the Theorem above then it follows $F$ is also real smooth and once $\Acal$-differentiable hence Theorem \ref{thm:fromonemany} provides the existence of arbitrarily many algebra derivatives of $F$. More to the point, from $\frac{dF}{d\zeta}=f$ we find $\frac{df}{d\zeta} =  \frac{d^2F}{d\zeta^2}$ hence $f$ is $\Acal$-differentiable on $U$. \\

\noindent
The exterior derivative of an $\Acal$-valued one form $\alpha = \sum_j \alpha_j v_j$ is given by component-wise exterior differentiation; $d\alpha = \sum_j d\alpha_j v_j$. Likewise, if $f = \sum_j u_jv_j$ then $df = \sum_j du_j v_j$. We say $\alpha$ is {\bf exact} if $\alpha = df$ for some function $f$ whereas $\alpha$ is {\bf closed} if $d\alpha=0$. Recall, for a simply connected subset of $\RN^n$, Poincare's Lemma states a closed form is necessarily exact. The identity $d^2=0$ holds for $\Acal$-valued forms hence every exact form is closed.

\begin{thm} \label{thm:exactAdiff}
Let $f: U \rightarrow \Acal$ be a function where $U$ is a simply connected subset of $\Acal$ and suppose $f$ is continuously differentiable in the real Frechet sense. The $\Acal$-valued one-form $f \star d \zeta$ is exact if and only if $f$ is $\Acal$-differentiable.  
\end{thm}

\noindent
{\bf Proof:} let $d \zeta = v_1dx_1+ \cdots + v_n dx_n$ and note $v_i \star v_j = v_j \star v_i$ whereas $dx_i \wedge dx_j = -dx_j \wedge dx_i$ for all $1 \leq i,j \leq n$ hence:
\begin{align} 
d\zeta \wedge d\zeta &= ( v_1dx_1+ \cdots + v_n dx_n) \wedge ( v_1dx_1+ \cdots + v_n dx_n) \\ \notag
&=\sum_{i=j} v_i \star v_j dx_i \wedge dx_j +\sum_{i<j}v_i \star v_j dx_i \wedge dx_j+\sum_{i>j}v_i \star v_j dx_i \wedge dx_j\\ \notag 
&= \sum_{i<j}(v_i \star v_j dx_i \wedge dx_j+v_j \star v_i dx_j \wedge dx_i) \\ \notag
&= 0.
\end{align} 
If $f = \sum_j u_jv_j$ is $\Acal$-differentiable then $\partial_j f = \partial_1 f \star v_j$ for $j=1,2, \dots, n$. Calculate,
\begin{align} 
 d( f \star d\zeta) &= df \wedge d\zeta \\ \notag
 &= (\partial_1f \, dx_1+\partial_2f \, dx_2+ \cdots \partial_nf \, dx_n )\ \wedge d\zeta \\ \notag
 &= (\partial_1f dx_1+\partial_1 f \star v_2dx_2+ \cdots +\partial_1f \star v_n dx_n )\ \wedge d\zeta \\ \notag
 &= \partial_1f \star d\zeta \wedge d\zeta \\ \notag
&=0.
\end{align}
Thus $f \star d\zeta$ is exact by Poincare's Lemma. 
Conversely, if $f \star d\zeta$ is exact then there exists $\phi: U \rightarrow \Acal$ for which $d\phi = f \star d\zeta$. However,
\begin{equation}
 f \star d\zeta = f \star (v_1 dx_1+ \cdots + v_n dx_n) = (f \star v_1)dx_1+ \cdots + (f \star v_n)dx_n
 \end{equation}
However, $d \phi = \partial_1 \phi \, dx_1+\partial_2 \phi \, dx_2+ \cdots + \partial_n \phi \, dx_n$ thus from $d \phi = f \star d\zeta$ we find $ f \star v_j = \partial_j \phi$ for $j=1, \dots , n$. Since $v_1 =\mathds{1}$ this gives $f = \partial_1 \phi$ and
\begin{equation}
 \partial_j f = \partial_j \partial_1 \phi = \partial_1 \partial_j \phi = \partial_1 (f \star v_j) = (\partial_1f)\star v_j. 
 \end{equation}
Therefore, $f$ is $\Acal$ differentiable on $U$ as it is a continuously differentiable function which satisfies the $\Acal$-CR equations on $U$. $\Box$  \\

\noindent
Finally, we arrive at the analog of Cauchy's Integral Theorem for $\Acal$:

\begin{coro} \label{coro:adiffloopszero}
If $U \subseteq \Acal$ is simply connected then $\int_C f \star d \zeta = 0$ for all loops $C$ in $U$ if and only if $f$ is $\Acal$-differentiable on $U$.
\end{coro}

\noindent
{\bf Proof:} let $U$ be simply connected. Observe $\int_C f \star d \zeta = 0$ for all loops $C$ in $U$ is equivalent to the existence of $F$ on $U$ for which $\frac{dF}{d\zeta} = f$ by Theorem \ref{thm:topologicalAintegral}. However,
\begin{align} 
 dF &= \partial_1 F \, dx_1+\partial_2 F \, dx_2 \cdots + \partial_nF \, dx_n \\ \notag
 &= \partial_1 F \, dx_1+\partial_1 F \star v_2 \, dx_2+ \cdots  + \partial_1 F \star v_n \, dx_n  \\ \notag
 &= \partial_1 F \star \left( dx_1+ v_2dx_2+ \cdots + v_ndx_n \right) \\ \notag
 &= f \star d \zeta
 \end{align}
whence we see $f \star d \zeta$ is exact on $U$. Then the Corollary follows by Theorem \ref{thm:exactAdiff}. $\Box$  \\

\noindent
In the case of complex analysis, Goursat's Theorem allows us to show that a function which is complex differentiable is necessarily continuously-complex-differentiable. In other words, Goursat's Theorem allows proof that the map $z \rightarrow f'(z)$ is continuous provided $f'(z)$ exists over some domain. The proof of Goursat's Theorem transfers nicely to the context of $\Acal$-differentiable functions. However, we omit the proof in this paper.

\begin{thm} \label{thm:FTCIforalg}
Let $C$ be a differentiable curve from $\zeta_o$ to $\zeta$ in $U \subseteq \Acal$ where $U$ is an open simply connected subset of $\Acal$. Assume $f$ is $\Acal$ differentiable on $U$ then 
$$ \frac{d}{d\zeta}\int_C f( \eta) \star d \eta = f(\zeta). $$
\end{thm}

\noindent
{\bf Proof:} If $U$ is simply connected and $f$ is $\Acal$-differentiable on $U$ then Corollary \ref{coro:adiffloopszero} provides $\oint_L f \star d\eta = 0$ for any loop in $U$. Thus, using the equivalence of (ii.) and (iii.) in Theorem \ref{thm:topologicalAintegral} we find there exists an $\Acal$-differentiable function $F$ on $U$ for which $f = \frac{dF}{d\zeta}$. Let $C$ be a differentiable curve from $\zeta_o$ to $\zeta$ and note by Theorem \ref{thm:FTCforalg} with $P=\zeta_o$ and $Q= \zeta$
\begin{equation}
 \int_C f( \eta) \star d \eta = F(\zeta) - F( \zeta_o). 
 \end{equation}
Therefore, as $\zeta_o$ is a constant and $f = \frac{dF}{d\zeta}$ we derive $  
 \frac{d}{d\zeta}\int_C f( \eta) \star d \eta = f( \zeta)$. $\Box$ \\
 
\begin{exa} \label{exa:logdefined}
Suppose $\Acal$ is a commutative, unital algebra and let $U$ be an open simply connected subset of $\Acalx$ which contains $\mathds{1}$. We can meaningfully write $ \ds g(\zeta) = \int_{\mathds{1}}^{ \zeta} \frac{d\eta}{\eta}$ and by Theorem \ref{thm:FTCIforalg} we have $\frac{dg}{d\zeta} = \frac{1}{\zeta}$ on $U$. This discussion is continued in \cite{bedellII} where it is shown $g$ defines a natural inverse function to the exponential on $\Acal$. We also recommend \cite{cookfreese} for details on how the exponential can be defined over many algebras. 
\end{exa}
 
 \section{Conclusions and future work}
 This paper primarily contains results which are analogs to the usual topics covered in first and second semester calculus. We found how to differentiate once or multiple times and we saw how integration is defined for functions over an algebra. Most of our results are quite general though, we have not discussed the exponential function, sine, cosine or their hyperbolic compatriots. \\
 
\noindent
 In the sequel to this paper with Daniel Freese \cite{cookfreese} we study the theory of power series over an algebra. Series give a natural method to generalize the usual elementary functions to an algebra. Then, we continue past calculus II to study differential equations over $\Acal$ with Nathan BeDell in \cite{cookbedell} where we prove the essential theory to frame linear $\Acal$-ODEs and we exhibit novel solution techniques which solve {\bf any} nondegenerate $\Acal$-ODE. \\
 
 \noindent
I should also mention, Nathan BeDell provided three independent papers to support the study of $\Acal$-Calculus. In particular, in \cite{bedellI} he studies zero-divisors and algebraic preliminaries,in\cite{bedellII} he develops logarithms over many algebras and in \cite{bedellIII} he attempts to find an analog of polar form for an algebra. Moreover, \cite{bedellIII} also continues the study of the $N$-Pythagorean Theorem found in \cite{cookfreese}. \\
  
\noindent
Another direction would be to further explore problems such as Examples \ref{exa:waveqn} and \ref{exa:AODEbackwards}. While we sketched a program for converting a real PDE into an $\Acal$-ODE we did not attempt to answer the most interesting and challenging question: {\it which $\Acal$ should guide the conversion ?}. It might be interesting to understand how the algebra substitution ties into the symmetries of a given PDE. \\

\noindent
We also do not understand the complete connection between the $\Acal$-Laplacian introduced by W.S. Leslie in \cite{cook} and the generalized $\Acal$-Laplace equations of Wagner. It might be possible to prove something interesting in terms of the Wirtinger calculus. \\

\noindent
 Our work here and in the sequels to this paper are by no means complete. Many theorems of complex analysis must have some analog in the $\Acal$-calculus. On the other hand, the exposition of calculus for just one choice of $\Acal$ has filled entire volumes already. The directions for future work here are endless.

\section{Acknowledgements}
The author is thankful to R.O. Fulp for helpful comments on a rough draft of this article. The author is also thankful for many conversations with M. L. Nguyen, W.S. Leslie, B. Zhang, D. Freese and N. BeDell which laid the foundation for this work.

\end{document}